\documentclass[12pt,reqno]{amsart}
\usepackage{amslfmt,amsthm}
\usepackage{graphicx}
\usepackage[svgnames]{xcolor}
\usepackage{hyperref}

\usepackage{lineno}

\begin{document}

\title[Vanishing adsorption limit of Riemann problem solutions for the polymer model]%
{Vanishing adsorption limit of Riemann problem solutions for the polymer model}

\author[Y. Petrova]{Yulia Petrova}%
\address{%
Pontifical Catholic University of Rio de Janeiro --- PUC-Rio, Rio de Janeiro, Brazil}%
\email{yu.pe.petrova@gmail.com}

\author[B. Plohr]{Bradley J. Plohr}%
\address{%
Los Alamos, New Mexico, USA}%
\email{bradley.j.plohr@gmail.com}

\author[D. Marchesin]{Dan Marchesin}%
\address{%
Instituto Nacional de Matem\'atica Pura e Aplicada, Rio de Janeiro, Brazil}%
\email{marchesi@impa.br}

\begin{abstract}
We examine the vanishing adsorption limit
of solutions of Riemann problems for
the Glimm-Isaacson model of chemical flooding of a petroleum reservoir.
A contact discontinuity is deemed admissible if it is the limit of
traveling waves or rarefaction waves for an augmented system
that accounts for weak chemical adsorption onto the rock.
We prove that this criterion justifies the admissibility criteria adopted
previously by Keyfitz-Kranzer, Isaacson-Temple, and de~Souza-Marchesin,
provided that the fractional flow function
depends monotonically on chemical concentration.
We also demonstrate that the adsorption criterion
selects the undercompressive contact discontinuities
required to solve the general Riemann problem
in an example model with non-monotone dependence.
\end{abstract}




\maketitle

\section{Introduction}
\label{sec:Introduction}
\begin{figure}[ht]
\centering \includegraphics[width=0.32\textwidth]
{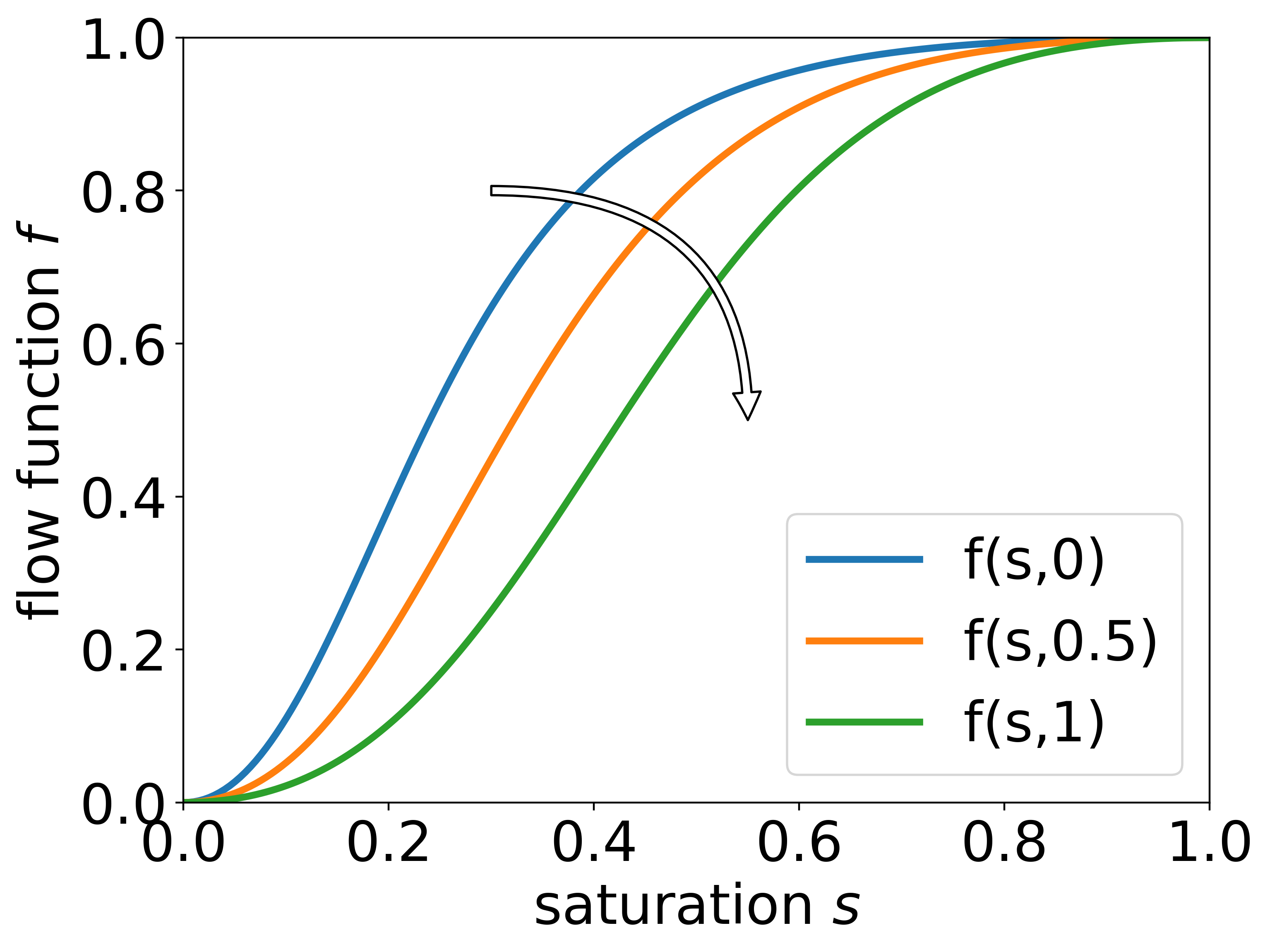}
\includegraphics[width=0.32\textwidth]{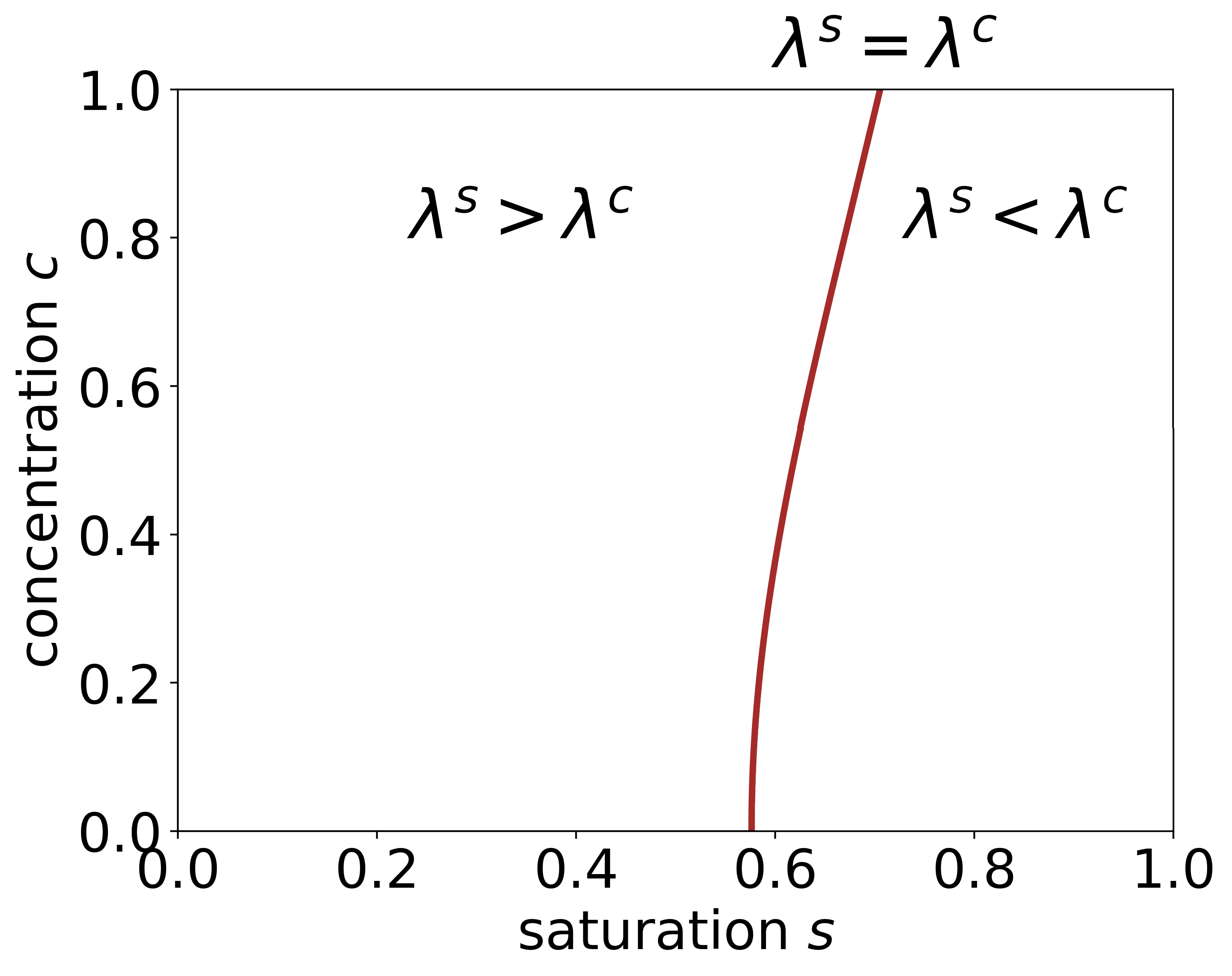}
\includegraphics[width=0.32\textwidth]{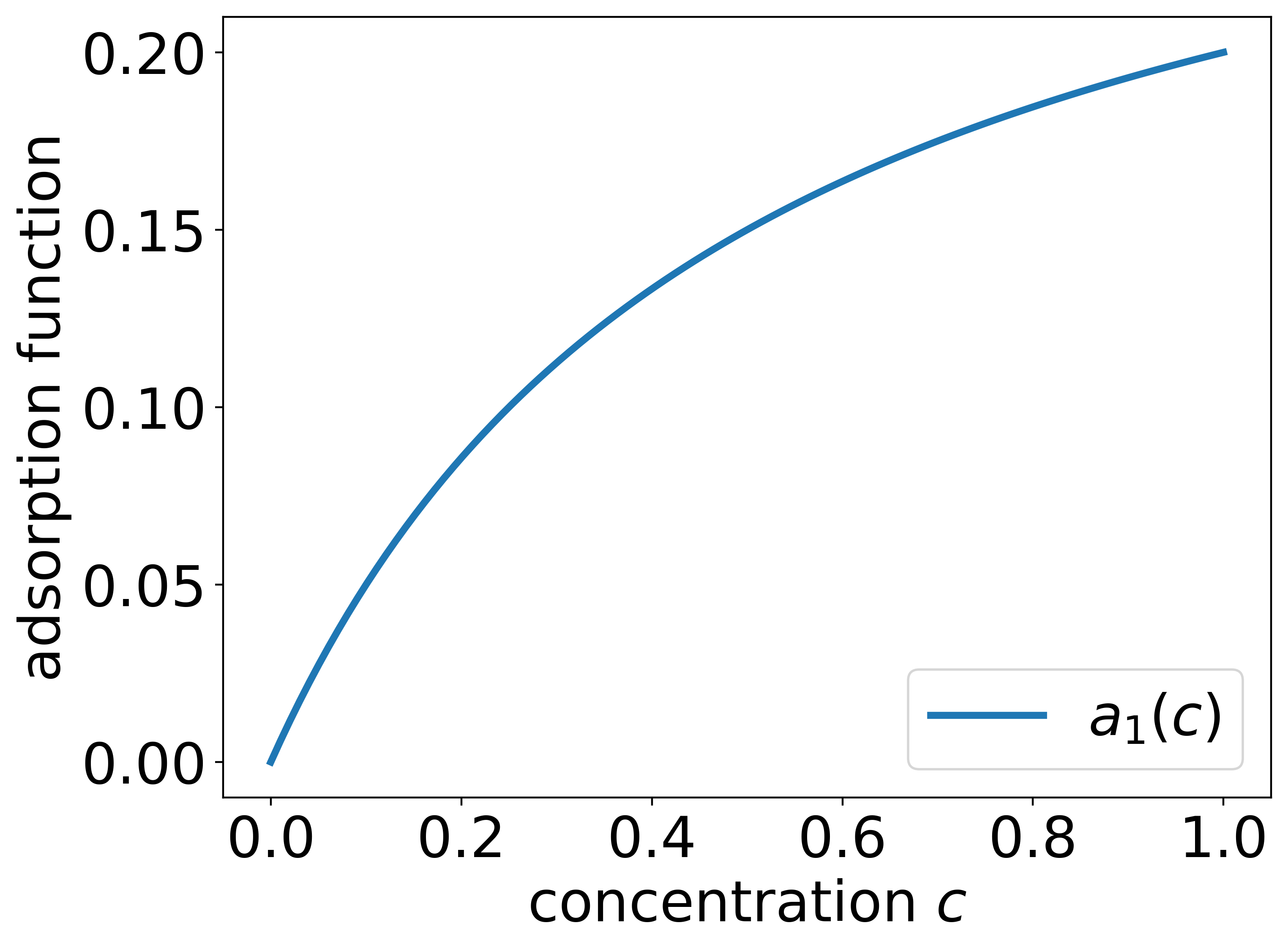} \\
\qquad\qquad
(a)\qquad\qquad\qquad\qquad\qquad\qquad
(b)\qquad\qquad\qquad\qquad\qquad\qquad
(c)\qquad\qquad
\caption{(a)~The fractional flow function
$f(s,c)$ for the monotone polymer model,
plotted \emph{vs}.~$s$ for some fixed values of $c$.
(b)~The $(s,c)$-space,
divided by the interior coincidence curve $\mathcal{C}$
into two regions of strict hyperbolicity,
$\{\lambda^s>\lambda^c\}$ and
$\{\lambda^s<\lambda^c\}$.
(c)~The adsorption function $a_1(c)$,
which appears in Eq.~\eqref{eq:conslaw_adsorption} below.}
\label{fig:f_a_t}
\end{figure}

In this paper,
we propose a new, physically-motivated admissibility criterion
for contact discontinuities
and analyze its relationships with existing criteria.
The context for developing this admissibility criterion is
the flow through porous rock
of oil, water, and a chemical agent within the water,
as modeled by the pair of conservation laws
\begin{equation}
\label{eq:conslaw}
\begin{split}
s_t + f(s, c)_x &= 0,\\
(c\,s)_t + \left[c\,f(s, c)\right]_x &= 0.
\end{split}
\end{equation}
Here $s=s(x,t) \in [0,1]$ is the saturation of the water phase
within the oil-water fluid mixture,
$c=c(x,t) \in [0,1]$ is the concentration
of the chemical agent in the water phase,
and $f=f(s,c)$ is the fractional flow of water
(the Buckley-Leverett flux function).
Following the classical paper of Buckley-Leverett~\cite{BL1942},
we assume that $f\in C^2([0,1]^2)$ is,
for every $c\in[0,1]$,
a monotone and S-shaped function of $s$
such that $f(0,c)=0$, $f(1,c)=1$, and $f_s(0,c)=f_s(1,c)=0$,
as illustrated in Fig.~\ref{fig:f_a_t}(a).

Our main result,
the equivalence of admissibility criteria
(Theorem~\ref{thm:main-thm}),
concerns models that satisfy the following \emph{monotonicity condition}:
\begin{align}
\label{eq:f-cond-monotone}
\text{$f_c<0$ for $s\in(0,1)$ and $c\in[0,1]$.}
\end{align}
Such monotone behavior is typical for chemical flooding
using a polymer as the chemical agent,
because the viscosity of water usually increases when polymer is added.
However,
in Sec.~\ref{sec:example},
we consider a model that does not obey the monotonicity condition.

We focus on finding solutions of the general Riemann initial-value problem:
\begin{align}
\label{eq:Riemann-problem}
\left.U\right\rvert_{t=0}=
\begin{cases}
U_L& \text{for $x<0$}, \\
U_R& \text{for $x>0$}.
\end{cases}
\end{align}
(Here $U$ denotes the state $(s, c)$.)
The solutions we construct are scale-invariant,
in that they depend on $(x,t)$ only through the combination $\xi = x/t$.
Along with smooth scale-invariant waves (centered rarefaction waves),
a solution can contain centered discontinuities,
which have the form
\begin{equation}
\label{eq:discontinuity}
U(x, t) = \begin{cases} U_- & \text{if $x < \sigma\,t$,} \\
U_+ & \text{if $x > \sigma\,t$.}
\end{cases}
\end{equation}
A discontinuity is required to satisfy the conservation laws
in the sense of distributions,
which amounts to satisfying the Rankine-Hugoniot condition
that constrains the left state $U_-$, right state $U_+$,
and propagation speed $\sigma$.

System~\eqref{eq:conslaw} was first considered
by Keyfitz and Kranzer~\cite{keyfitz-kranzer1980}
and by Isaacson~\cite{eli1980}.
It is known by various names,
including the Glimm-Isaacson model,
the KKIT model (referring to Keyfitz, Kranzer, Isaacson, and Temple),
and the chemical flooding model.
It is closely related to several models of multi-phase flow in porous media,
including chromatography~\cite{rheearis1986},
gas flooding~\cite{orr2007}, CO$_2$-flooding~\cite{wahanik2011},
expansive flow~\cite{lam2009}, thermal flow~\cite{damota1992},
chemical flooding~\cite{souza1995}, 
vertical flow~\cite{rodriguez2013}, 
and other examples described in books on enhanced oil
recovery~\cite{bedrikovetsky1993,bruining2022,lake2014}.
We shall call it the polymer model.

The system~\eqref{eq:conslaw} has the following two characteristic speeds
(see, \emph{e.g.},\ Temple~\cite{temple1982}):
\begin{itemize}
\item The saturation (or $s$-type)
characteristic speed is $\lambda^s=f_s$,
which has associated right eigenvector $r^s = (1, 0)^T$;
this characteristic family is genuinely nonlinear except where $f_{ss} = 0$.
For an $s$-type solution of system~\eqref{eq:conslaw},
$c$ has a fixed value and $s$ is governed by the first equation
in system~\eqref{eq:conslaw},
which is the Buckley-Leverett scalar conservation law.
Such a solution can be an $s$-rarefaction wave,
an $s$-shock wave, or a contiguous group of such waves,
which we call an $s$-wave group.

\item The concentration (or $c$-type)
characteristic speed is $\lambda^c=f/s$; this characteristic family is linearly degenerate
(meaning that $(D\lambda^c)\,r^c \equiv 0$,
where $r^c$ is a corresponding right eigenvector field).
As a consequence,
its rarefaction curves and Hugoniot loci coincide
and every $c$-type solution is a contact discontinuity,
\emph{i.e.},\  a weak solution~\eqref{eq:discontinuity} such that
\begin{equation}
\label{eq:double-sonic}
\lambda^c(U_-) = \sigma = \lambda^c(U_+),
\end{equation}
which we call a $c$-wave.
\end{itemize}

A challenging feature of system~\eqref{eq:conslaw}
is its loss of strict hyperbolicity in the interior of state space:
the coincidence locus,
where the characteristic speeds coincide,
comprises not only the boundary line $s = 0$ but also the set
\begin{equation}
\label{eq:coincidence}
\mathcal{C}:=\{\,U\,:\,\text{$\lambda^s(U)=\lambda^c(U)$ and $s \ne 0$}\,\},
\end{equation}
which we call the interior coincidence locus.
This set is a smooth curve along which $s$ is parametrized by $c$,
as illustrated in Fig.~\ref{fig:f_a_t}(b).
To the left and right of the interior coincidence locus
are the separate regions of strict hyperbolicity
\begin{align*}
\{\,\lambda^s>\lambda^c\,\}
&:=\{U\,:\,\lambda^s(U)>\lambda^c(U)\,\}, \\
\{\,\lambda^s<\lambda^c\,\}
&:=\{U\,:\,\lambda^s(U)<\lambda^c(U)\,\}.
 \end{align*}

Relative to more complete models of fluid flow,
the conservation laws~\eqref{eq:conslaw}
have been simplified by neglecting some physical effects.
To account for these effects,
the conservation laws must be supplemented by an admissibility criterion
that selects the discontinuities that appear in a Riemann solution,
just as the conservation laws governing gas dynamics
is amended by an entropy condition~\cite{lax57}.

Each $s$-wave group solves the Buckley-Leverett scalar conservation law
for a fixed value of $c$,
so we require it to satisfy the standard
Ole\u{\i}nik admissibility criterion~\cite{Ole57,Ole63}. For a system of equations, a general and powerful admissibility criterion derives from adding a diffusive term to the equations and requiring a shock wave to be the limit of traveling wave solutions as the diffusion coefficient vanishes~\cite{CouFri48,Gel59,Gel63}. 

\begin{remark}
For a system of equations, there can be several diffusion coefficients.
An example is the compressible Navier-Stokes system,
which involves fluid viscosity and heat conductivity coefficients.
In some circumstances,
the vanishing viscosity limit depends on the limiting values
of the ratios of the diffusion coefficients.
For instance, in combustion theory,
weak deflagration waves depend on the Prandtl number~\cite[p. 234]{CouFri48}.
We encounter this phenomenon in the polymer model; see Sec.~\ref{sec:example}.
Accordingly, we do not require the vanishing viscosity limit to be independent
of the ratios of diffusion coefficients;
rather, we require the specification of a vanishing viscosity criterion 
to include these ratios.
\end{remark}

However, a $c$-wave is a contact discontinuity,
for which there is no standard admissibility criterion.
In gas dynamics,
the archetype of conservation laws,
every contact discontinuity is admissible,
but blanket admissibility of contact discontinuities
does not carry over to all systems.
The vanishing viscosity admissibility criterion does not apply to contact discontinuities,
which lack nonlinear wave focusing to balance viscous spreading,
so that traveling waves are not supported.

Keyfitz and Kranzer~\cite{keyfitz-kranzer1980}
and Isaacson~\cite{eli1980} recognized that,
for system~\eqref{eq:conslaw},
treating all contact discontinuities as admissible
leads to nonuniqueness of solutions of Riemann problems.
Instead,
they both adopted the following admissibility criterion
for contact discontinuities,
calling it the generalized Lax entropy condition,
and deduced the existence and uniqueness of solutions
of the general Riemann problem~\eqref{eq:Riemann-problem}
for system~\eqref{eq:conslaw}.

\begin{definition}[Criterion of Keyfitz and Kranzer~\cite{keyfitz-kranzer1980}
and of Isaacson~\cite{eli1980}]
\label{def:lax}
A $c$-wave,
\emph{i.e.},\  a contact discontinuity,
is admissible if one of the following conditions holds:
\begin{align}
\label{eq:1-contact}
\text{either } &\text{$\sigma<\lambda^s(U_-)$ and $\sigma<\lambda^s(U_+)$}; \\
\label{eq:2-contact}
\text{or } &\text{$\sigma>\lambda^s(U_-)$ and $\sigma>\lambda^s(U_+)$}.
\end{align}
In addition,
a $c$-wave is deemed admissible if it is the limit of
a sequence of $c$-waves satisfying~\eqref{eq:1-contact}
or a sequence satisfying~\eqref{eq:2-contact}.
\end{definition}

Isaacson and Temple~\cite{isaacsontemple1986}
used an equivalent admissibility criterion.

\begin{definition}[Criterion of Isaacson and Temple~\cite{isaacsontemple1986}]
\label{def:IT}
A $c$-wave is admissible if either
(1)~both $U_-$ and $U_+$ belong to $\{\,\lambda^s\ge\lambda^c\,\}$
or (2)~both $U_-$ and $U_+$ belong to $\{\,\lambda^s\le\lambda^c\,\}$.
\end{definition}

A generalization of these criteria was proposed by
de~Souza and Marchesin~\cite{cido-dan1998},
who named it the KKIT criterion.

\begin{definition}[Criterion of de Souza and Marchesin~\cite{cido-dan1998}]
\label{def:KKIT}
A $c$-wave is admissible if $c$ varies continuously and monotonically
along a sequence of contact curves leading from $U_-$ to $U_+$.
\end{definition}

In this paper,
we show that the foregoing admissibility criteria
relate to a physical phenomenon,
namely adsorption of the chemical agent onto the surface of the porous rock,
in the limit when the adsorption coefficient tends to zero.
In particular,
we consider the following pair of conservation laws
(first analyzed by Johansen and Winther~\cite{johansen1988}):
\begin{equation}
\label{eq:conslaw_adsorption}
\begin{split}
s_t + f(s, c)_x &= 0, \\
\left[c\,s+\alpha\,a_1(c)\right]_t + \left[c\,f(s, c)\right]_x &= 0,
\end{split}
\end{equation}
supplemented by a vanishing viscosity admissibility criterion.
This system is distinguished from system~\eqref{eq:conslaw}
by the term $\alpha\,a_1(c)$,
which is the equilibrium amount of chemical agent adsorbed onto the rock surface
when the concentration is $c$.
Here $\alpha > 0$ and the function $a_1(c)\in C^2([0,1])$ is assumed to be
strictly increasing ($a_1'>0$) and strictly concave ($a_1''<0$),
as depicted in Fig.~\ref{fig:f_a_t}(c).
We will refer to system~\eqref{eq:conslaw_adsorption} as model $M_\alpha$.
In particular, system~\eqref{eq:conslaw} is model~$M_0$.

Although adsorption of a chemical agent onto a porous rock is commonplace,
the adsorption term $\alpha\,a_1(c)$ is usually small compared to $c\,s$.
Often, the effect of adsorption is expected to be negligible
and $\alpha$ is set to zero.
However, as we shall see,
even small, but nonzero, adsorption
serves to distinguish the admissible contact discontinuities.
Thus motivated, we propose the following admissibility criterion.

\begin{definition}[Vanishing adsorption admissibility criterion
for contact discontinuities]
\label{def:AD}
A contact discontinuity for system~\eqref{eq:conslaw}
is admissible provided that it is the $L^1_\text{loc}$ limit,
as $\alpha \to 0^+$,
of a sequence of admissible solutions
of the Riemann problem~\eqref{eq:Riemann-problem}
for system~\eqref{eq:conslaw_adsorption}.
\end{definition}

\begin{remark}
We emphasize that a contact discontinuity is admissible
according to Def.~\ref{def:AD}
not only if it is a limit of shock waves
but also if it is a limit of rarefaction waves, composite waves,
or more general Riemann solutions for system~\eqref{eq:conslaw_adsorption}.
(In contrast, for the vanishing viscosity admissibility criterion,
a shock wave must be the limit of traveling wave solutions.)
\end{remark}

Consider a polymer model~\eqref{eq:conslaw} that satisfies the monotonicity
condition~\eqref{eq:f-cond-monotone}.
As observed in Refs.~\cite{isaacsontemple1986,cido-dan1998},
the set of admissible Riemann solutions
is the same under each of the following criteria for
contact discontinuities:
\begin{itemize}
\item Keyfitz-Kranzer/Isaacson criterion of Def.~\ref{def:lax};
\item Isaacson-Temple criterion of Def.~\ref{def:IT};
\item de~Souza-Marchesin criterion of Def.~\ref{def:KKIT}.
\end{itemize}
\noindent
The main result of this paper is the following
equivalence of admissibility criteria:
\begin{theorem}
\label{thm:main-thm}
For a polymer model~\eqref{eq:conslaw} that satisfies the monotonicity
condition~\eqref{eq:f-cond-monotone},
the set of admissible Riemann solutions
under the vanishing adsorption criterion of Def.~\ref{def:AD}
is the same as under the criteria of Defs.~\ref{def:lax}-\ref{def:KKIT}.
\end{theorem}

As a consequence of the equivalence of admissibility criteria,
the proof of existence and uniqueness~\cite{eli1980,keyfitz-kranzer1980}
yields the following corollary:

\begin{corollary}
For a model~\eqref{eq:conslaw}
that conforms to the monotonicity condition~\eqref{eq:f-cond-monotone},
if the vanishing adsorption admissibility criterion
for contact discontinuities is adopted,
then the solution of the general Riemann problem~\eqref{eq:Riemann-problem}
exists and is unique.
\end{corollary}

\begin{remark}
We emphasize that Theorem~\ref{thm:main-thm}
applies to the specific model~\eqref{eq:conslaw}
with monotone variation of $f$ in $c$.
We will see that this equivalence of criteria
does not hold for non-monotone models.
\end{remark}

Contact discontinuities obeying the Keyfitz-Kranzer/Isaacson
inequalities~\eqref{eq:1-contact} and~\eqref{eq:2-contact}
are not the only ones that occur in Riemann solutions.
In Ref.~\cite{cido-dan1998},
de~Souza and Marchesin made use of contact discontinuities that satisfy
\begin{equation}
\label{eq:crossing}
\text{$\lambda^s(U_-) < \sigma < \lambda^s(U_+)$}
\end{equation}
along with their admissibility criterion of Def.~\ref{def:KKIT}.
They showed that these discontinuities
(which they called transitional contact discontinuities)
are needed to solve certain Riemann problems.

These discontinuities are analogues,
for the case of a linearly degenerate characteristic family,
of non-classical shock waves that occur in various systems of conservation laws.
Adopting the terminology used for such waves to the present context,
we say that a contact discontinuity has
a \emph{crossing} configuration of characteristic paths~\cite{IsaMarPlo88a}
if inequalities~\eqref{eq:crossing} hold,
and we will call it \emph{undercompressive} if,
in addition, it is admissible in an appropriate sense. 
\begin{remark}
The term ``undercompressive'' was introduced for shock waves in Ref.~\cite{SheSchMar87}. An undercompressive shock wave has also  been called a \emph{transitional} shock wave because it is a prime example of the more general concept of \emph{transitional wave group}~\cite{transitional1990}.
\end{remark}

\begin{example}
\label{ex:undercompressive-contact}
An example~\cite{nonmonotonicity2021} of a polymer model~\eqref{eq:conslaw},
presented in Sec.~\ref{sec:example},
has a contact discontinuity with a crossing configuration of characteristics
that satisfies the de~Souza-Marchesin criterion
but neither the Keyfitz-Kranzer/Isaacson nor the Isaacson-Temple criteria.
(Necessarily,
the monotonicity condition~\eqref{eq:f-cond-monotone}
does not hold for this model.)
Therefore, the set of admissible Riemann solutions
is different for these different criteria.
\end{example}

In Sec.~\ref{sec:example},
we demonstrate the adsorption admissibility
of undercompressive contacts for a particular non-monotone polymer model:

\begin{theorem}
\label{thm:undercompressive-admissible}
The undercompressive contact discontinuities
found by de~Souza and Marchesin~\cite{cido-dan1998}
satisfy the vanishing adsorption admissibility criterion.
\end{theorem}

The paper is organised as follows.
In Sec.~\ref{sec:review},
we present some background on conservation laws
and review properties of the chemical flooding models~\eqref{eq:conslaw}
and~\eqref{eq:conslaw_adsorption}.
In Sec.~\ref{sec:proof},
we prove the equivalence of admissibility criteria
(Theorem~\ref{thm:main-thm}).
Finally, in Sec.~\ref{sec:example},
we describe Example~\ref{ex:undercompressive-contact}
and prove that undercompressive contact discontinuities
obey the adsorption admissibility criterion
(Theorem~\ref{thm:undercompressive-admissible}).

\section{Background on conservation laws and chemical flooding models}
\label{sec:review}
In this section,
we recall important properties of hyperbolic conservation laws in general
and of the polymer models~$M_\alpha$ with $\alpha\geq0$ in particular.
More details can be found in textbooks,
such as Refs.~\cite{dafermos,smoller2012},
in previous work~\cite{cido-dan1998,eli1980,eliblake1986,keyfitz-kranzer1980}
on the polymer model without adsorption,
and in Ref.~\cite{johansen1988} for the polymer model with adsorption.

\subsection{Basic notions for hyperbolic conservation laws}
Systems~\eqref{eq:conslaw} and~\eqref{eq:conslaw_adsorption}
take the form of a conservation law
\begin{equation}
G(U)_t + F(U)_x = 0
\label{eq:conslaw-general}
\end{equation}
governing the evolution of the state $U\in\Omega\subseteq \mathbb{R}^2$.
Here $G(U)$ is a vector of conserved quantities
and $F(U)$ is a vector of flux functions.
We intend to solve the Riemann initial-value
problem~\eqref{eq:Riemann-problem}.
Because the governing equations and initial data are scale-invariant,
we seek scale-invariant solutions,
meaning that $U$ depends on $(x, t)$ solely through $\xi := x/t$.

\subsubsection{Rarefaction waves}
In a space-time region where $U$ is continuously differentiable,
system~(\ref{eq:conslaw-general}) requires that
\begin{equation}
DG(U)\,U_t + DF(U)\,U_x = 0.
\label{eq:differentiable_solution}
\end{equation}
When $U$ is scale invariant,
this equation reduces to the ordinary differential equation (ODE)
\begin{equation}
\left[-\xi\,DG(U) + DF(U)\right] U_\xi = 0.
\label{eq:differentiable_scale-invariant}
\end{equation}
We distinguish two kinds of solutions of this equation.
\begin{itemize}
\item A constant state:
in an open interval of $\xi$,
$U \equiv U^{\text{const}}$ is constant,
so that $U_\xi \equiv 0$.
\item A rarefaction wave:
in an open interval of $\xi$,
$U_\xi$ is a right eigenvector of the characteristic matrix
\begin{equation}
\label{eq:matrix-A}
A(U) := [DG(U)]^{-1}\,DF(U),
\end{equation}
and $\xi$ equals the corresponding eigenvalue $\lambda(U)$,
which is called a characteristic speed.
\end{itemize}

State space can be decomposed into three parts.
The strictly hyperbolic region
is the set of states $U$ such that $A(U)$ has two distinct real eigenvalues.
In a sufficiently small neighborhood of a state in this region,
each eigenvalue is a smooth function $\lambda(U)$,
and a corresponding smooth right eigenvector field $r(U)$ can be constructed.
In contrast, in the elliptic region,
the two eigenvalues of $A(U)$ are not real,
so that the null space of the matrix in brackets
in Eq.~\eqref{eq:differentiable_scale-invariant} is empty.
Finally,
on the coincidence locus,
the matrix $A(U)$ has a repeated eigenvalue,
\emph{i.e.},\ its discriminant is zero.

Starting from state $U_-$ in the strictly hyperbolic region,
a rarefaction wave can be constructed
for the characteristic speed $\lambda(U)$ by
(i)~solving the initial-value problem $U_\eta = r(U)$ with $U(0) = U_-$,
(ii)~inverting the relationship $\lambda(U(\eta)) = \xi$, and
(iii)~replacing the parameter $\eta$ by $\xi$
in $U(\eta)$ to obtain $U^{\text{raref}}(\xi)$.
Step~(ii) requires $U(\eta)$ to lie entirely within a connected component
of the region of state space where $(D\lambda)\,r \ne 0$.
Therefore, the subset of state space where $(D\lambda)\,r = 0$
plays an important role when constructing rarefaction waves
associated to $\lambda(U)$.
Two situations are typical:
either this subset is a curve,
called the inflection locus for $\lambda(U)$;
or $(D\lambda)\,r \equiv 0$ throughout an open region
in which the characteristic family  $\lambda(U)$ is said to be linearly degenerate.

\subsubsection{Shock waves}
A shock wave is a scale-invariant discontinuous solution
of system~(\ref{eq:conslaw-general}).
Such a solution, which takes the form~\eqref{eq:discontinuity},
separates left and right states $U_-$ and $U_+$
and propagates with speed $\sigma$.
The weak form of the conservation law~\eqref{eq:conslaw-general}
requires the Rankine-Hugoniot condition
\begin{equation}
- \sigma\,\left[G(U_+) - G(U_-)\right] + F(U_+) - F(U_-) = 0.
\label{eq:Rankine-Hugoniot}
\end{equation}

The Hugoniot locus $\mathcal{H}(U_-)$ for reference state $U_-$
consists of all states $U_+$ such that
there exists a speed $\sigma$ for which
the Rankine-Hugoniot condition~(\ref{eq:Rankine-Hugoniot}) is satisfied.
Usually,
the Hugoniot locus for $U_-$ is a smooth curve except at $U_-$,
where it crosses itself.

\subsection{Polymer models, with and without adsorption}
\label{subsec:zero-adsorption}

\subsubsection{Characteristic speeds}
For the system~\eqref{eq:conslaw_adsorption},
the characteristic matrix~\eqref{eq:matrix-A} is
\begin{equation}
A = \begin{pmatrix}
f_s & f_c \\
0 & f/(s+\alpha\,a_1')
\label{eq:characteristic-matrix}
\end{pmatrix}.
\end{equation}
The eigenvalues of $A$,
\emph{i.e.},\ the characteristic speeds for the system,
are
\begin{equation}
\label{eq:lambda}
\lambda^s = f_s
\qquad \text{and} \qquad
\lambda^c_\alpha := f/(s+\alpha\,a_1').
\end{equation}
We choose the right eigenvectors
corresponding to $\lambda^s$ and $\lambda^c_\alpha$ to be
\begin{equation}
\label{eq:eigenvectors}
r^s := \begin{pmatrix} 1 \\ 0 \end{pmatrix}
\qquad \text{and} \qquad
r^c_\alpha := \begin{pmatrix} -f_c \\
\lambda^s - \lambda^c_\alpha \end{pmatrix}.
\end{equation}
As both characteristic speeds are real, this model is hyperbolic.
However,
the characteristic speeds coincide,
not only on the boundary line $s = 0$,
but also
\begin{align}
\label{eq:coincidence-alpha}
\mathcal{C}_{\alpha}
:= \{\, U \,:\, \text{$\lambda^s(U) = \lambda^c_\alpha(U)$ and $s \ne 0$} \,\},
\end{align}
which generalizes the interior coincidence locus~\eqref{eq:coincidence}.
As is easily seen by graphing $f$ \emph{vs.}~$s$ with $c$ fixed,
$\mathcal{C}_{\alpha}$ is a smooth curve
along which $s$ is parametrized by $c$,
which we write as $s = s_\text{coinc}(c)$.

\subsubsection{Rarefaction waves}
For a rarefaction wave associated to the $s$-type characteristic speed,
the formula for $r^s$ entails that $c$ remains constant.
Thus, it is a solution of the Buckley-Leverett scalar conservation law
(the first of Eqs.~\eqref{eq:conslaw}) for fixed~$c$.
Such a solution cannot cross the $s$-type inflection locus,
which is the curve in the $(s,c)$-plane where $f_{ss} = 0$.

The $c$-type characteristic family is genuinely nonlinear except where
the following quantity vanishes:
\begin{align}
\label{eq:dlambda_c-ads}
\left(D\lambda^c_\alpha\right) r^c_\alpha
=-\alpha\,\frac{fa_1''(\lambda^s-\lambda^c_\alpha)}{(s+\alpha\,a_1')^2}.
\end{align}
When $\alpha = 0$,
it vanishes identically,
\emph{i.e.},\ the $c$-type family is linearly degenerate:
$\lambda^c_0(U)$ remains constant along each integral curve
of the ODE $\dot U = r^c_0(U)$.
We denote the integral curve containing $U_-$,
which is the level curve $\lambda^c_0(U) = \lambda^c_0(U_-)$,
by $\mathcal{I}_0(U_-)$;
such curves foliate state space.
These curves cannot be used to construct $c$-type rarefaction waves;
rather, they are loci of contact discontinuities.

\begin{figure}[ht]
    \centering
        \includegraphics[width=0.45\textwidth]{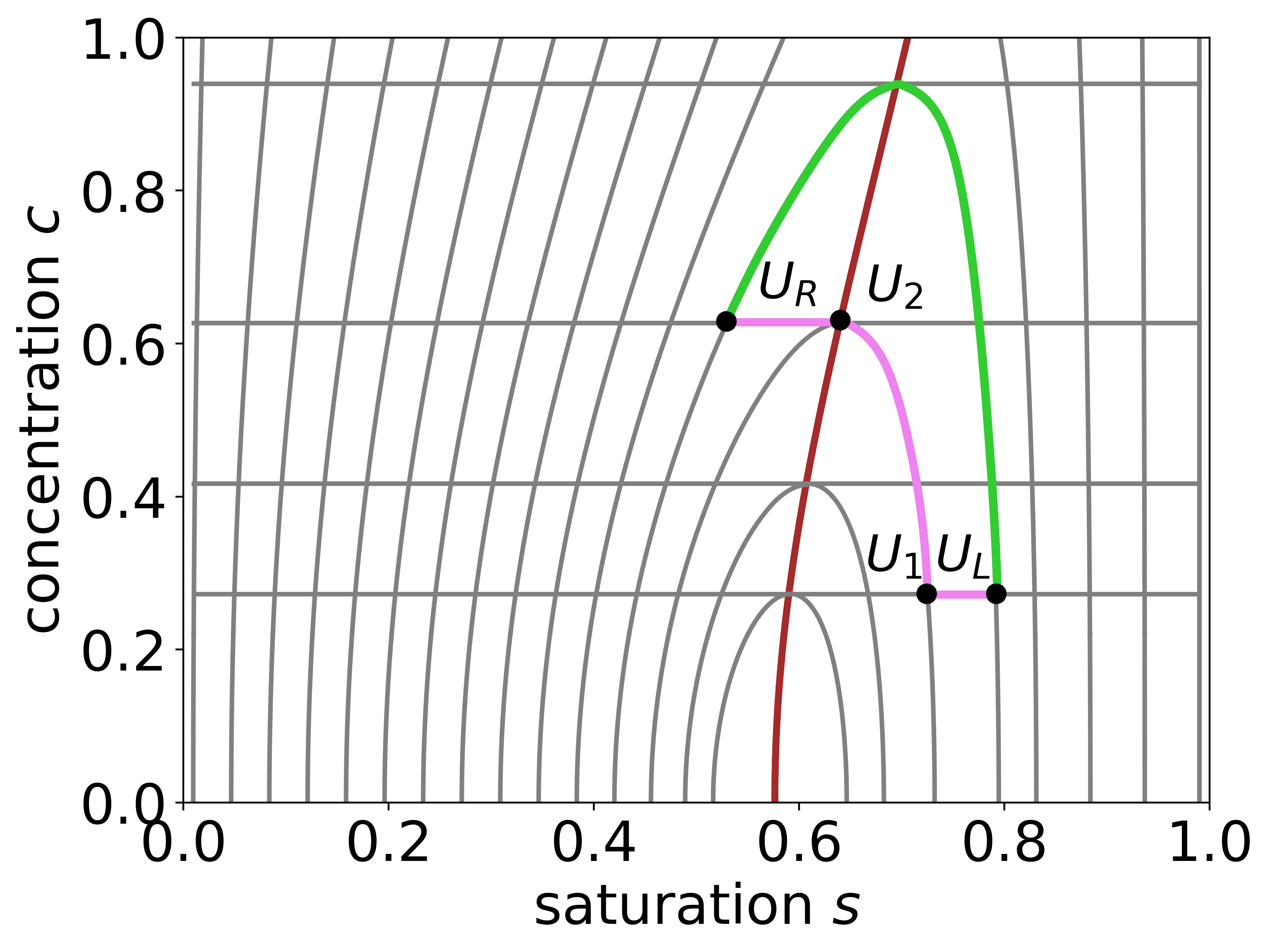}
        \includegraphics[width=0.45\textwidth]{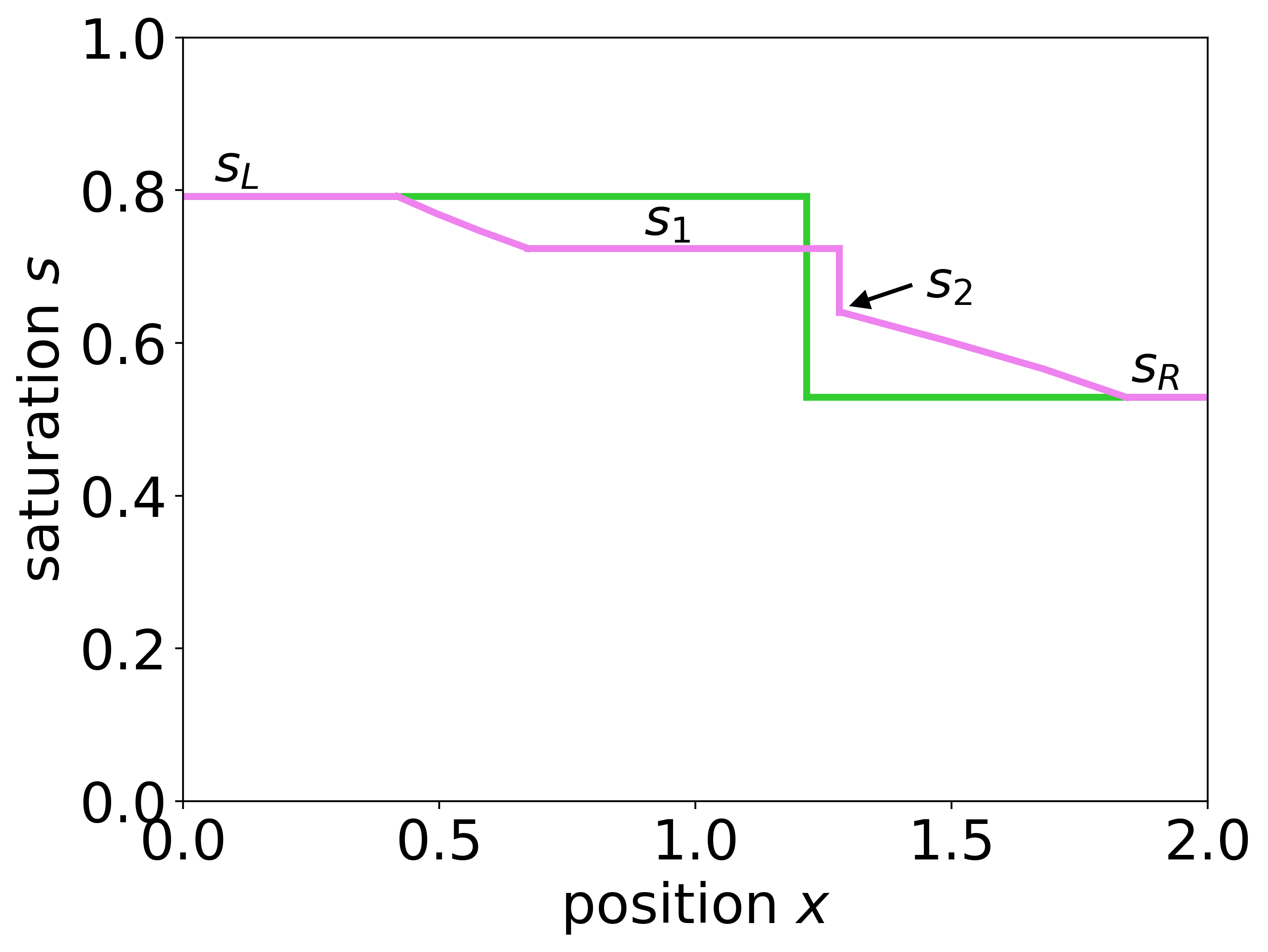}
        \\
    \qquad(a)\qquad\hfil \qquad\qquad\qquad(b)\hfil
    \caption{ (a) Horizontal lines are $s$-type branches of Hugoniot loci,
    and arches are $c$-type branches.
    The brown curve is the interior coincidence locus.
    The violet and green curves on both (a) and (b) represent
    two possible solutions of a particular
    Riemann problem, which exemplify multiplicity of solutions
    in the absence of an admissibility condition for contacts.}
    \label{fig:example-non-uniqueness}
\end{figure}

When $\alpha>0$,
the vanishing of expression~\eqref{eq:dlambda_c-ads}
defines the $c$-type inflection locus. Notice that the coincidence locus is contained in the inflection locus, as observed in Ref.~\cite{wahanik2011}. 
The $c$-type inflection locus is the union of the interior coincidence locus $\mathcal{C}_\alpha$
with the boundary line $s=0$ (because $f = 0$ only if $s = 0$).  
The $c$-type integral curve $\mathcal{I}_\alpha(U_-)$ through a state $U_-$
is the solution of the ODE
\begin{equation}
\label{eq:integral-U}
U_\eta=r^c_\alpha(U)
\end{equation}
that emanates from $U_-$ at $\eta = 0$ in both directions,
\emph{i.e.},\ for negative as well as positive $\eta$.
A $c$-type rarefaction wave is constructed,
as previously described,
from any segment within an integral curve
that does not cross the inflection locus
(meaning, in practice, the interior coincidence locus).
\begin{remark}
\label{rm:rare-continuous}
By virtue of the monotonicity assumption~\eqref{eq:f-cond-monotone},
$\mathcal{I}_\alpha(U_-)$ can be constructed
by solving the ODE
\begin{align}
\frac{dc}{ds} = \frac{\lambda^s-\lambda^c_\alpha}{-f_c}
\label{eq:arch-ODE}
\end{align}
to find $c$ as a function of $s$.
Thus, $\mathcal{I}_\alpha(U_-)$ takes the form of an arch
with horizontal tangent at the interior coincidence locus,
as illustrated in Fig.~\ref{fig:example-non-uniqueness}.
Because the right-hand side of Eq.~\eqref{eq:arch-ODE}
is smooth in $\alpha$ as well as $(s, c)$,
$\mathcal{I}_\alpha(U_-)$ varies smoothly in $\alpha$,
and it approaches $\mathcal{I}_0(U_-)$ as $\alpha \to 0^+$.
\end{remark}

\subsubsection{Shock waves}
For system~\eqref{eq:conslaw_adsorption},
the Rankine-Hugoniot condition~(\ref{eq:Rankine-Hugoniot})
for states $U_- = (s_-, c_-)$ and $U_+ = (s_+, c_+)$
is the pair of equations
\begin{align}
- \sigma\,\Delta s + \Delta f &= 0, \\
- \sigma \left[\overline{s}\,\Delta c + \alpha\,\Delta a_1\right]
+ \overline{f}\,\Delta c + \left[- \sigma\,\Delta s
+ \Delta f \right] \overline{c} &= 0.
\end{align}
The notation used here is $\Delta J := J_+ - J_-$
and $\overline{J} := (J_- + J_+)/2$,
and we have invoked the identity
$\Delta(J\,K) = (\Delta J)\,\overline{K} + \overline{J}\,\Delta K$.

Let us define the smooth function
\begin{equation}
\langle a'_1 \rangle
:= \int_0^1 a'_1\left((1-\theta)\,c_- + \theta\,c_+\right) d\theta
\end{equation}
of $c_-$ and $c_+$,
which satisfies $\Delta a_1 = \langle a'_1 \rangle \, \Delta c$.
Using this identity and the first Rankine-Hugoniot equation,
the second one becomes
\begin{align}
\left\{- \sigma \left[\overline{s} + \alpha\,\langle a'_1 \rangle\right]
+ \overline{f}\right\} \Delta c = 0.
\end{align}
We conclude that either $c_+ = c_-$ or the quantity in braces vanishes.
In the former case,
the first Rankine-Hugoniot equation
is the jump condition of the Buckley-Leverett scalar conservation law
for a fixed value of $c$.
In the latter case,
adding and subtracting $1/2$ times the first equation
to the vanishing of the braces
gives the equivalent pair of equations
\begin{align}
- \sigma \left[s_\pm + \alpha\,\langle a'_1 \rangle\right]
+ f(s_\pm, c_\pm) = 0.
\end{align}
Recalling our assumptions that $\alpha \ge 0$ and $a'_1 > 0$,
we observe that the quantity in brackets is positive when $s_\pm > 0$.

Therefore, the Rankine-Hugoniot condition has two types of solutions:
\begin{itemize}
\item a saturation (or $s$-type) discontinuity,
for which $c_+ = c_- =: c$ and
\begin{equation}
- \sigma\,(s_+ - s_-) + f(s_+, c) - f(s_-, c) = 0;
\end{equation}
\item a concentration (or $c$-type) discontinuity, which satisfies
\begin{equation}
\label{eq:c-type-hugoniot}
\frac{f(s_+, c_+)}{s_+ + \alpha\,\langle a'_1 \rangle}
= \sigma = \frac{f(s_-, c_-)}{s_- + \alpha\,\langle a'_1 \rangle}.
\end{equation}
\end{itemize}
The branches of the Hugoniot locus of a state $U_-$
comprising $s$-type and $c$-type discontinuities
are denoted $\mathcal{H}_\alpha^s(U_-)$ and $\mathcal{H}_\alpha^c(U_-)$,
respectively.

\begin{remark}
\label{rm:hugoniot-continuous}
When $\alpha = 0$,
a $c$-type discontinuity is a contact discontinuity~\eqref{eq:double-sonic},
so that $\mathcal{H}_0^c(U_-)$ coincides with
the integral curve $\mathcal{I}_0(U_-)$.
When $\alpha > 0$,
$\mathcal{H}_\alpha^c(U_-)$ is the zero-set of
\[
\mathcal{F} := f(s_+, c_+)/\left[s_+ + \alpha\,\langle a'_1 \rangle\right]
- f(s_-, c_-)/\left[s_- + \alpha\,\langle a'_1 \rangle\right],
\]
which is a smooth function of $\alpha$ as well as $(s_+, c_+)$.
Moreover, when evaluated at $\alpha = 0$,
the derivative $\mathcal{F}_{c_+}$ equals $f_c(s_+, c_+) / s_+$,
which is nonzero for all $s_+ > 0$.
Therefore, the zero-set is a smooth curve,
with $c$ parametrized by $s$,
so long as $\alpha > 0$ is sufficiently small.
Thus, $\mathcal{H}_\alpha^c(U_-)$ varies smoothly with $\alpha$
and approaches $\mathcal{H}^c_0(U_-)$ as $\alpha \to 0^+$.
Also, for each $c_+$,
$\mathcal{F}_{s_+} = \left[f_s(s_+, c_+) - \sigma\right]
/\left[s_+ + \alpha\,\langle a'_1 \rangle\right]$
vanishes at a unique value of $s_+$
(which exceeds $s_\text{coinc}(c_+)$ because $a''_1 < 0$),
so that $\mathcal{H}_\alpha^c(U_-)$ takes the form of an arch
with horizontal tangent near to the interior coincidence locus.
\end{remark}

\subsubsection{Admissibility}
A solution of a Riemann problem is an assembly
of $s$-wave groups, $c$-waves,
and constant states ordered by speed.
We shall use the notation $U\xrightarrow{s}U'$
(respectively, $U\xrightarrow{c}U'$)
to denote an $s$-wave group (resp., a $c$-wave) leading,
in the direction of increasing speed,
from a state $U$ to another state $U'$.

As is well-known for the Buckley-Leverett conservation law,
an $s$-wave group can be a rarefaction wave, a shock wave,
or a composite wave (\emph{i.e.},\ a sonic shock wave adjoining a rarefaction fan).
An $s$-wave group is required to
satisfy the Ole\u{\i}nik admissibility criterion~\cite{Ole57,Ole63};
equivalently,
each $s$-type shock wave in the solution must satisfy
the vanishing viscosity criterion~\cite{Hop69}.

Likewise, if $\alpha > 0$,
the vanishing-viscosity admissibility criterion can be
applied to $c$-type shock waves~\cite{johansen1988}.
In contrast, when $\alpha = 0$,
the admissibility of a $c$-wave
(\emph{i.e.},\ a contact discontinuity)
does not have a foundation based on a physical phenomenon.
Resolving this issue is the main purpose of this paper.

\subsubsection{Nonuniqueness of Riemann solutions}
Unless we impose some admissibility criterion
on $c$-waves in model $M_0$,
Riemann problems can have multiple solutions,
as illustrated by the following example
depicted in Fig.~\ref{fig:example-non-uniqueness}.
Take two points $U_L=(s_L,c_L)\in\{\lambda^s<\lambda^c\}$
and $U_R=(s_R,c_R)\in\{\lambda^s>\lambda^c\}$
that lie on the same contact curve
but are on opposite sides of the interior coincidence locus.
There exist at least two solutions to the Riemann problem:
\begin{align*}
    U_L\xrightarrow{c}U_R\qquad\text{ and }\qquad
    U_L\xrightarrow{s}U_1\xrightarrow{c}U_2\xrightarrow{s}U_R.
\end{align*}
(In fact, there are infinitely many solutions.)
Here $U_2=(s_2,c_R)$ is a point on $\mathcal{C}_0$
and $U_1=(s_1,c_L)\in\{\lambda^s<\lambda^c\}$ is
the unique point on the same contact curve as $U_2$ with concentration $c_L$.

\subsubsection{Layout of solutions to Riemann problems}
\label{subsubsec:RP-polymer}
The types of waves that appear in a Riemann solution,
as well as their relative speeds,
depend on the locations of both $U_L$ and $U_R$
with respect to certain curves.
Let us recount how solutions
of the general Riemann problem~\eqref{eq:Riemann-problem}
are laid out in the $(s, c)$-plane.

\begin{figure}[ht]
    \centering
\includegraphics[width=0.37\textwidth]{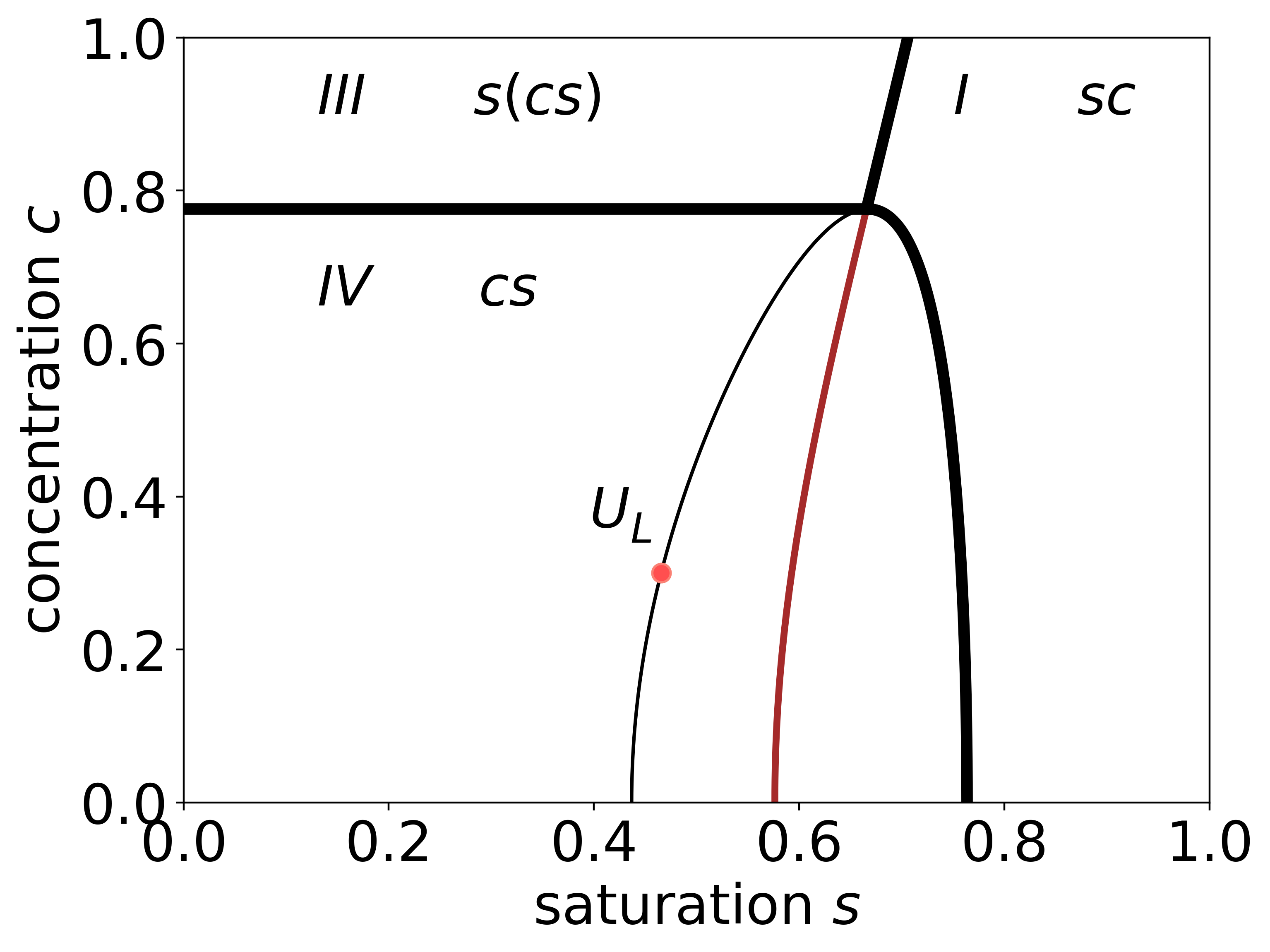}
\includegraphics[width=0.37\textwidth]{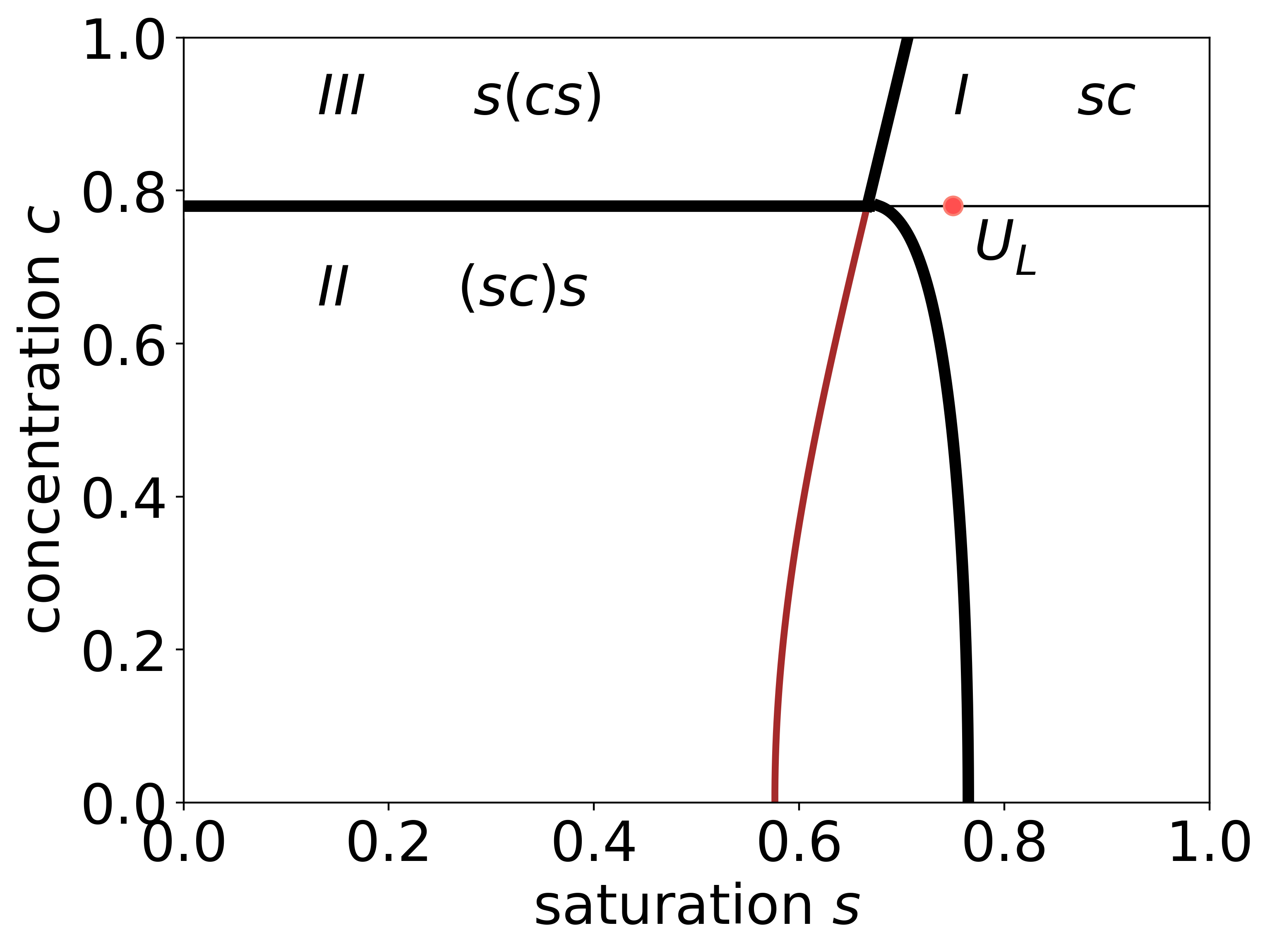} \\
    (a)\qquad\qquad\qquad\qquad\quad\quad\qquad\qquad(b)\hfil
    \caption{Polymer model~$M_0$:
    decomposition into $U_R$-regions~I, II, \emph{etc}.,
    along with the solution structure,
    depending on the location of the left state~$U_L$:
    (a) $U_L\in\{\lambda^s>\lambda^c\}$;
    (b) $U_L\in\{\lambda^s<\lambda^c\}$.}
    \label{fig:scs-structure}
\end{figure}

\begin{figure}[ht]
    \centering
\includegraphics[width=0.37\textwidth]{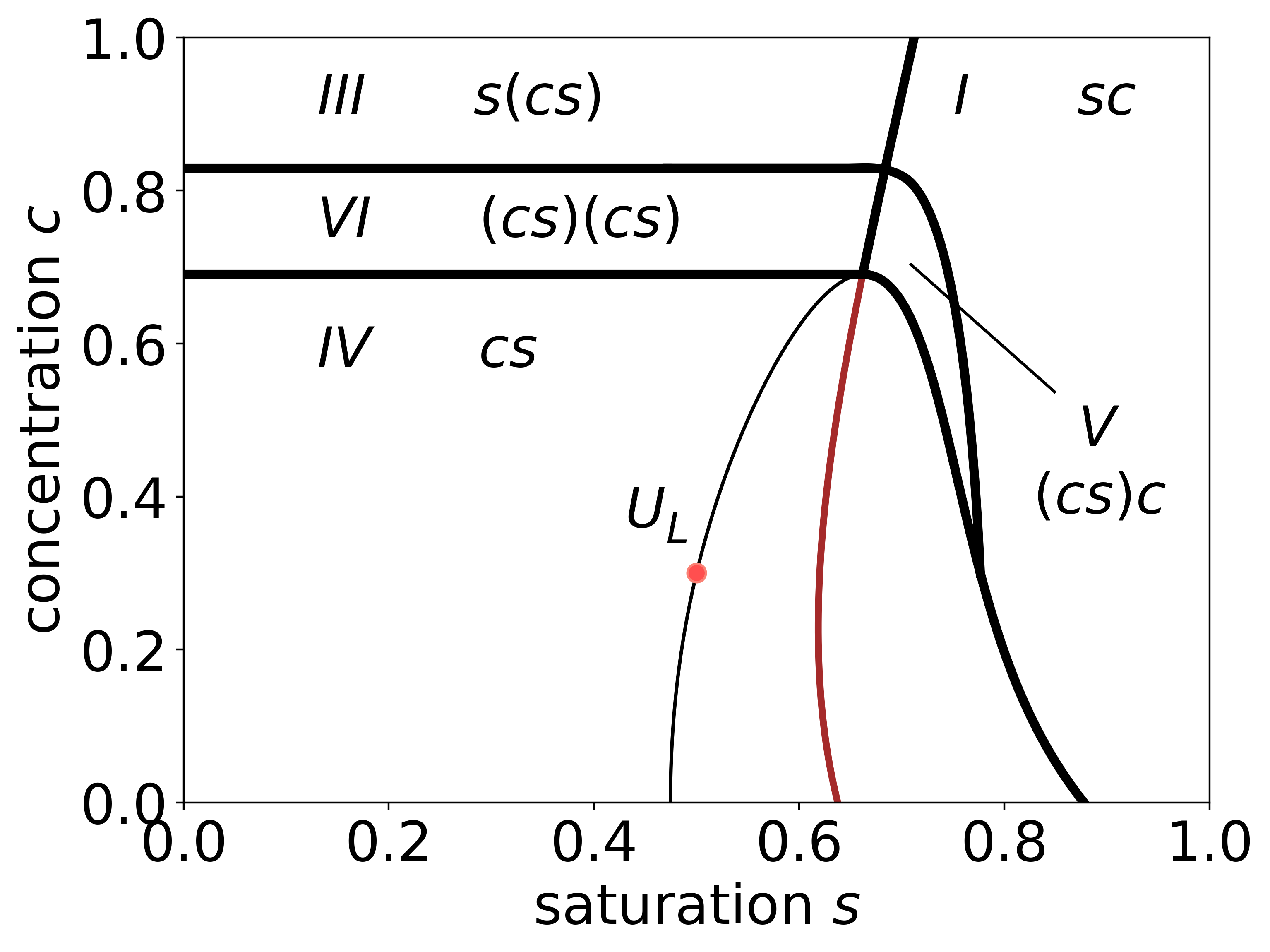}
\includegraphics[width=0.37\textwidth]{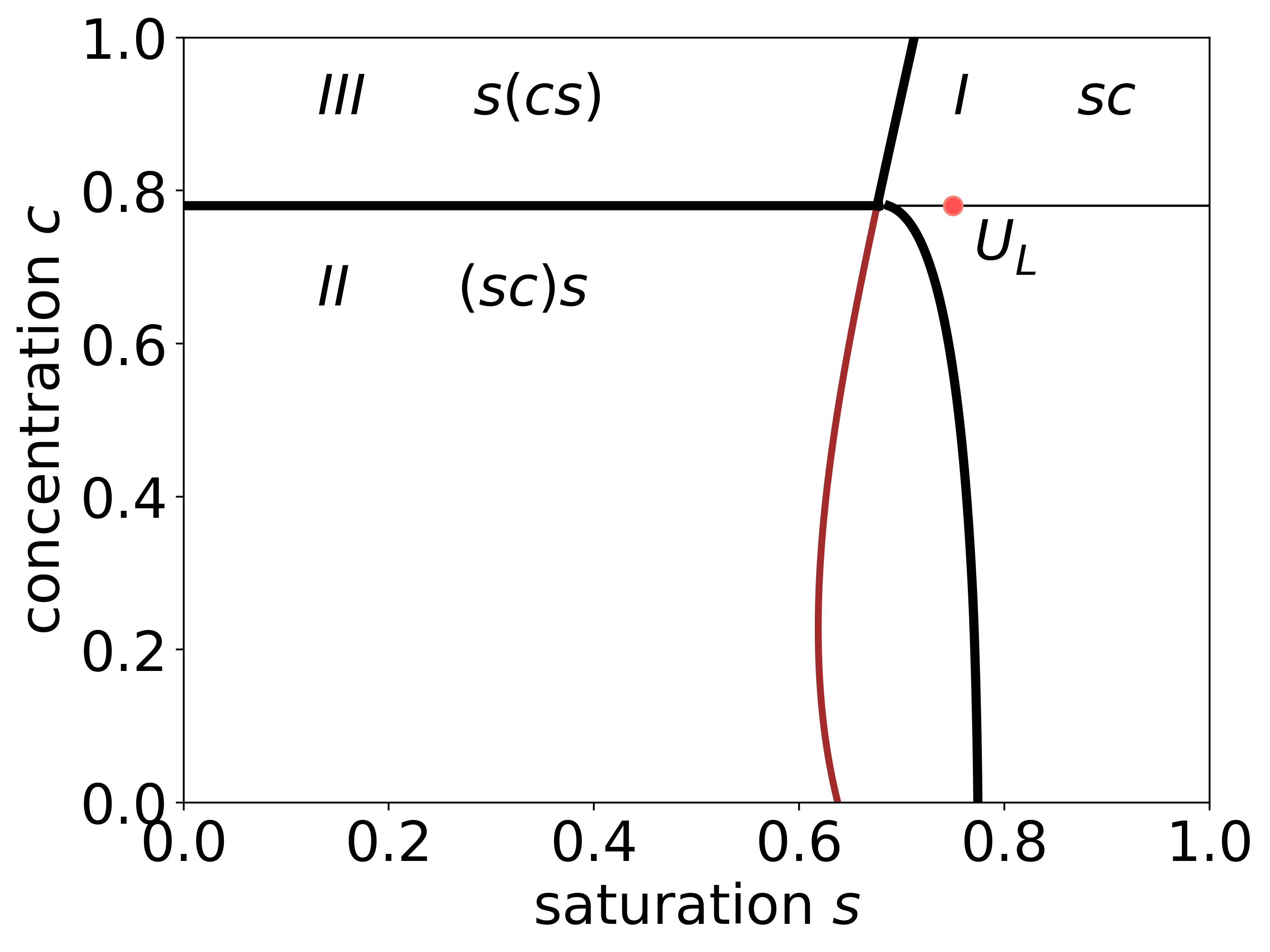} \\
    (a)\qquad\qquad\qquad\qquad\quad\quad\qquad(b)\hfil
    \caption{Polymer model~$M_\alpha$ with $\alpha>0$:
    decomposition into $U_R$-regions~I, II, \emph{etc}.,
    along with the solution structure,
    depending on the location of the left state~$U_L$:
    (a) $U_L\in\{\lambda^s>\lambda^c\}$;
    (b) $U_L\in\{\lambda^s<\lambda^c\}$.}
    \label{fig:scs-structure-ads}
\end{figure}

In Figs.~\ref{fig:scs-structure} and~\ref{fig:scs-structure-ads}
we present the Riemann solutions when
(a)~$U_L\in \{\lambda^s>\lambda^c\}$
and (b)~$U_L\in\{\lambda^s<\lambda^c\}$.
(If $U_L$ lies on the interior coincidence locus,
Riemann solutions can be obtained
by taking appropriate limits of solutions in these cases.)
Figures~\ref{fig:scs-structure}(a) and~(b)
are for the polymer model without adsorption ($\alpha=0$),
whereas Figs.~\ref{fig:scs-structure-ads}(a) and~(b)
correspond to the model with adsorption ($\alpha>0$).

Let us illustrate how to interpret these diagrams
by explaining Fig.~\ref{fig:scs-structure}(b) for $M_0$ in more detail.
For a fixed left state $U_L\in\{\lambda^s<\lambda^c\}$,
there are three different sequences of waves that appear in Riemann solutions,
depending on the location of the right state $U_R$:
\begin{itemize}
    \setlength\itemsep{0.4em}
    \item $U_R$ in region~I:
    the solution is $U_L\xrightarrow{s} U_M \xrightarrow{c} U_R$;

    \item $U_R$ in region~II:
    the solution is
    $U_L\xrightarrow{s}U_1\xrightarrow{c}U_{M}\xrightarrow{s}U_R$,
    where the $s$-wave group leads from $U_L$ to the unique state
    $U_1=(s_1,c_L)\in\mathcal{C}_0$
    and the adjoining contact discontinuity
    leads from $U_1$ to the unique state
    $U_M=(s_2,c_R)\in\{\lambda^s>\lambda^c\}\cap\mathcal{H}_0^c(U_1)$;

    \item $U_R$ in region~III:
    the solution is
    $U_L\xrightarrow{s}U_M\xrightarrow{c}U_2\xrightarrow{s}U_R$,
    where the $s$-wave group leading to $U_R$
    begins at the unique state $U_2=(s_2,c_R)\in\mathcal{C}_0$
    and the adjoining contact leading to $U_2$ begins at the unique state
    $U_M=(s_1,c_L)\in\{\lambda^s<\lambda^c\}\cap\mathcal{H}_0^c(U_2)$.
\end{itemize}

The boundary between regions~I and~III
is the part of the interior coincidence locus
$\mathcal{C}_0$ with $c\geq c_L$,
whereas regions~II and~III are separated by the horizontal line $c = c_L$
in $\{\lambda^s>\lambda^c\}$
and the boundary between regions~I and~II is the part of the contact curve
containing the state $U_0=(s,c_L)\in\mathcal{C}_0$
that lies within $\{\lambda^s<\lambda^c\}$.

Figure~\ref{fig:scs-structure}(a) for $\alpha=0$
and Figs.~\ref{fig:scs-structure-ads}(a) and~(b) for $\alpha > 0$
are interpreted in a similar fashion.
See Ref.~\cite{johansen1988} for details.

\begin{figure}[ht]
    \centering
        \includegraphics[width=0.042\textwidth]{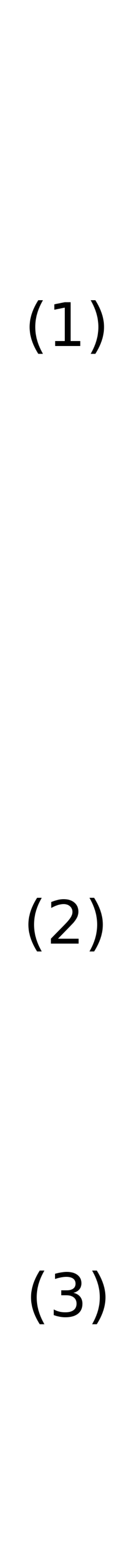}
        \includegraphics[width=0.31\textwidth]{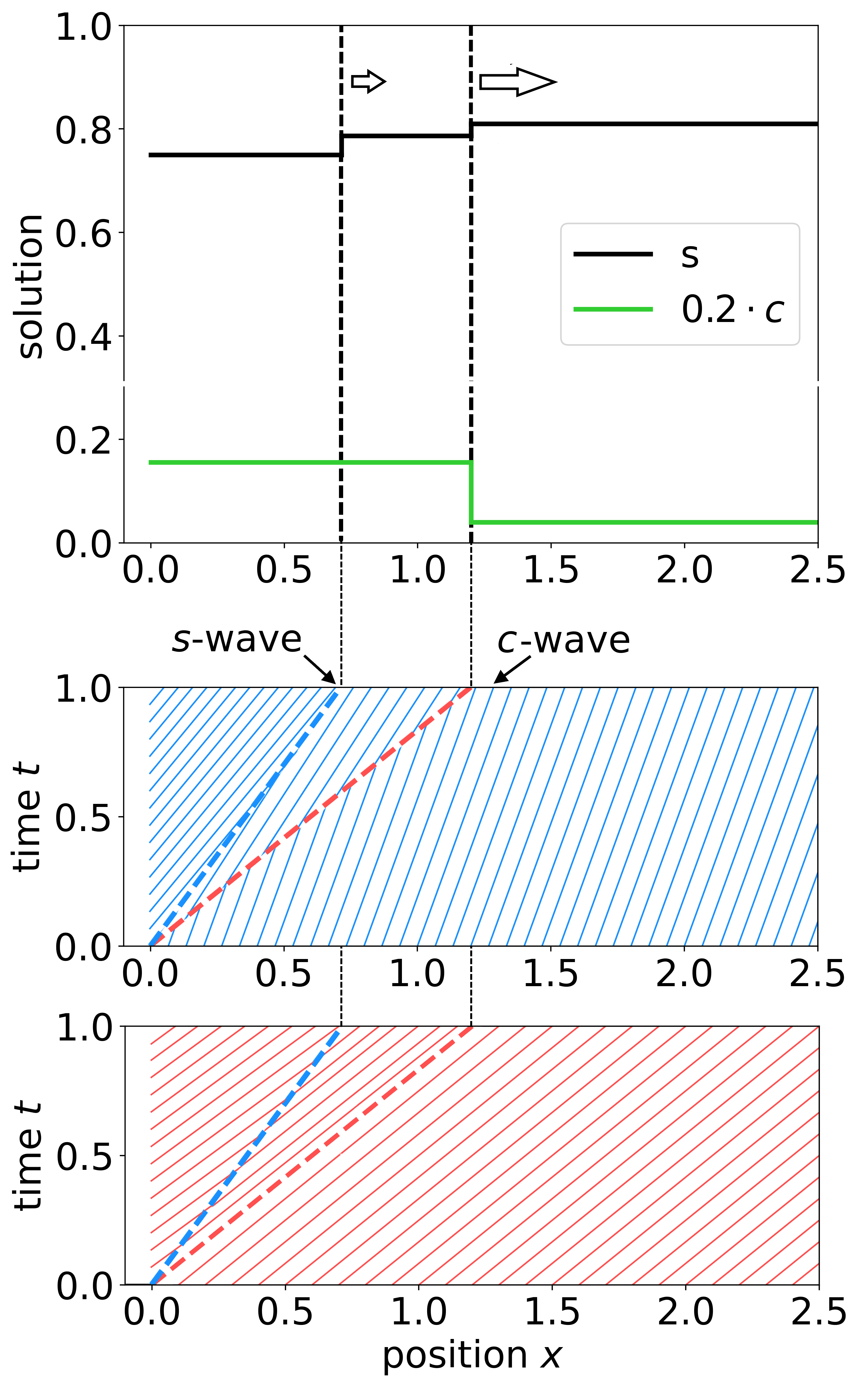}
        \includegraphics[width=0.31\textwidth]{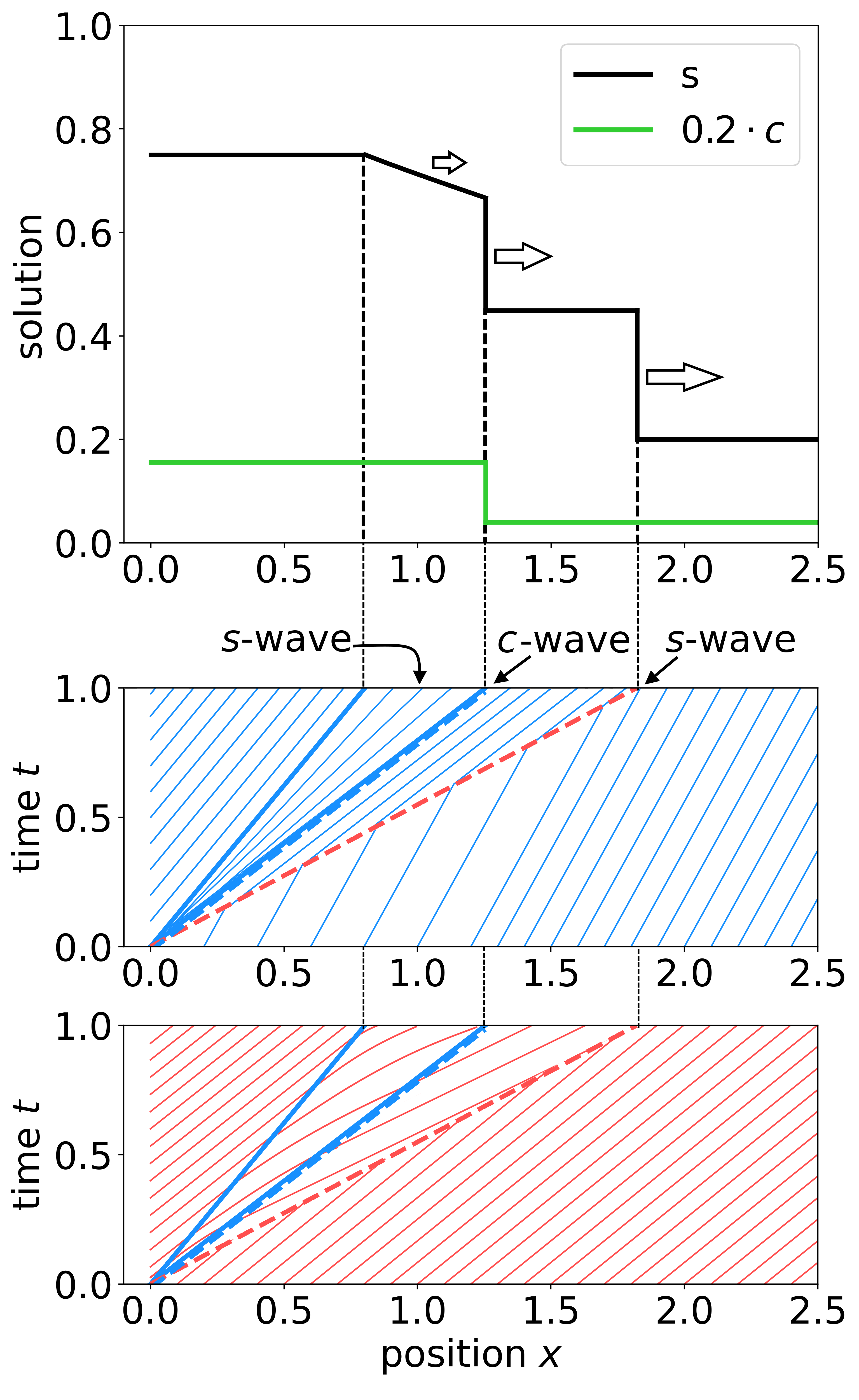}
        \includegraphics[width=0.31\textwidth]{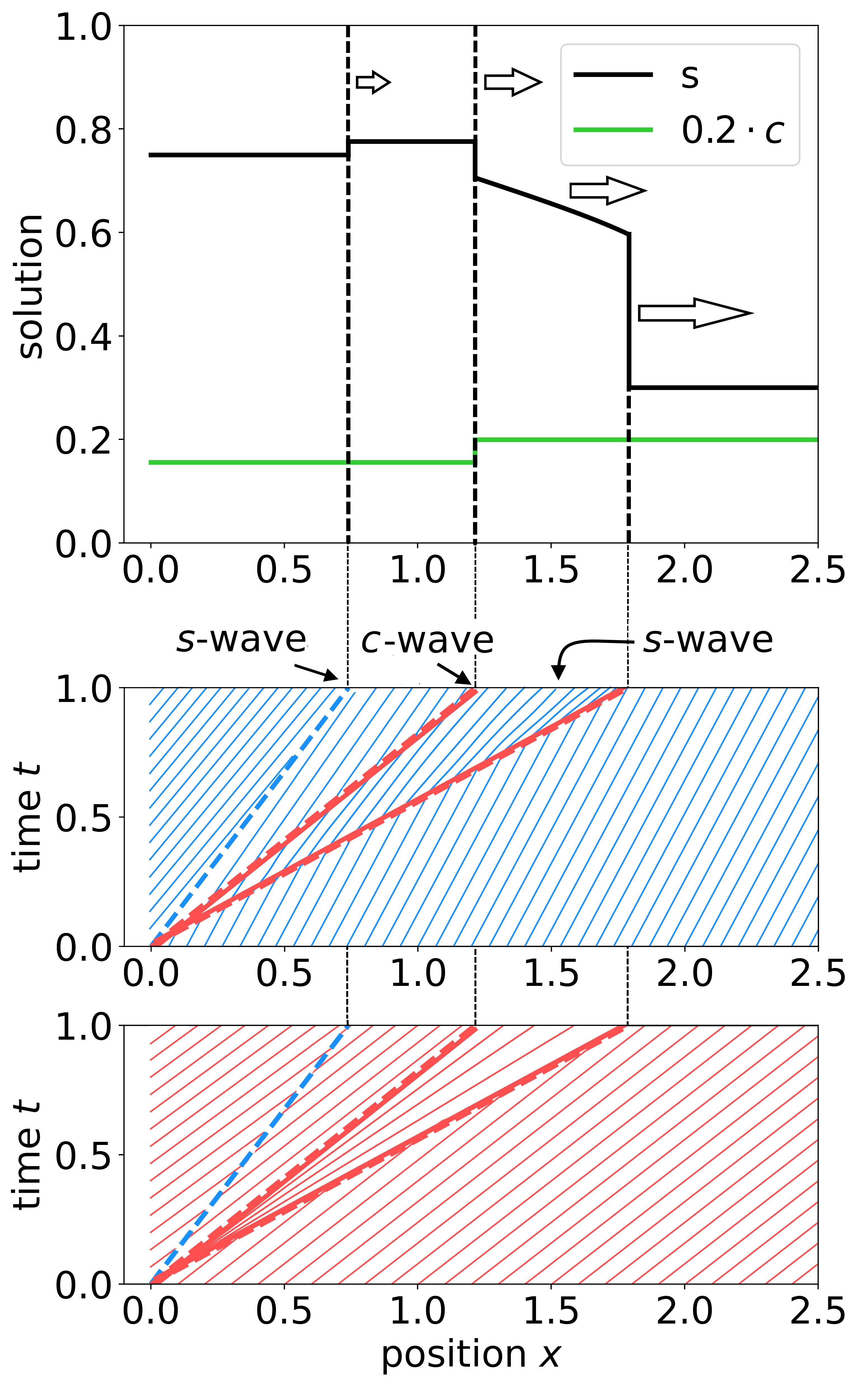}
        \\
    (a) $U_R$ in region~I\qquad\qquad
    (b) $U_R$ in region~II\qquad\qquad
    (c) $U_R$ in region~III
    \caption{Examples of three different sequences of waves
    that are solutions to a Riemann Problem
    for $U_L\in\{\lambda^s<\lambda^c\}$
    and $U_R$ in the respective regions~I, II, and~III
    from Fig.~\ref{fig:scs-structure}(b).}
    \label{fig:RP-sx-cx-xt}
\end{figure}

Another way of representing these three solutions
is to draw graphs of the functions $s(x/t)$ and $c(x/t)$ for fixed time $t > 0$.
These graphs appear in Figs.~\ref{fig:RP-sx-cx-xt}~(1a), (1b), and~(1c).
In addition, Figs.~\ref{fig:RP-sx-cx-xt}~(2a)-(2c) and~(3a)-(3c)
depict the trajectories, in the $(x,t)$-plane,
of shock waves
(drawn as dashed lines along which $x = \sigma\,t$)
and of characteristics
(\emph{i.e.},\ solutions of $dx/dt = \lambda(U(x/t))$, drawn as solid curves).
The blue and red coloring indicates whether the speed of the characteristic
is the smaller or larger eigenvalue, respectively,
as explained at the beginning of the next section.

\section{Equivalence of admissibility criteria}
\label{sec:proof}
Recall the standard ordering of the eigenvalues
of the characteristic matrix $A(U)$
in \eqref{eq:characteristic-matrix}
as $\lambda_1(U)<\lambda_2(U)$,
called the 1-family and 2-family characteristic speeds.
In polymer models,
$\lambda_1(U)$ equals $\lambda^c(U)$
when $U\in\{\lambda^s>\lambda^c\}$,
but equals $\lambda^s(U)$
when $U\in\{\lambda^s<\lambda^c\}$.

\begin{figure}[ht]
    \centering
    \includegraphics[width=0.24\textwidth]{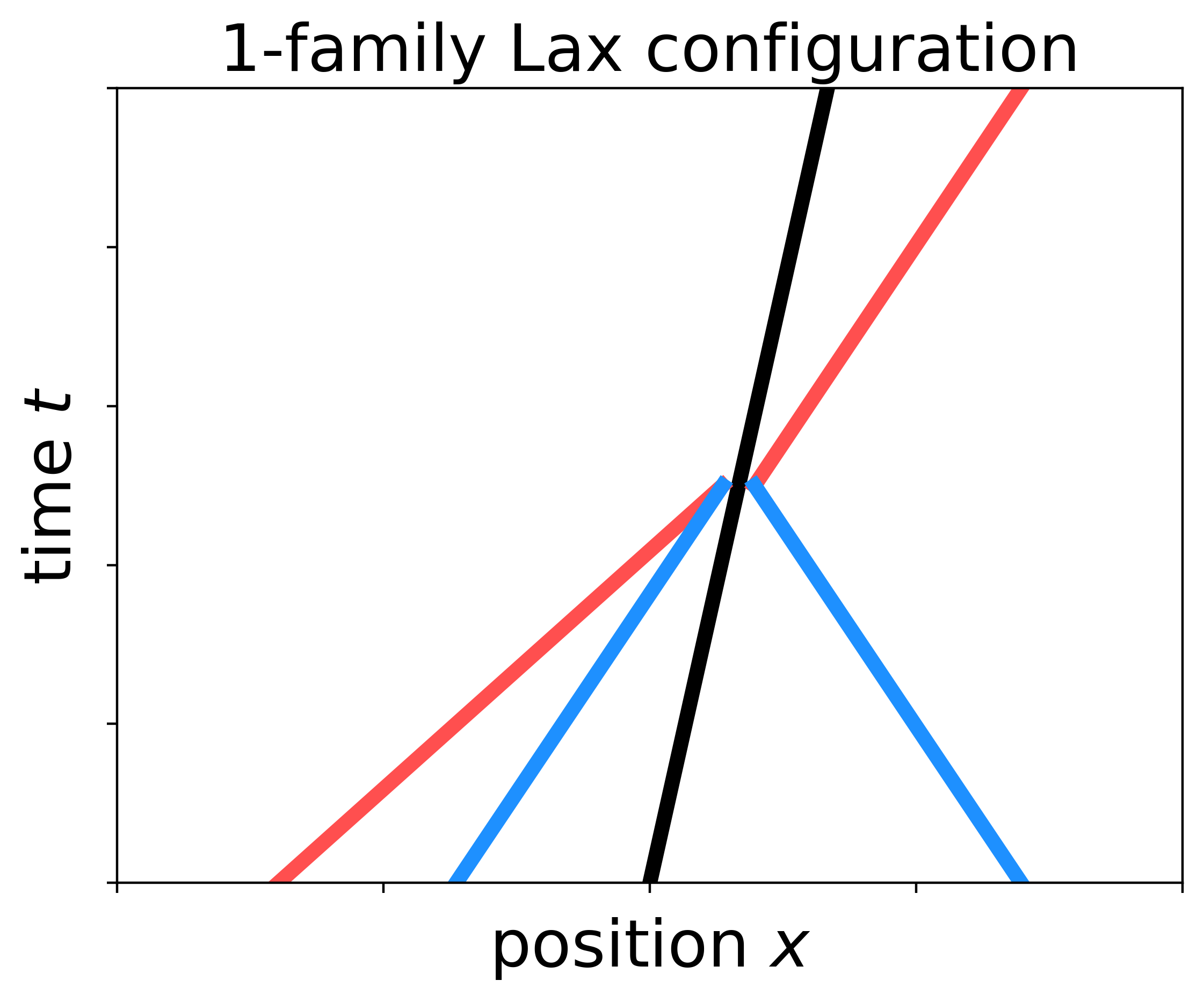}
    \includegraphics[width=0.24\textwidth]{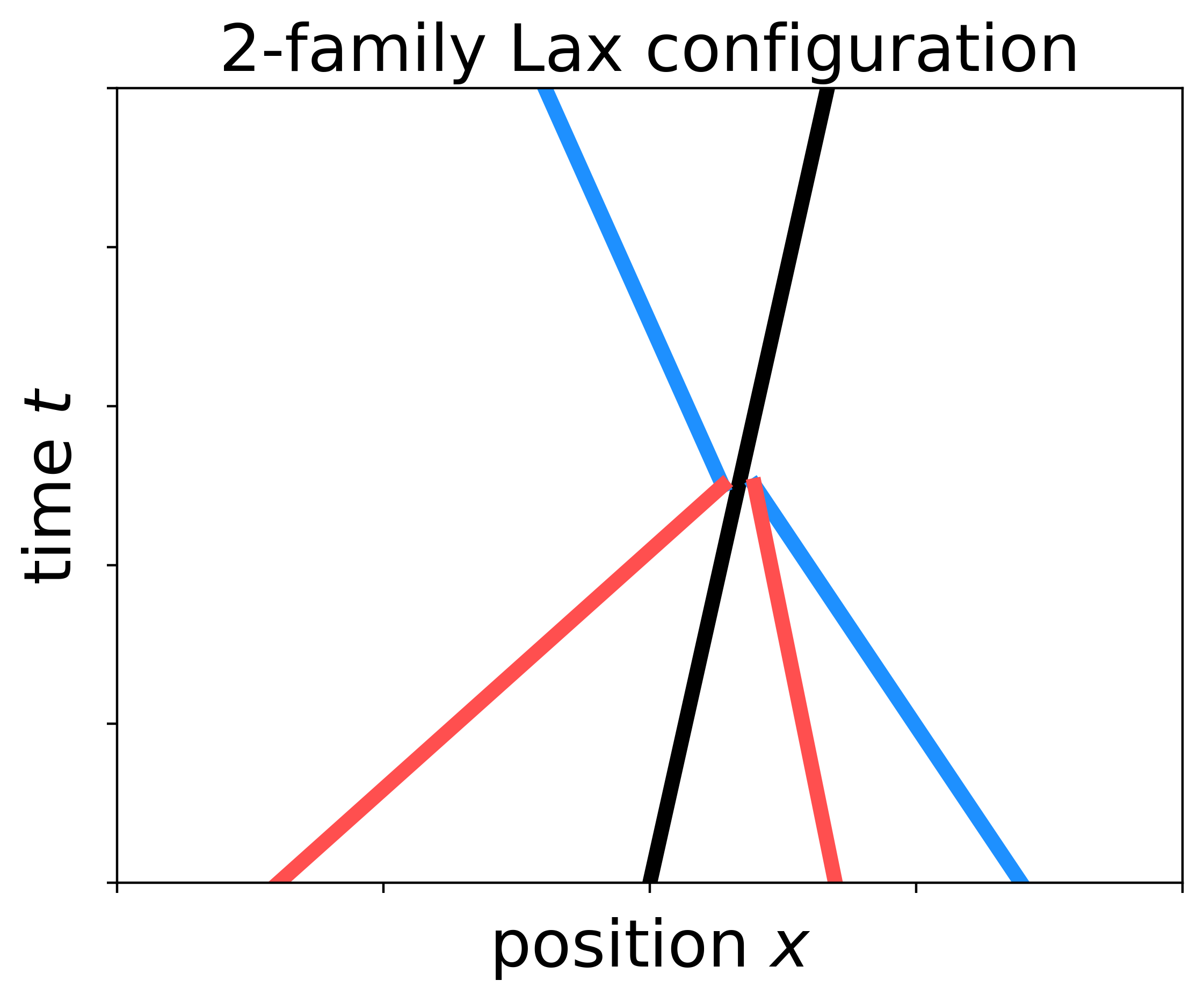}
    \includegraphics[width=0.24\textwidth]{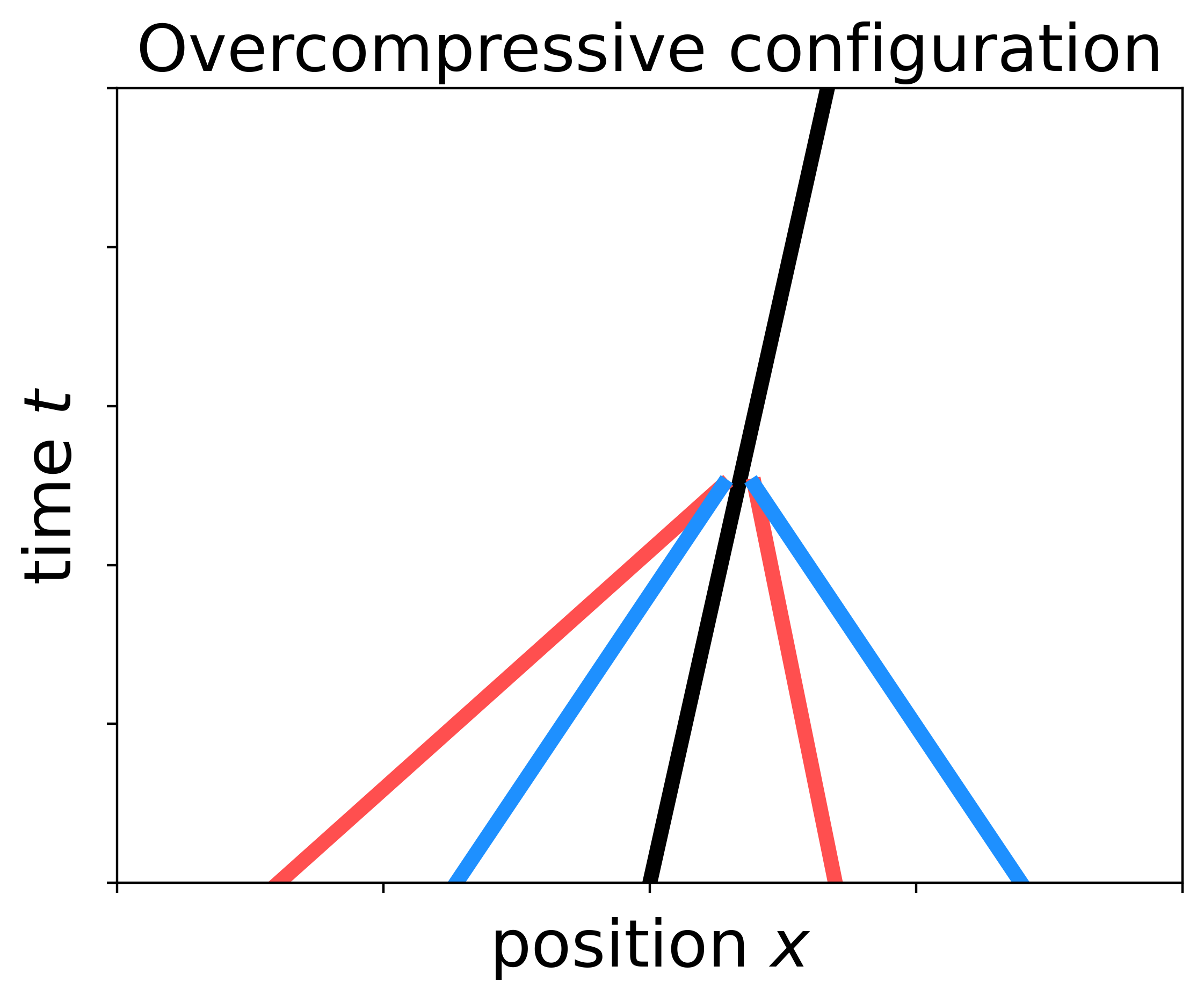}
    \includegraphics[width=0.24\textwidth]{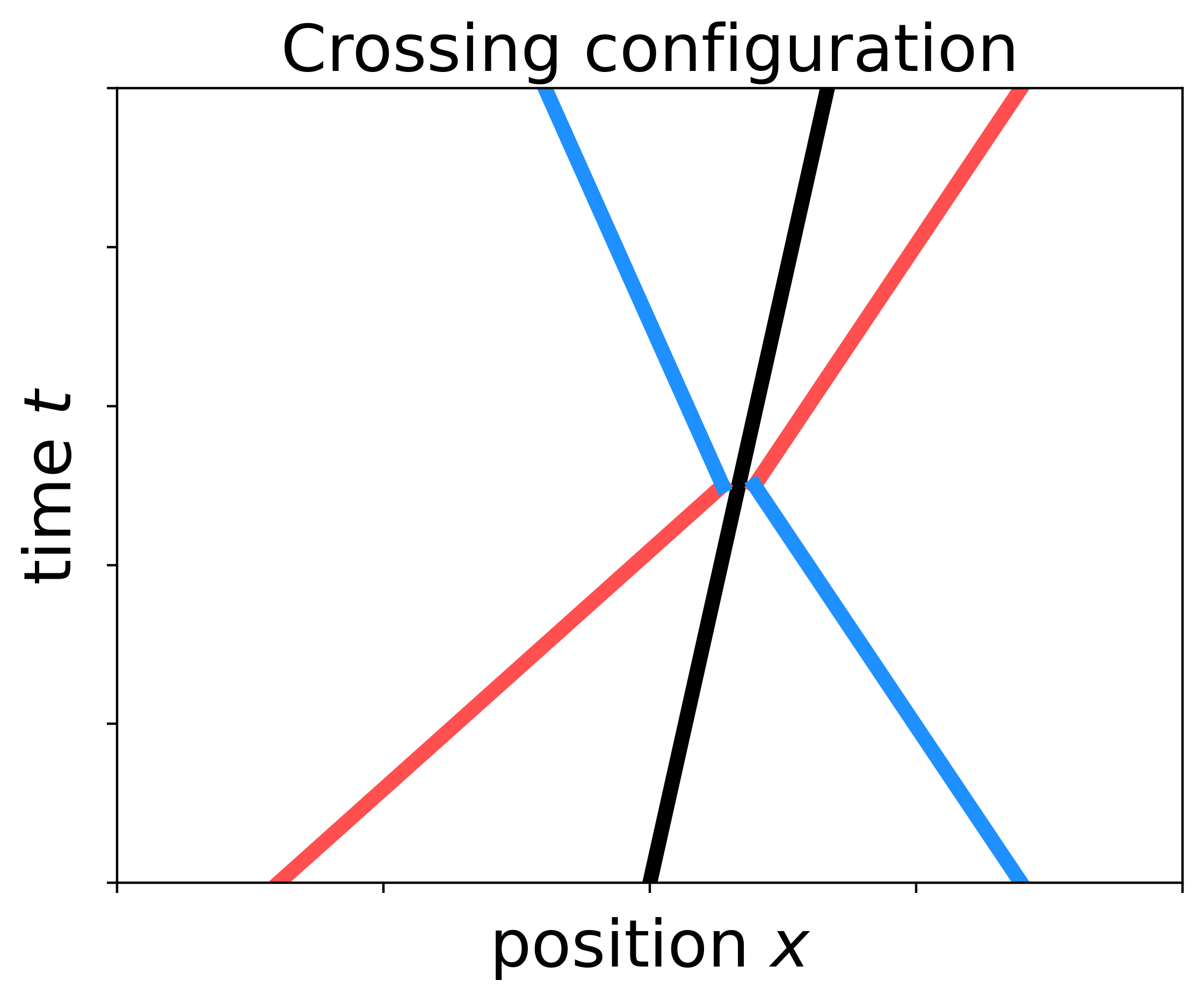}
    \caption{Four possible configurations of characteristic paths
    drawn in $(x,t)$-space: 1-family Lax, 2-family Lax, overcompressive,
    and crossing. The black line is the shock trajectory with
    speed $\sigma$. Blue curves are characteristic paths for the 1-family,
    and red curves are for the 2-family.}
    \label{fig:types_shocks}
\end{figure}

A standard way to classify a discontinuity
is based on the ordering of the characteristic speeds
on its two sides relative to its propagation speed $\sigma$,
\emph{i.e.},\ $\lambda_i(U_-) - \sigma$ and $\lambda_i(U_+) - \sigma$ for $i = 1, 2$.
See, for example,
Refs.~\cite{lax57}, \cite{transitional1990}, and~\cite[Chapter 8]{dafermos}.
Four of the possibilities, which we call
the 1-family Lax, 2-family Lax, overcompressive, and crossing
configurations of characteristic paths,
are depicted in Fig.~\ref{fig:types_shocks}:
\begin{itemize}
\item 1-family Lax: $\lambda_1(U_-)>\sigma>\lambda_1(U_+)$,
$\sigma <\lambda_2(U_-)$, and $\sigma<\lambda_2(U_+)$;
\item 2-family Lax: $\lambda_2(U_-)>\sigma>\lambda_2(U_+)$,
$\sigma >\lambda_1(U_-)$, and $\sigma>\lambda_1(U_+)$;
\item overcompressive: $\lambda_1(U_-)>\sigma>\lambda_1(U_+)$
and $\lambda_2(U_-)>\sigma>\lambda_2(U_+)$;
\item crossing:
$\lambda_2(U_-)>\sigma>\lambda_1(U_-)$
and $\lambda_1(U_+)<\sigma<\lambda_2(U_+)$.
\end{itemize}

\begin{figure}[ht]
    \centering
    \includegraphics[width=0.24\textwidth]{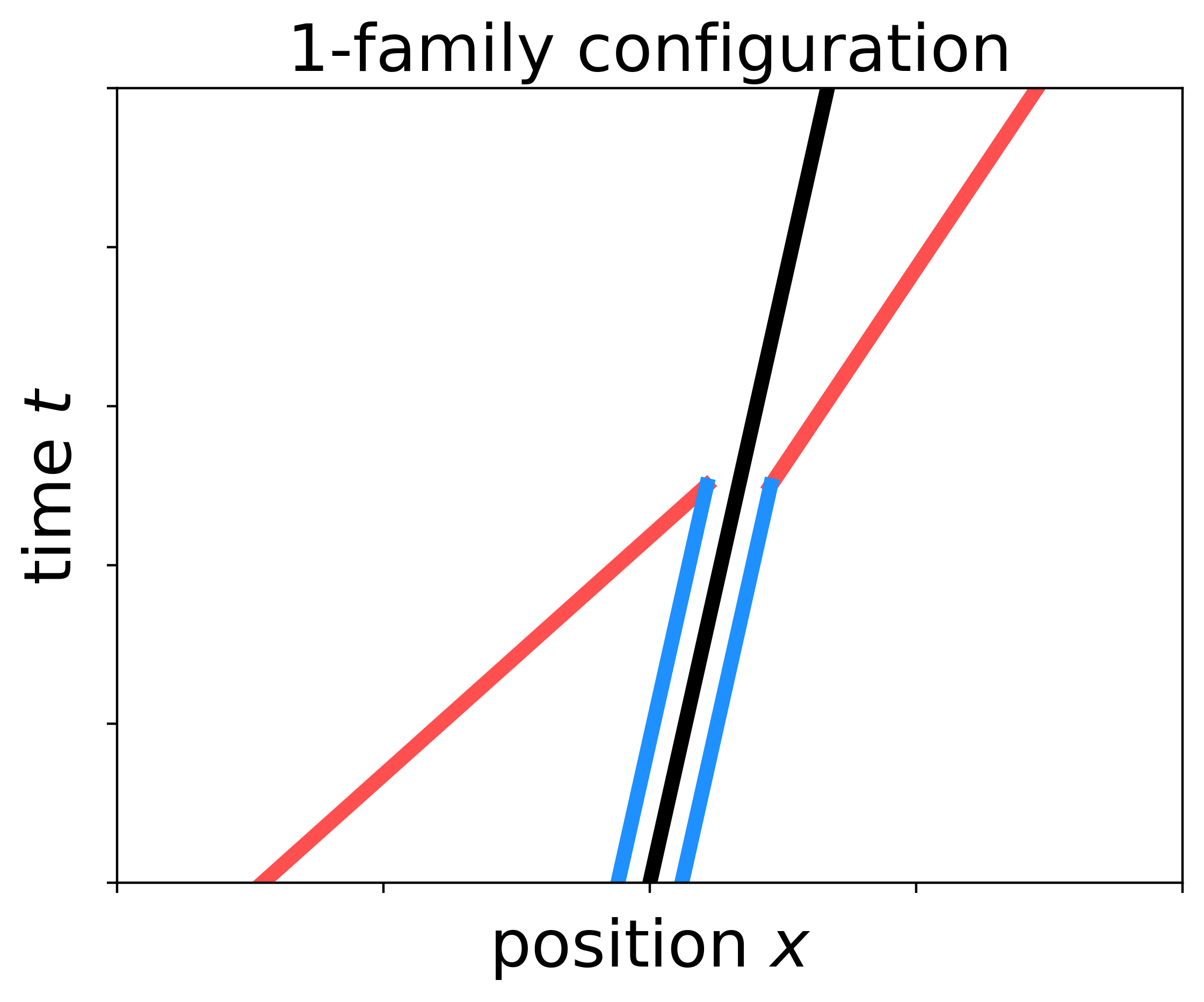}
    \includegraphics[width=0.24\textwidth]{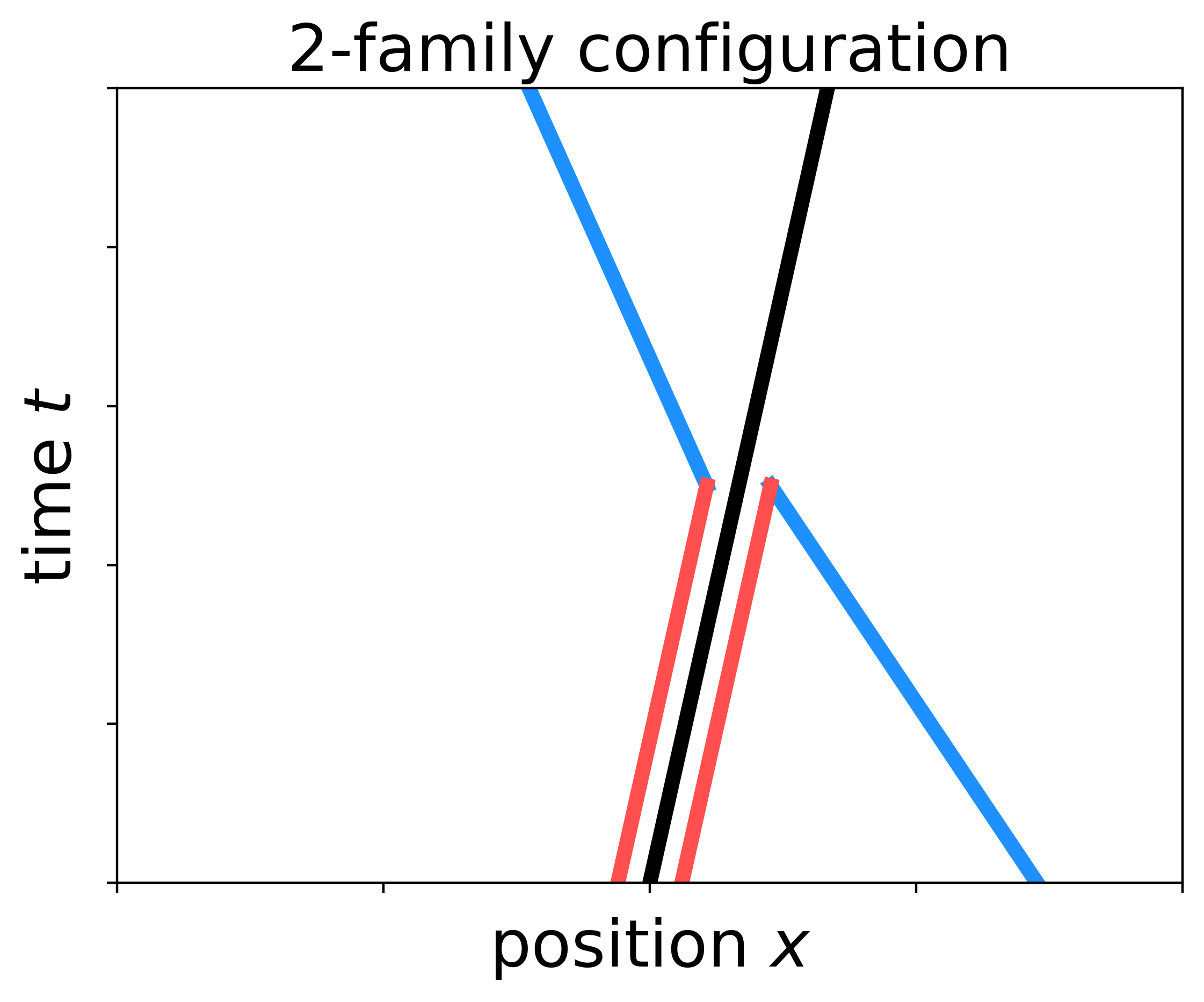}
    \includegraphics[width=0.24\textwidth]{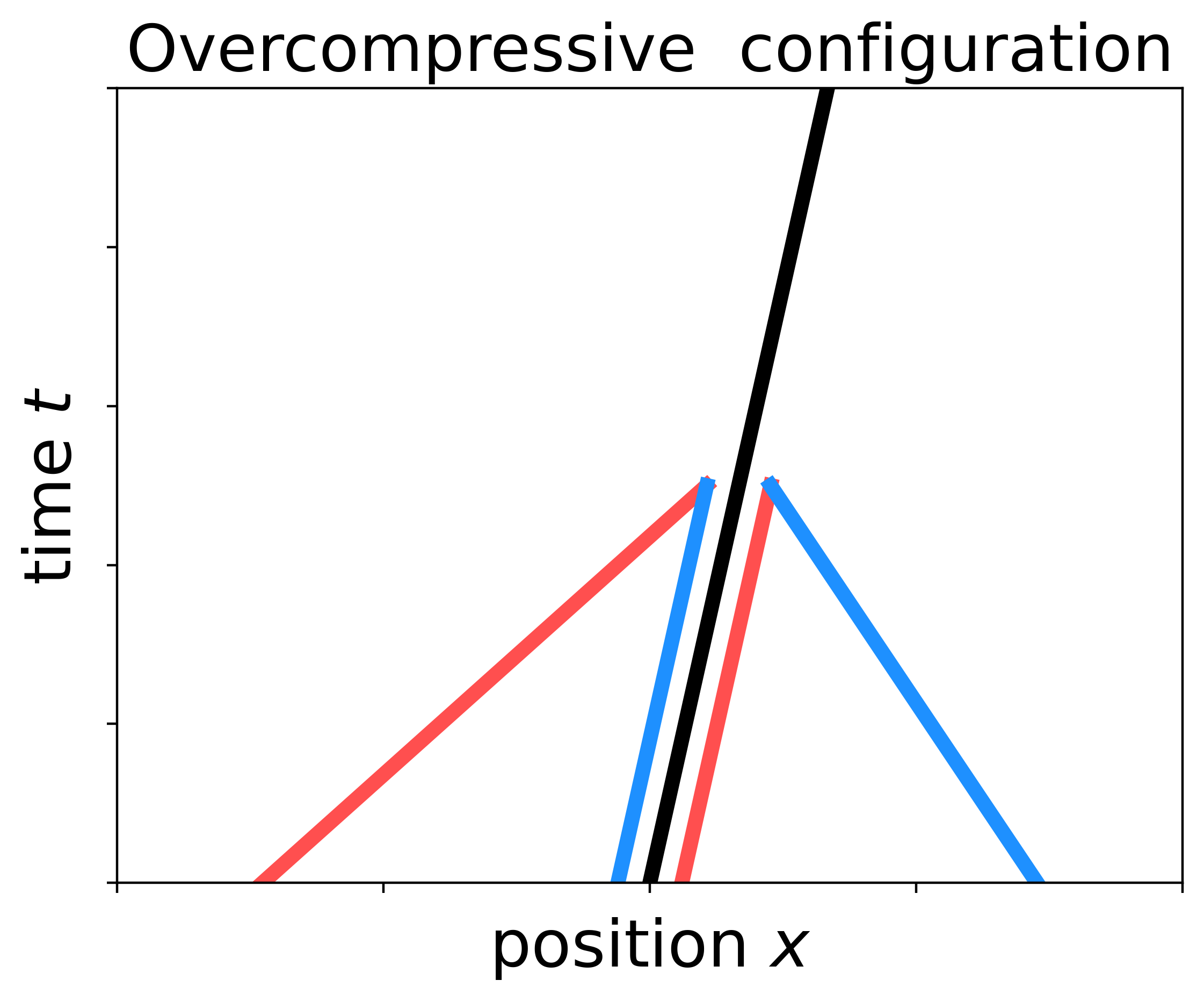}
    \includegraphics[width=0.24\textwidth]{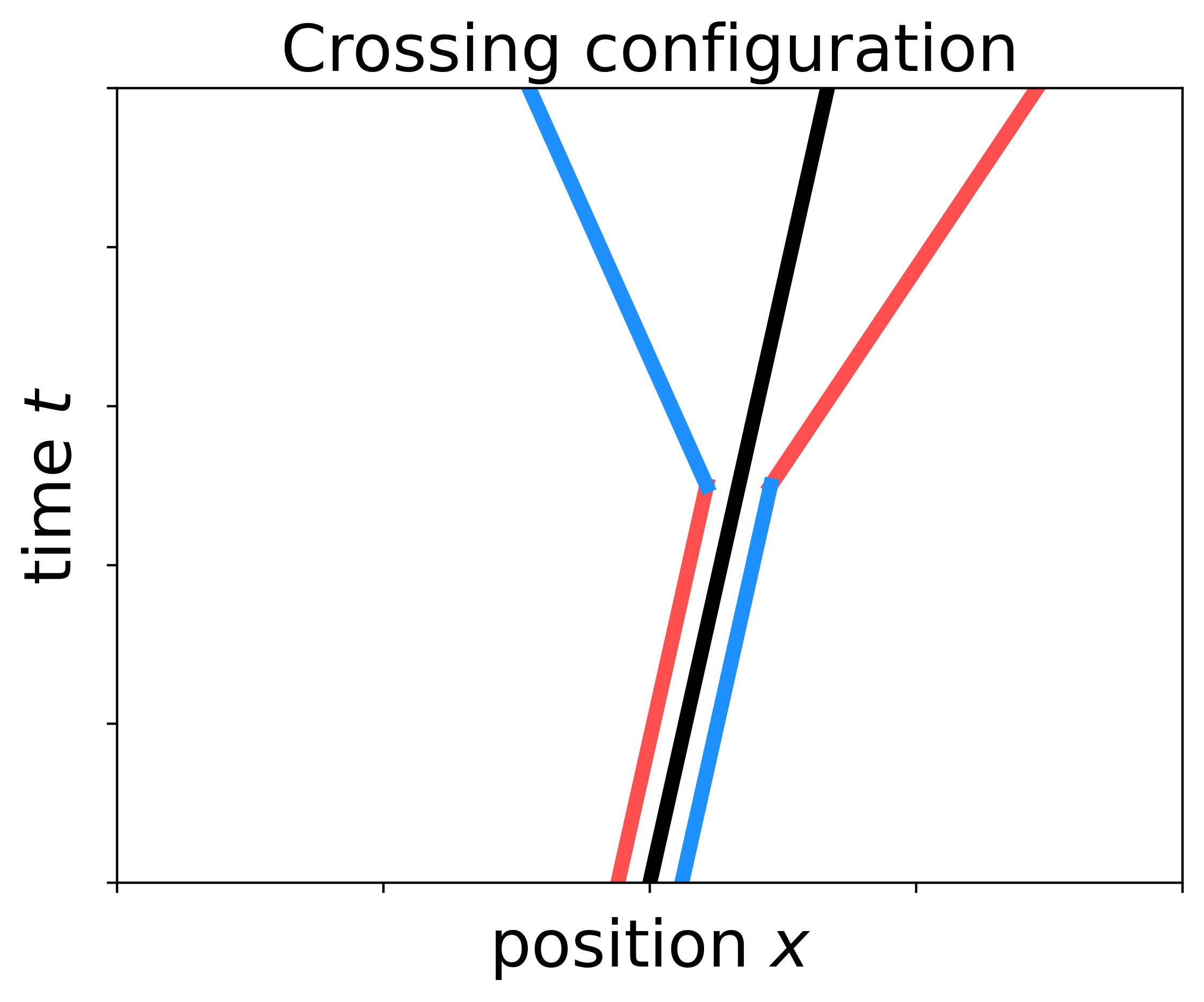}
    \caption{Four possible configurations of characteristic paths
    drawn in $(x,t)$-space for contact discontinuities:
    1-family, 2-family, overcompressive, and crossing.}
    \label{fig:types_contacts}
\end{figure}

It will be useful in what follows
to extend this classification to contact discontinuities.
\begin{definition}
\label{def:contact-char-configs}
We adopt the following terminology for configurations of characteristic paths
for contact discontinuities,
depicted in Fig.~\ref{fig:types_contacts}:
\begin{itemize}
\item 1-family: $\lambda_1(U_-)=\sigma=\lambda_1(U_+)$,
$\sigma <\lambda_2(U_-)$, and $\sigma<\lambda_2(U_+)$;
\item 2-family: $\lambda_2(U_-)=\sigma=\lambda_2(U_+)$,
$\sigma >\lambda_1(U_-)$, and $\sigma>\lambda_1(U_+)$;
\item overcompressive: $\lambda_1(U_-)=\sigma>\lambda_1(U_+)$
and $\lambda_2(U_-)=\sigma>\lambda_2(U_+)$;
\item crossing:
$\lambda_2(U_-)>\sigma=\lambda_1(U_-)$
and $\lambda_1(U_+)=\sigma<\lambda_2(U_+)$.
\end{itemize}
\end{definition}

\begin{figure}[ht]
    \centering
    \includegraphics[width=0.4\textwidth]{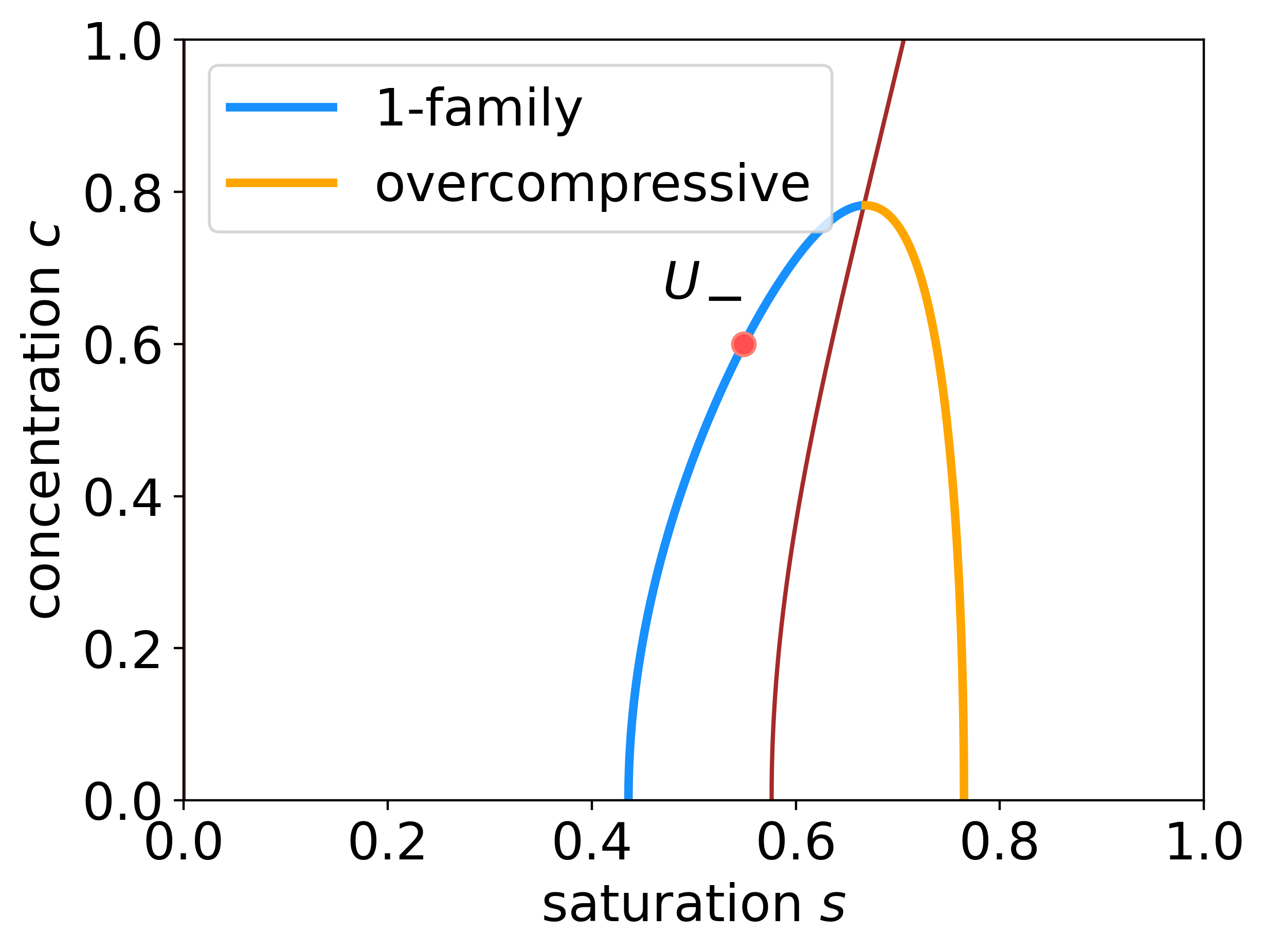}
    \includegraphics[width=0.4\textwidth]{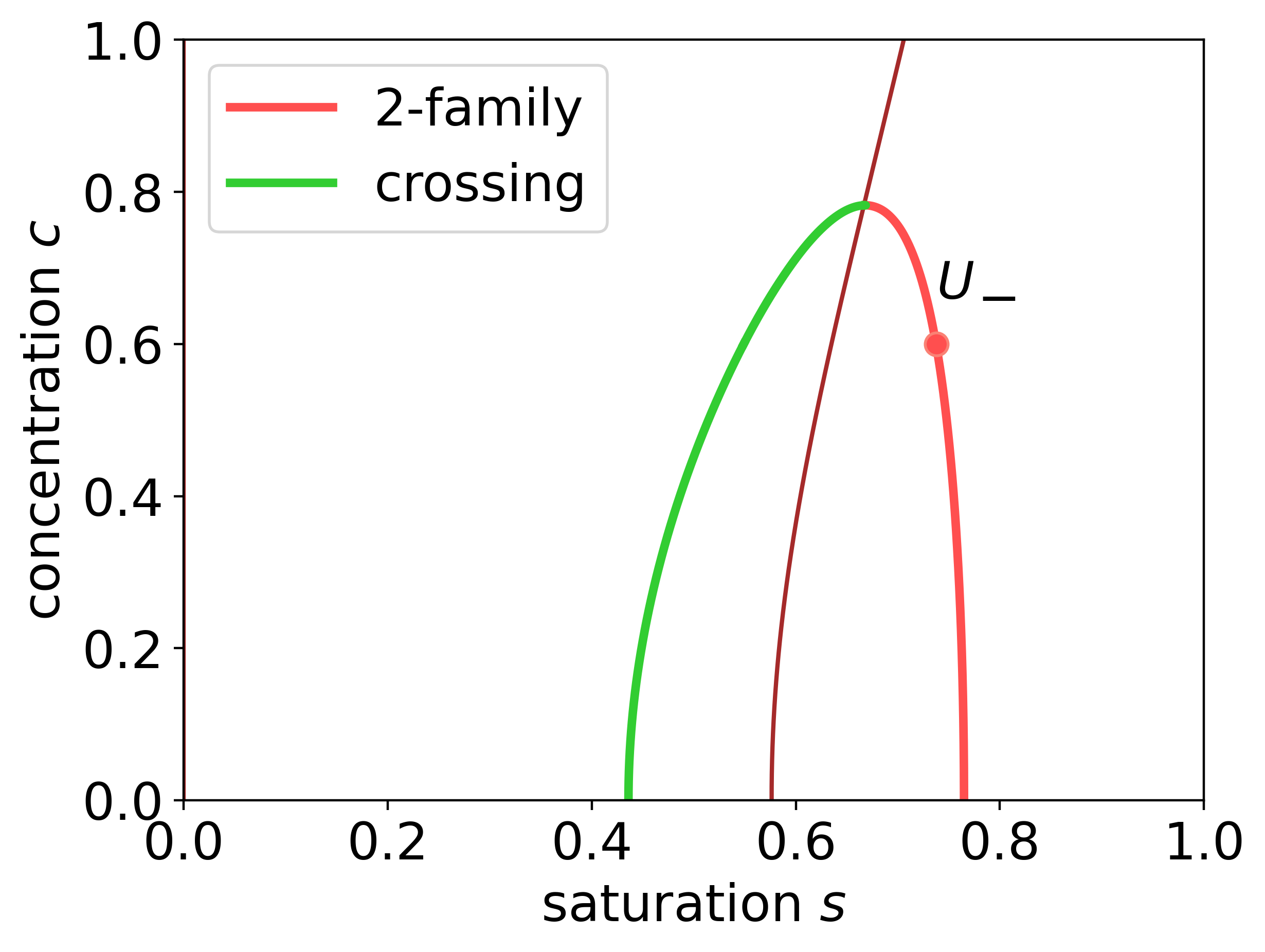} \\
    (a)\hfil
    \qquad\qquad(b)\hfil
    \caption{Contact discontinuities classified according to
    configuration of characteristics:
    (a) if $U_-\in\{\lambda^s>\lambda^c\}$,
    the blue (respectively, orange) curve corresponds to
    contacts with 1-family (resp., overcompressive) configuration;
    (b) if $U_-\in\{\lambda^s<\lambda^c\}$,
    the red (respectively, green) curve corresponds to
    contacts with 2-family (resp., crossing) configuration.}
    \label{fig:monotone-type-of-contact}
\end{figure}

For the polymer model without adsorption $M_0$,
a contact discontinuity with
\begin{itemize}
\item 1-family configuration has
$U_{\pm}\in \{\lambda^s>\lambda^c\}$;
\item 2-family configuration has
$U_{\pm}\in \{\lambda^s<\lambda^c\}$;
\item overcompressive configuration has
$U_{-}\in \{\lambda^s>\lambda^c\}$
and $U_{+}\in \{\lambda^s<\lambda^c\}$;
\item crossing configuration has
$U_{-}\in \{\lambda^s<\lambda^c\}$ and
$U_{+}\in \{\lambda^s>\lambda^c\}$.
\end{itemize}
See Fig.~\ref{fig:monotone-type-of-contact}.
Analogous statements are true for each type of $c$-shock wave
in the polymer model with adsorption~\eqref{eq:conslaw_adsorption}.

The main step in proving Theorem~\ref{thm:main-thm},
is the following Lemma~\ref{lm:main-lm}.
Theorem~\ref{thm:main-thm} easily follows from it
(for details, see the proof at the end of this section).

\begin{lemma}
\label{lm:main-lm}
For the polymer model~\eqref{eq:conslaw}
that satisfies the monotonicity condition~\eqref{eq:f-cond-monotone},
contact discontinuities with 1-family, 2-family, and overcompressive
configurations of characteristics
satisfy the vanishing adsorption admissibility criterion,
whereas contact discontinuities with crossing configuration do not.
\end{lemma}

\noindent
\emph{Proof of Lemma~\ref{lm:main-lm}}.
To prove this result,
we consider each type of configuration separately.
For 1-family, 2-family, and overcompressive configurations
we will give explicit sequences of solutions
of the system with adsorption~\eqref{eq:conslaw_adsorption}
that tend in $L^1_\text{loc}$ to the contact under consideration.
For crossing contact discontinuities,
we will show that there is no sequence of solutions
of Eq.~\eqref{eq:conslaw_adsorption} that tends in $L^1_\text{loc}$ to the contact.

\subsubsection*{1-family and 2-family contact discontinuities.}
Let us first prove that 1-family contact discontinuities are admissible under the adsorption criterion. An analogous argument works for the 2-family contact discontinuities. Fix $U_{\pm}=(s_{\pm},c_{\pm})\in\{\lambda^s>\lambda^c\}$.

\begin{figure}[ht]
    \centering
    \includegraphics[width=0.45\textwidth]{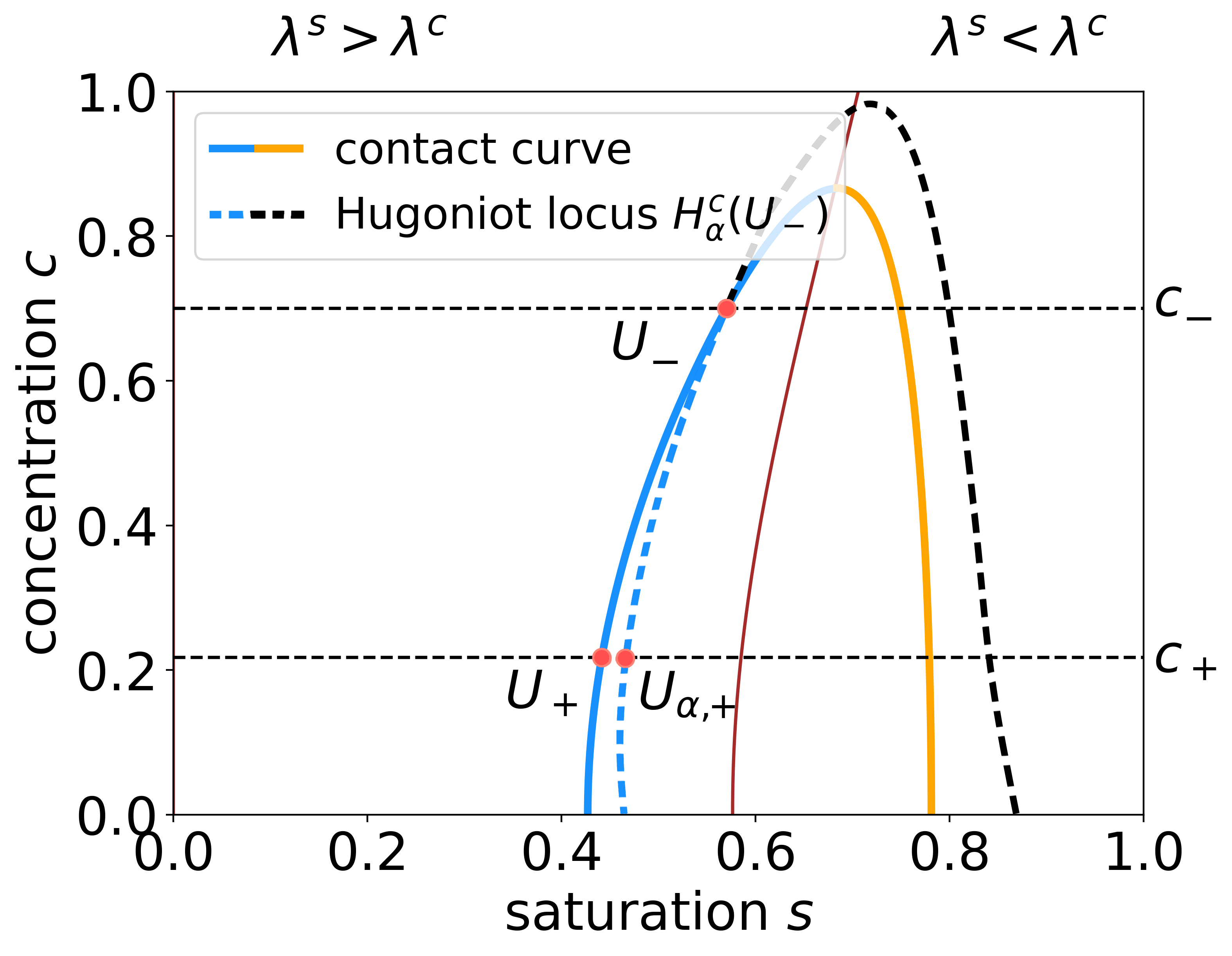}
    \quad
    \includegraphics[width=0.45\textwidth]{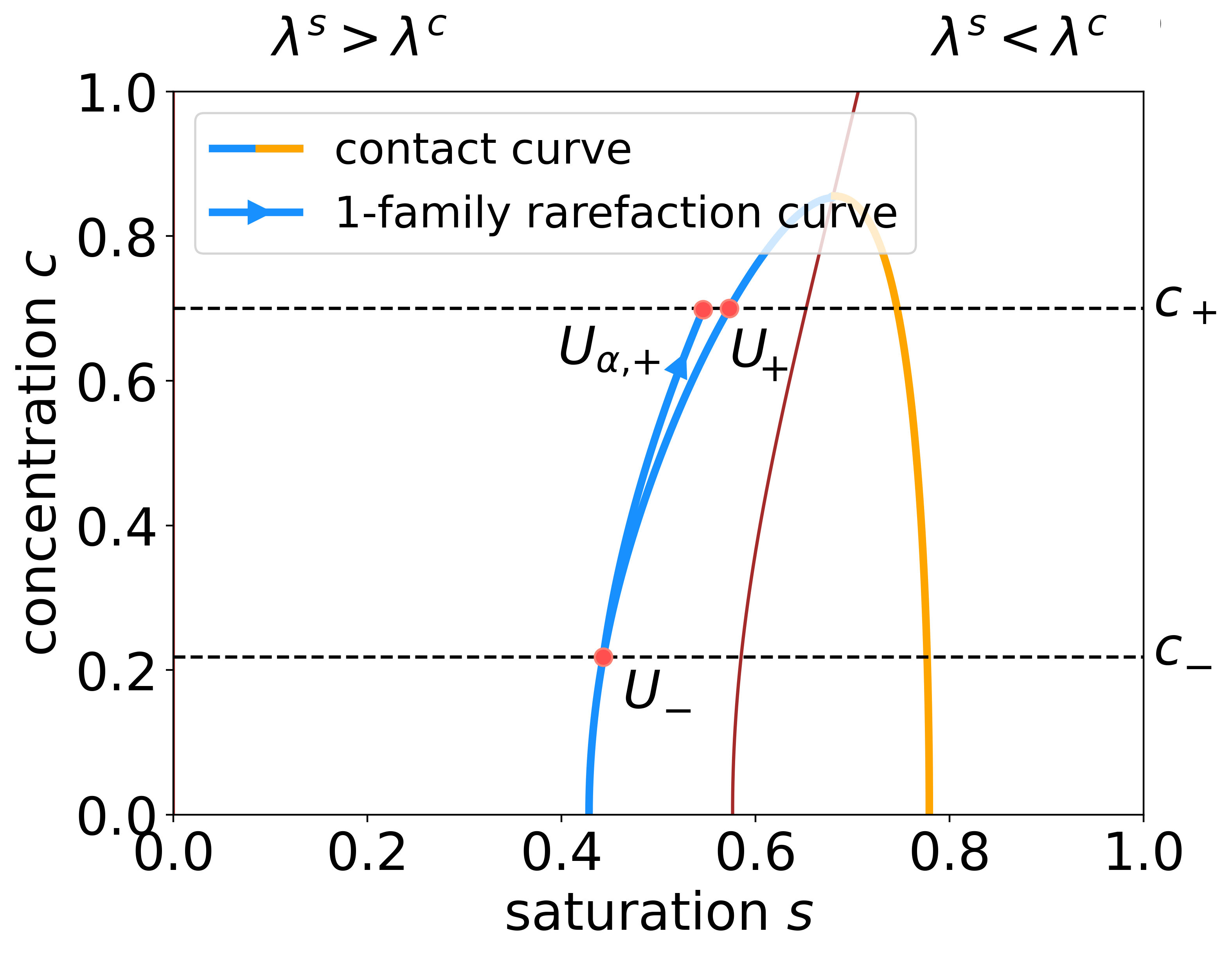}
    \\
    (\emph{i})\hfil
    \qquad\qquad\qquad\qquad(\emph{ii})\hfil
    \caption{Schematic illustration for the proof for the 1-family
    contact discontinuities.}
    \label{fig:i-ii}
\end{figure}

(\emph{i}) Let $c_->c_+$, as in Fig.~\ref{fig:i-ii}~(\emph{i}).
Then the contact between $U_-$ and~$U_+$
is the limit of 1-family $c$-shock waves
for the polymer model with adsorption~\eqref{eq:conslaw_adsorption}.
Indeed,
by Remark~\ref{rm:hugoniot-continuous},
the $c$-type branch of the Hugoniot locus $\mathcal{H}_\alpha^c(U_-)$
limits to the contact curve $\mathcal{H}_0^c(U_-)$ as $\alpha \to 0^+$.
Consequently,
we can choose
$U_{\alpha,+}=(s_{\alpha,+},c_{\alpha,+})
\in\mathcal{H}^c_\alpha(U_-)\cap\{\lambda^s>\lambda^c\}$
such that $U_{\alpha,+}$ tends to $U_+$ as $\alpha \to 0^+$.
Because both $U_-$ and $U_{\alpha,+}$
lie in the region $\{\lambda^s>\lambda^c\}$,
a shock wave between $U_-$ and $U_{\alpha,+}$
has the 1-family configuration of characteristics,
and thus is admissible for model~\eqref{eq:conslaw_adsorption}.
Moreover,
by the Rankine-Hugoniot condition~\eqref{eq:Rankine-Hugoniot},
the speed of this shock wave, $\sigma_\alpha$,
tends to the speed of a contact, $\sigma$, as $\alpha \to 0^+$.
From these observations we obtain the convergence
of the shock wave between $U_-$ and $U_{\alpha,+}$
to the contact between $U_-$ and $U_+$ in $L^1_{loc}$,
so that this contact is admissible by the vanishing adsorption criterion.

(\emph{ii}) Let $c_-<c_+$, as in Fig.~\ref{fig:i-ii}~(\emph{ii}).
We can approximate the contact between $U_-$ and~$U_+$
by $c$-rarefaction waves for the polymer model
with adsorption~\eqref{eq:conslaw_adsorption}.
Indeed, according to Remark~\ref{rm:rare-continuous},
for small enough $\alpha>0$,
the $c$-type integral curve $\mathcal{I}_{\alpha}(U_-)$
is close to $c$-type integral curve  $\mathcal{I}_0(U_-)$.
As a consequence, we can choose
$U_{\alpha,+}\in\mathcal{I}_{\alpha}(U_-)\cap\{\lambda^s>\lambda^c\}$
such that $U_{\alpha,+}$ tends to $U_+$ as $\alpha \to 0^+$.
Moreover, the characteristic speeds
$\lambda^c_\alpha(U_-)$ and $\lambda^c_\alpha(U_{\alpha,+})$
approach each other as $\alpha \to 0^+$ because of Eq.~\eqref{eq:dlambda_c-ads}:
    \begin{align*}
        |\lambda^c_\alpha(U_{\alpha,+})-\lambda^c_\alpha(U_-)|
	\leq \mathrm{const}_1 \cdot |(D\lambda^c_\alpha)\,r^c_\alpha|
	\leq \mathrm{const}_2 \cdot \alpha \to 0^+,
    \end{align*}
where the constants $\mathrm{const}_1$ and $\mathrm{const}_2$
are independent of $\alpha$.
The characteristic speed $\lambda^c_\alpha(U_-)$
also tends to the contact speed $\sigma=\lambda^c_0(U_-)$
because of Eq.~\eqref{eq:lambda}.
We deduce the convergence of the $c$-rarefaction wave
between $U_-$ and $U_{\alpha,+}$ to the contact
between $U_-$ and $U_+$ in $L^1_\text{loc}$,
and thus conclude that such a contact is admissible
by the vanishing adsorption criterion.

\subsubsection*{Overcompressive contact discontinuities. }
Let us prove that overcompressive contact discontinuities are admissible under the vanishing adsorption criterion. Fix $U_-\in\{\lambda^s>\lambda^c\}$ and $U_+\in \mathcal{H}_0(U_-)\cap \{\lambda^s<\lambda^c\}$. The interesting feature of the overcompressive contact is that it can be represented as a sequence of two shock waves  with equal speeds: $s$-wave between $U_-$ and $U_1$, and $c$-wave between $U_1$ and $U_+$. Here $U_1=(s_1,c_-)\in\{\lambda^s<\lambda^c\}$ is the uniquely defined state on the contact curve for $U_-$ and the speed of the shocks coincides with the characteristic speed $\lambda^c_0(U_-)$. The properties we assume on $f$ (monotone  and  $S$-shape in $s$, $f(0,c)=f'(0,c)=0$) guarantee that the $s$-wave between $U_-$ and $U_1$ is a shock wave.  Notice that the $c$-wave between $U_1$ and $U_+$ is a 2-family contact (see Def.~\ref{def:contact-char-configs}).

\begin{figure}[ht]
    \centering
    \includegraphics[width=0.45\textwidth]{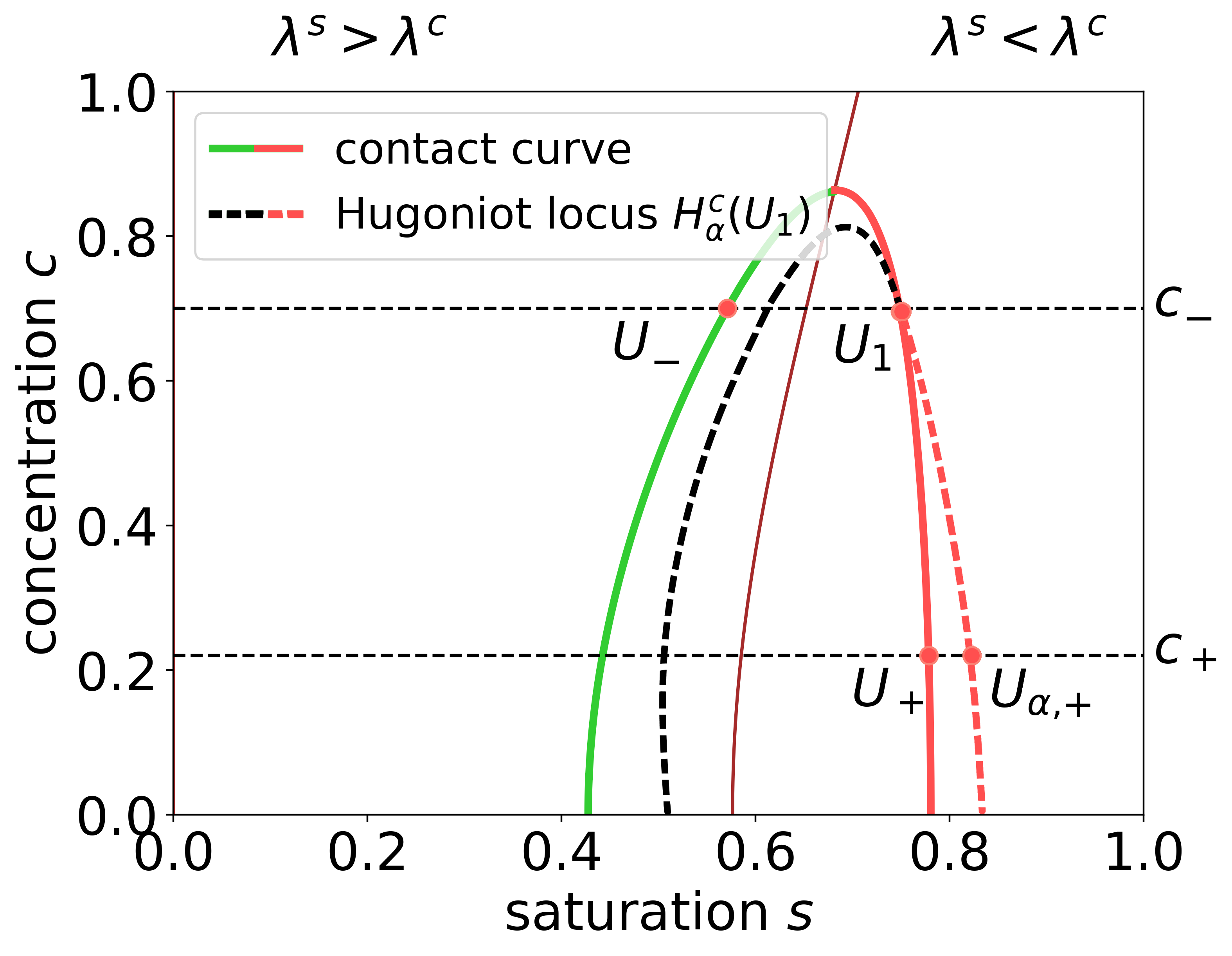}
    \quad
    \includegraphics[width=0.45\textwidth]{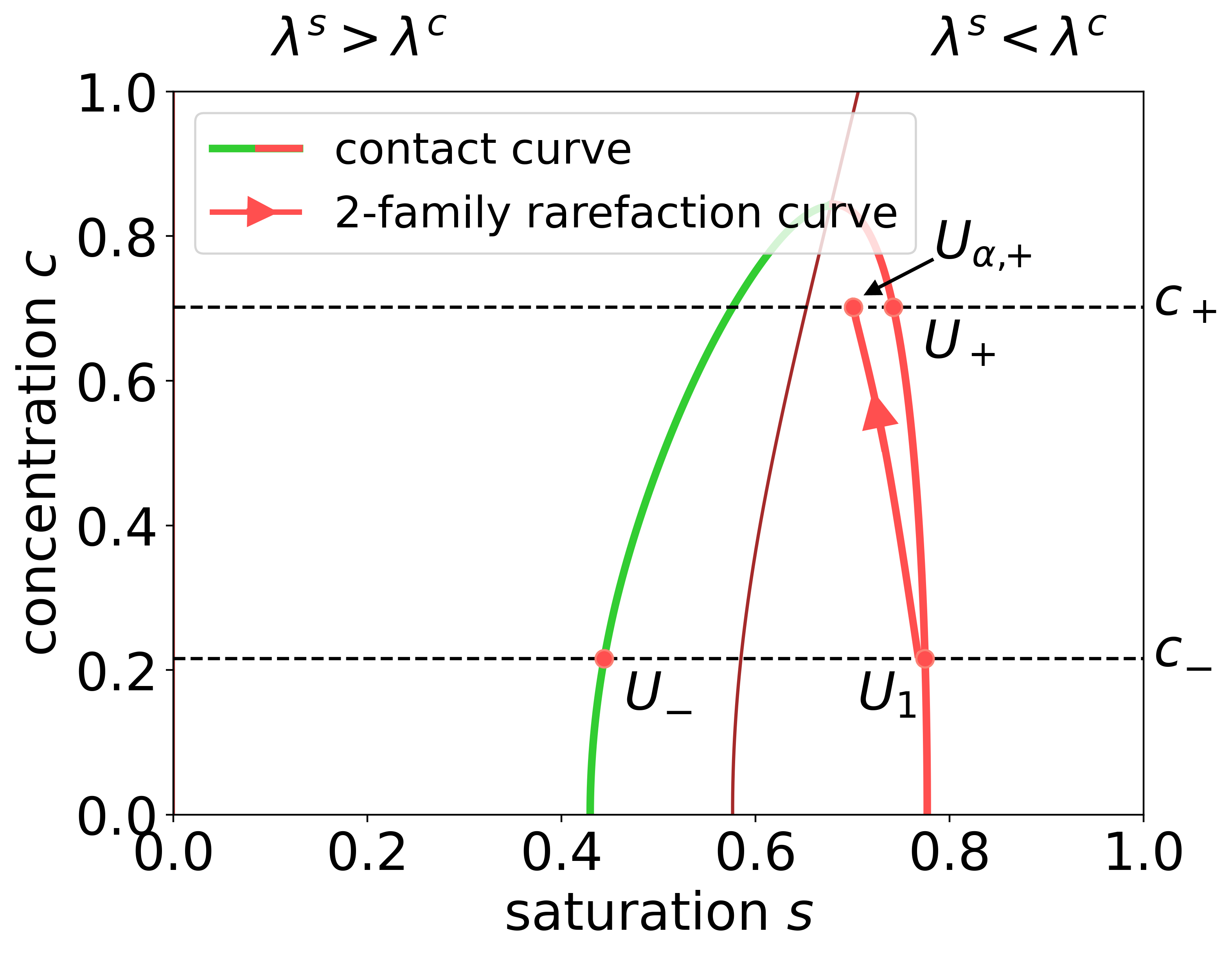}
    \\
    (\emph{iii})\hfil
    \qquad\qquad\qquad\qquad(\emph{iv})\hfil
    \caption{Schematic illustration for the proof for overcompressive
    contact discontinuities.}
    \label{fig:iii-iv}
\end{figure}

(\emph{iii}) Let $c_->c_+$, see Fig.~\ref{fig:iii-iv}~(\emph{iii}). By case (\emph{i}), the 2-family contact between $U_1$ and $U_+$ can be approximated by a sequence of $c$-shock waves between $U_1$ and $U_{\alpha,+}$. As a result we get that the sequence of the two waves: $s$-wave between $U_-$ and $U_1$, and
$c$-shock wave between $U_1$ and $U_{\alpha,+}$ tends in $L^1_\text{loc}$ to a contact connecting states $U_-$ and~$U_+$.

(\emph{iv}) Let $c_-<c_+$, see Fig.~\ref{fig:iii-iv}~(\emph{iv}). By case (\emph{ii}), the 2-family contact between $U_1$ and $U_+$ can be approximated by a sequence of $c$-rarefaction waves between $U_1$ and $U_{\alpha,+}$. As a result we get that the sequence of the two waves: $s$-wave between $U_-$ and $U_1$, and
$c$-rarefaction wave between $U_1$ and $U_{\alpha,+}$ tends in $L^1_\text{loc}$ to a contact connecting states $U_-$ and $U_+$.

\subsubsection*{Crossing contact discontinuities.}
Let us prove that crossing contact discontinuities
are not admissible under the vanishing adsorption criterion.
Let us fix $U_-\in\{\lambda^s<\lambda^c\}$
along with $U_+\in \mathcal{H}_0(U_-)\cap \{\lambda^s>\lambda^c\}$.

\begin{figure}[ht]
    \centering
    \includegraphics[width=0.65\textwidth]{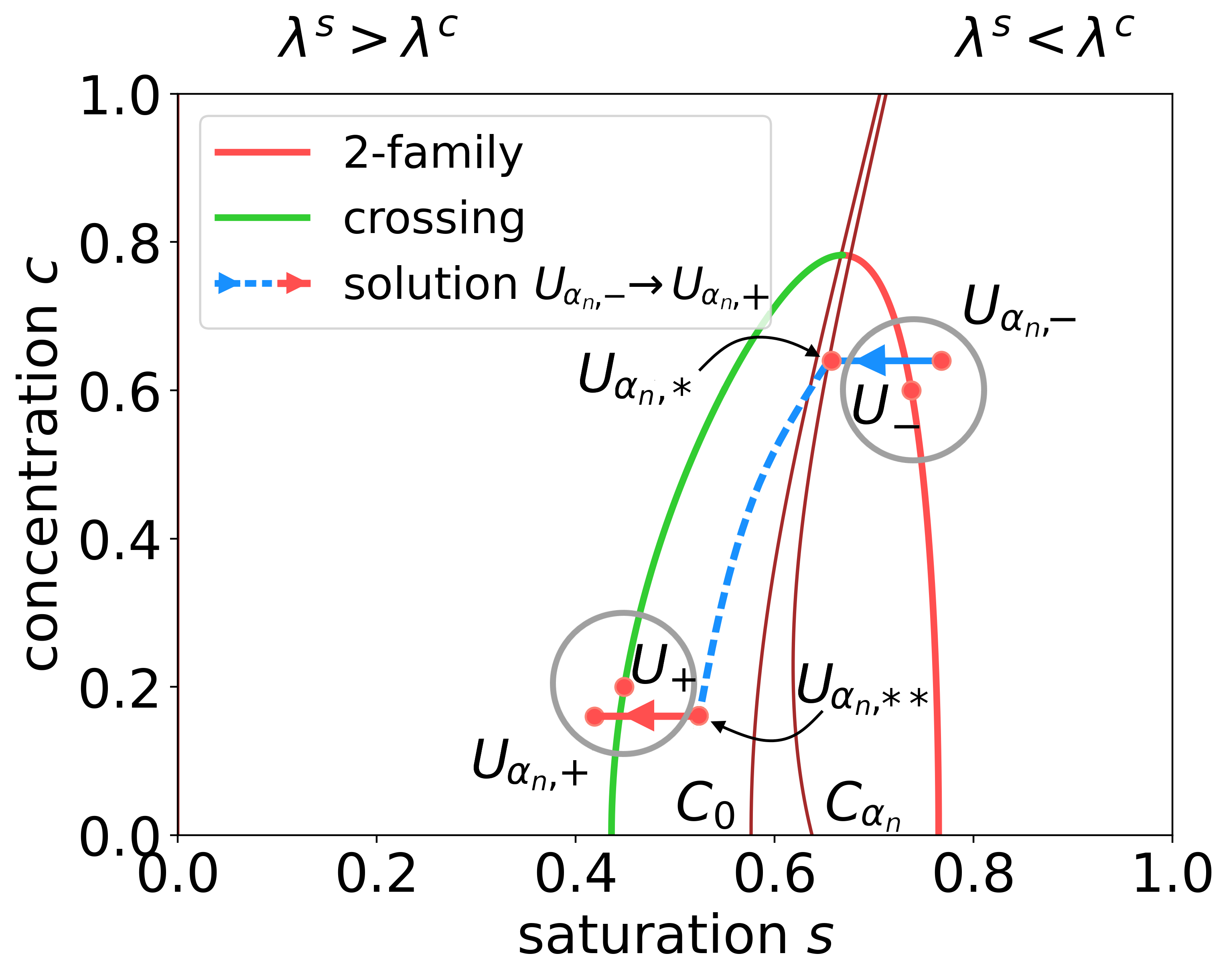}
    \\
    \qquad\textit{(v)}\hfil
    \caption{Schematic illustration for the proof that
    crossing contact discontinuities are not admissible
    under the vanishing adsorption criterion.}
    \label{fig:c+smallerc-crossing}
\end{figure}

(\emph{v}) Assume that $c_->c_+$, as in Fig.~\ref{fig:c+smallerc-crossing}.
(The case $c_-<c_+$ can be proved in a similar fashion.)
Let $U_*=(s_{*},c_-)$ be the unique point
on the interior coincidence locus $\mathcal{C}_0$ for $c=c_-$.
Define
\begin{align*}
    \sigma_*: =\frac{f(U_*)}{s_*}=f_s(U_*),\qquad
    \delta: = |f_s(U_*) - f_s(U_-)|>0.
\end{align*}
The proof is by contradiction.
Suppose that there exists a sequence of solutions
to a Riemann problem for the polymer model~\eqref{eq:conslaw_adsorption}
with left state $U_{\alpha_n,-}=(s_{\alpha_n,-},c_{\alpha_n,-})$
and right state $U_{\alpha_n,+}=(s_{\alpha_n,+},c_{\alpha,+})$
that tends in $L^1_\text{loc}$ to a contact
that connects states $U_-$ and $U_+$ as $\alpha_n \to 0^+$ and $n\in\mathbb{N}$.
Then necessarily $U_{\alpha_n,-}\to U_-$ and $U_{\alpha_n,+}\to U_+$
as $\alpha_n \to 0^+$.
Fix $\varepsilon\in(0,(c_--c_+)/2)$;
then there exists $\bar{\alpha}>0$ such that,
for all $\alpha_n\in(0,\bar{\alpha})$,
the state $U_{\alpha_n,-}$ lies in the $\varepsilon$-neighborhood of $U_-$
and the state $U_{\alpha_n,+}$ lies in the $\varepsilon$-neighborhood of $U_+$.
We take small enough $\varepsilon>0$ such that
these $\varepsilon$-neighborhoods
do not intersect the interior coincidence locus $\mathcal{C}_{\alpha_n}$
if $\alpha_n\in[0,\bar{\alpha})$.

From Sec.~\ref{subsubsec:RP-polymer} we obtain that any solution to a Riemann problem with left state $U_{\alpha_n,-}$ and right state $U_{\alpha_n,+}$  has the following structure (see Fig.~\ref{fig:c+smallerc-crossing}):
\begin{align*}
    U_{\alpha_n,-}\xrightarrow{s} U_{\alpha_n,*}\xrightarrow{c}U_{\alpha_n,**}\xrightarrow{s} U_{\alpha_n,+}.
\end{align*}
Here $U_{\alpha_n,*}=(s_{\alpha_n,*},c_{\alpha_n,-})$
is a state on the interior coincidence locus~$\mathcal{C}_{\alpha_n}$
for $c=c_{\alpha_n,-}$,
and $U_{\alpha_n,**}=(s_{\alpha_n,**},c_{\alpha_n,+})\in\{\lambda^s>\lambda^c\}$
is a state on the $c$-branch of the Hugoniot locus
$\mathcal{H}^c_{\alpha_n}(U_{\alpha_n,*})$ for $c=c_{\alpha_n,+}$.
Notice that the wave $U_{\alpha_n,-}\xrightarrow{s} U_{\alpha_n,*}$
is an $s$-rarefaction wave
because the flow function $f(s,c_{\alpha_n,-})$ is concave for $s>s_{\alpha_n}$.

To arrive at a contradiction,
it is enough to show that the characteristic speeds at the end states
$U_{\alpha_n,-}$ and $U_{\alpha_n,*}$ of the $s$-rarefaction wave
do not approach each other as $\alpha_n \to 0^+$.
To this end,
notice that the characteristic speed $f_s(U_{\alpha_n,*})$
coincides with the speed of a $c$-shock wave
$\sigma_{\alpha_n,*}:=\frac{f(U_{\alpha_n,*})}{s_{\alpha_n,*}}$.
For $\alpha_n$ small enough we have
\begin{align*}
    \lvert f_s(U_{\alpha_n,-}) - f_s(U_-)\rvert<\delta/4
    \text{ and }
    \lvert f_s(U_{\alpha_n,*})-f_s(U_*)\rvert
    = \lvert \sigma_{\alpha_n,*}-\sigma_{*}\rvert<\delta/4.
\end{align*}
Thus
\begin{align*}
    \lvert f_s(U_{\alpha_n,-})&-f_s(U_{\alpha_n,*})\rvert \\
    &\ge \lvert f_s(U_{-})-f_s(U_{*})\rvert
    - \lvert f_s(U_{\alpha_n,-})-f_s(U_{-})\rvert
    - \lvert f_s(U_{\alpha_n,*})-f_s(U_{*})\rvert \\
    &> \delta/2.
\end{align*}
The constant $\delta$ does not depend on $\alpha_n$,
so we find a contradiction.
\emph{Lemma~\ref{lm:main-lm} is proved.}

To complete the proof of Theorem~\ref{thm:main-thm},
we observe that a contact for which one of the states $U_-$ or $U_+$
lies on the interior coincidence locus $\mathcal{C}_0$
is also admissible by the vanishing adsorption admissibility criterion.

\section{Adsorption admissibility of undercompressive contact discontinuities}
\label{sec:example}
In Sec.~\ref{sec:proof},
we demonstrated the equivalence of four admissibility criteria
for contact discontinuities,
summarized in Defs.~\ref{def:lax}--\ref{def:AD},
assuming that the model obeys the monotonicity condition~\eqref{eq:f-cond-monotone}.
We now present an example of a non-monotone model
with undercompressive contact discontinuities
that are admissible under the de~Souza-Marchesin and the adsorption criteria
but not under the Keyfitz-Kranzer/Isaacson or Isaacson-Temple criteria.
Therefore, if the monotonicity condition is violated,
these criteria can be distinct.

\begin{figure}[ht]
    \centering
     \includegraphics[width=0.43\textwidth]{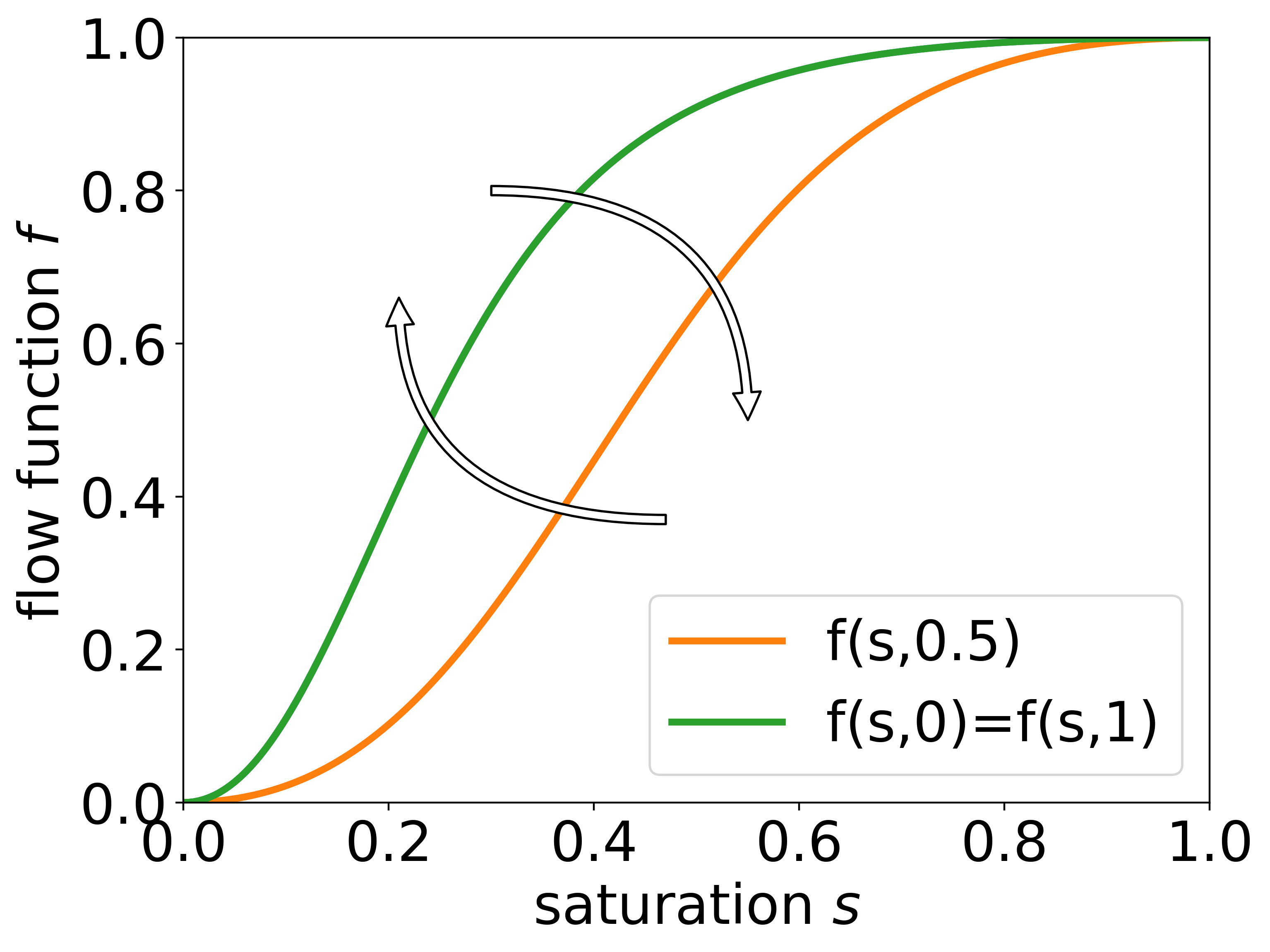}
     \includegraphics[width=0.44\textwidth]{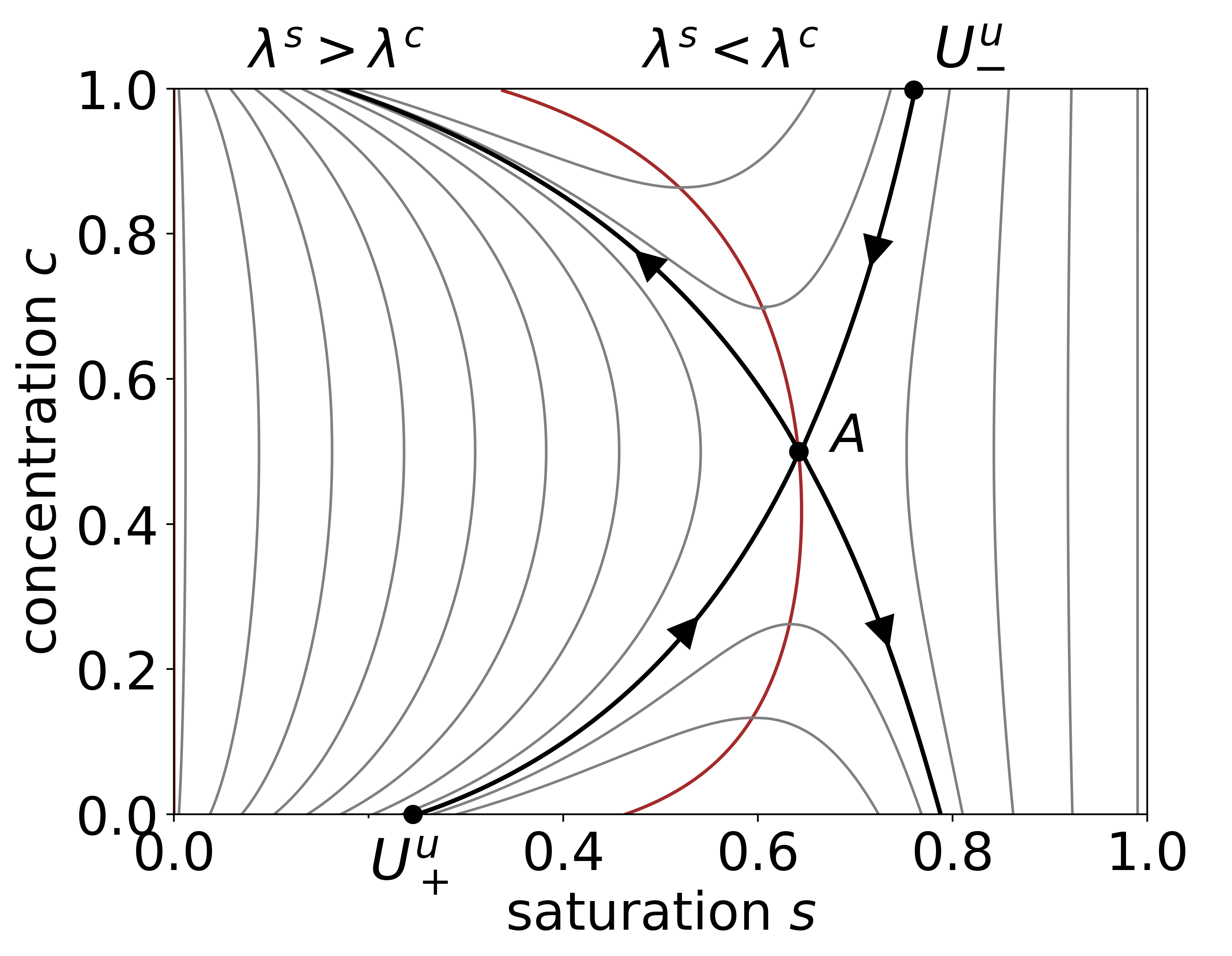}
       \\
    (a)\hfil
    \qquad\qquad\qquad(b)\hfil
    \caption{(a) Fractional flow function $f(s,c)$
    defined by formula~\eqref{eq:f-nonmonotone}. \\
    (b) Contact curves for the boomerang model.
    }
    \label{fig:f_nonmonotone}
\end{figure}

Consider a chemical flooding model~\eqref{eq:conslaw}
with a simple non-monotone flow function,
namely the ``boomerang'' model of Ref.~\cite{nonmonotonicity2021}.
For each $s \in (0,1)$,
the fractional flow $f(s, c)$
decreases as $c$ ranges from $0$ to some value $c^M\in(0,1)$,
and it increases as $c$ ranges from $c^M$ to $1$,
returning to the same value $f(s,1) = f(s,0)$.
An example, shown in Fig.~\ref{fig:f_nonmonotone}(a), is
\begin{align}
\label{eq:f-nonmonotone}
\text{$f(s,c) := \frac{s^2}{s^2+\mu(c)\,(1-s)^2}$
with $\mu(c):=0.1+1.8\,c\,(1-c)$,}
\end{align}
for which $c^M=0.5$.

The contact curves for system~\eqref{eq:conslaw}
are the level sets of the eigenvalue $\lambda^c_0 = f/s$;
equivalently, they are orbits of the ODE
$\dot U = r^c_0(U)$, where $r^c_0$ is given by Eq.~\eqref{eq:eigenvectors}.
This equation has an equilibrium when $f_c = 0$ and $\lambda^s = \lambda^c_0$,
such as at the point $A$ on $\mathcal{C}_0$ where $c = c^M$.
The nature of such an equilibrium is determined by the derivative $Dr^c_0$
evaluated at the equibrium.
A straightforward calculation finds that its determinant
equals $-(f_{sc})^2 + f_{ss}\,f_{cc}$.
For the boomerang model,
this quantity is negative,
so that its equilibrium is a saddle-point.

The contact curves for the boomerang model
are depicted in Fig.~\ref{fig:f_nonmonotone}(b).
The structure of solutions to local Riemann problems
near this equilibrium was considered by de~Souza and Marchesin~\cite{cido-dan1998}.
The authors proved that there exists an undercompressive contact
\emph{i.e.},\ a contact with crossing characteristic configuration
that is admissible under de~Souza-Marchesin admissibility criterion.
In this section we want to show that this undercompressive contact
is also admissible under the adsorption criterion
(Theorem~\ref{thm:undercompressive-admissible}).

Fix two points in the stable manifold of the equilibrium $A$,
\begin{align}
\text{$U_-^u=(s_-^u,1)\in\{\lambda^s<\lambda^c\}$
and $U_+^u=(s_+^u,0)\in\{\lambda^s>\lambda^c\}$},
\end{align}
as in Fig.~\ref{fig:f_nonmonotone}(b).
The corresponding speed $\sigma^u$ can be calculated
from the Rankine-Hugoniot condition~\eqref{eq:c-type-hugoniot}.
As the contact curve connecting points $U_-^u$ and $U_+^u$
is continuous and monotone in $c$,
it represents contacts that are admissible
by the de~Souza-Marchesin admissibility criterion;
see Def.~\ref{def:KKIT}.
However, it fails to satisfy the Keyfitz-Kranzer/Isaacson
and Isaacson-Temple admissibility criteria
because $U_+^u$ and $U_-^u$ lie on opposite sides
of the interior coincidence locus.

To apply the adsorption criterion,
one needs to adopt an admissibility criterion
for weak solutions of the system~\eqref{eq:conslaw_adsorption}.
We will use the following viscosity criterion
of Ref.~\cite{nonmonotonicity2021}.

\begin{definition}[Bakharev \emph{et~al.}~\cite{nonmonotonicity2021}]
\label{def:vanishingviscosity}
Fix $\kappa>0$.
We say that a $c$-shock is admissible
if it is the limit of traveling wave solutions of the system
\begin{equation}
\label{eq:main_system_smooth_cap_diff}
\begin{cases}
s_t + f(s, c)_x = \varepsilon_c\,s_{xx}, \\
\left[cs + \alpha\,a(c)\right]_t + \left[cf(s,c)\right]_x
= \varepsilon_c\,(c s_x)_x + \varepsilon_d\,c_{xx}
\end{cases}
\end{equation}
as $\varepsilon_{c} \to 0^+$ and $\varepsilon_d \to 0^+$
with fixed ratio $\varepsilon_d/\varepsilon_c=:\kappa$.
\end{definition}

\begin{remark}
If we do not assume $\varepsilon_d/\varepsilon_c$ to be a fixed constant,
then, as shown in Ref.~\cite{nonmonotonicity2021},
there are multiple solutions to some Riemann problems,
depending on the value $\kappa$.
Fixing the value of the parameter~$\kappa$ resolves this problem.
Moreover, as one can ascertain from the proof of
Theorem~\ref{thm:undercompressive-admissible},
the precise value of $\kappa$ does not affect its conclusion.
\end{remark}

Now we are ready to prove Theorem~\ref{thm:undercompressive-admissible}.

\subsection{Proof of Theorem~\ref{thm:undercompressive-admissible}}
The main idea of the proof is to show that the undercompressive contact can be
approximated by a sequence of undercompressive $c$-shock waves for the
system~\eqref{eq:conslaw_adsorption} satisfying the  viscosity criterion. The
existence of undercompressive $c$-shock waves is shown
in Ref.~\cite{nonmonotonicity2021}:

\begin{proposition}[Bakharev \emph{et~al.}~\cite{nonmonotonicity2021}]
    \label{Theorem-nonmonotonicity}
Consider the family conservation laws~\eqref{eq:conslaw_adsorption},
parameterized by $\alpha>0$,
with flux function $f$ defined in Eq.~\eqref{eq:f-nonmonotone}.
There exist functions $\sigma_{\min}(\alpha)$ and $\sigma_{\max}(\alpha)$
of $\alpha$ with $0<\sigma_{\min}(\alpha)<\sigma_{\max}(\alpha)<\infty$
such that, for every $\kappa=\varepsilon_d/\varepsilon_c\in(0, +\infty)$,
there exist functions
\begin{itemize}
    \item $s_-^{\kappa}(\alpha)\in[0,1]$ and $s_+^{\kappa}(\alpha)\in[0,1]$ and
    \item $\sigma^{\kappa}(\alpha)
    \in[\sigma_{\min}(\alpha),\sigma_{\max}(\alpha)]$
\end{itemize}
of $\alpha$ such that the undercompressive $c$-shock wave
from $U_-^{\kappa}(\alpha):=(s_-^{\kappa}(\alpha),1)$
to $U_+^{\kappa}(\alpha):=(s_+^{\kappa}(\alpha),0)$
with speed $\sigma^{\kappa}(\alpha)$
is admissible by the viscosity criterion of Def.~\ref{def:vanishingviscosity}.
\end{proposition}

\begin{figure}[ht]
    \centering
    \includegraphics[width=0.565\textwidth]{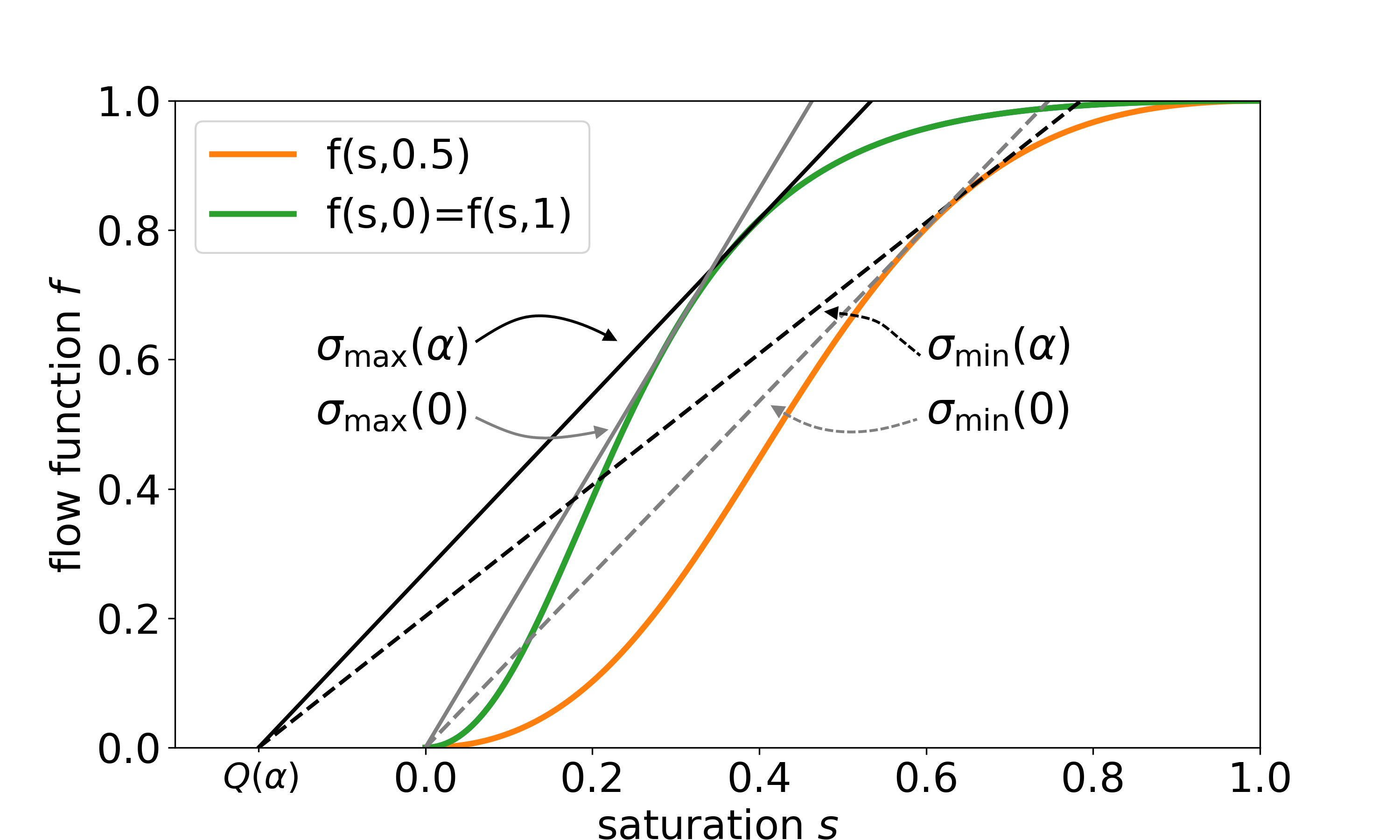}
    \includegraphics[width=0.425\textwidth]{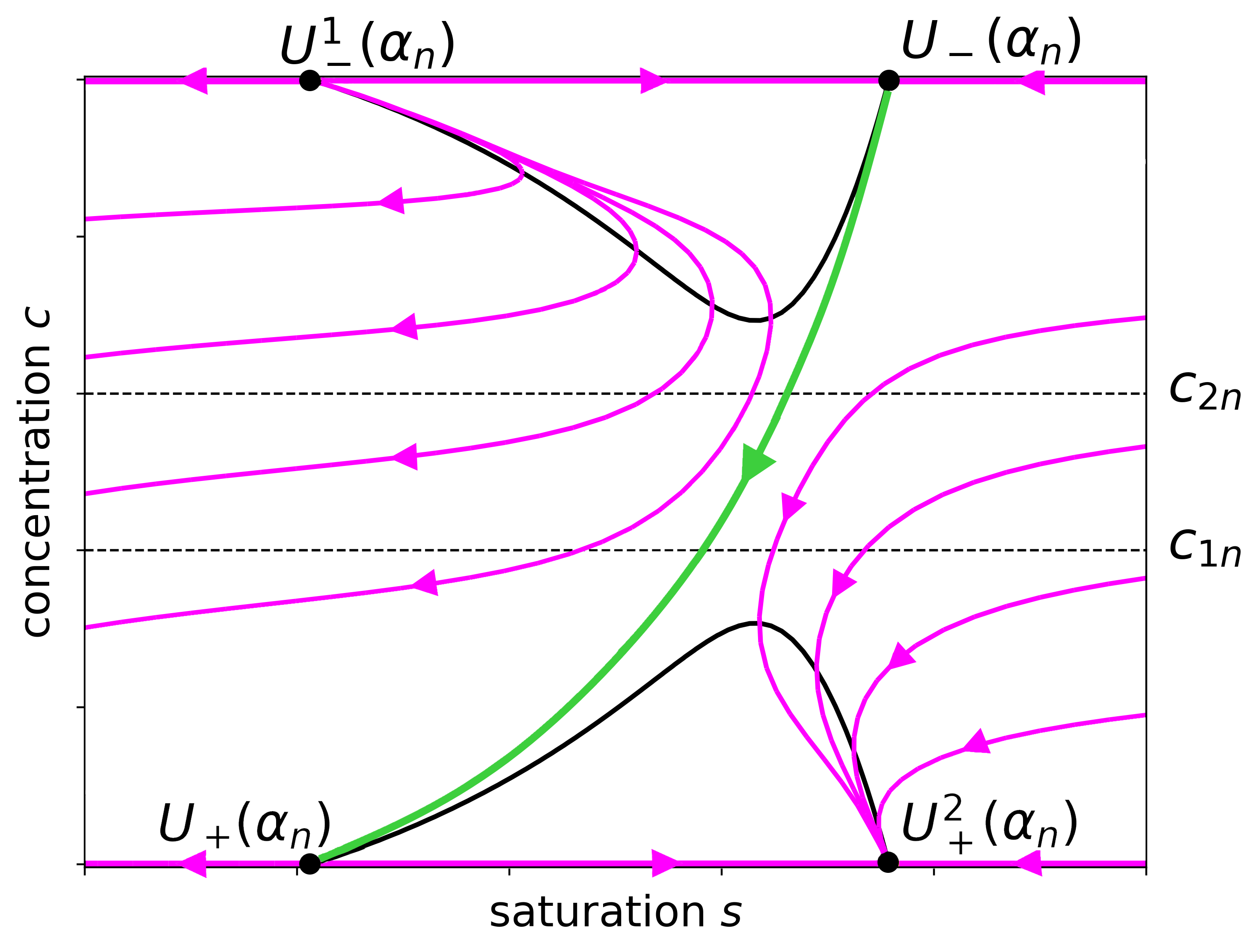}\\
    \qquad(a)\hfil\qquad\qquad\qquad\qquad\qquad(b)\hfil
    \caption{(a)~The geometric determination of $\sigma_{\min}(\alpha)$
    and $\sigma_{\max}(\alpha)$ for $\alpha\geq0$.
    The green curve is the graph of $f(\cdot,0)=f(\cdot,1)$;
    the orange curve is the graph of $f(\cdot,0.5)$.
    (b)~Phase portrait for the dynamical system~\eqref{eq:dyn_sys_cap_diff}.
    The green curve is the heteroclinic orbit joining
    the two saddle points $U_-(\alpha_n)$ and $U_+(\alpha_n)$,
    whereas the purple curves are other representative orbits.
    The black curves are where $s_\xi=0$.}
    \label{fig:sigma_max_sigma_min}
\end{figure}

\begin{remark}
\label{rmk:sigmamin-max}
Using the results of Ref.~\cite[Section 4.2.2]{nonmonotonicity2021},
one obtains the explicit expressions for $\sigma_{\min}(\alpha)$
and $\sigma_{\max}(\alpha)$;
see Fig.~\ref{fig:sigma_max_sigma_min}~(a) for the geometric representation.
In particular,
\begin{itemize}
    \item $\sigma_{\min}(\alpha)$ equals the minimum, with respect to $c\in[0,1]$,
    of the slope of the tangent to the graph of $f(\cdot,c)$
    that starts at $Q(\alpha):=(-\alpha\,a_1(1),0)$.
    By formula~\eqref{eq:f-nonmonotone},
    this value is obtained for $c = 0.5$,
    \emph{i.e.},\ for the orange curve $f(\cdot,0.5)$.
    Notice that
    \begin{align}
    \text{$\sigma_{\min}(\alpha)\to\sigma^u$ as $\alpha\to0^+$}.
    \end{align}

    \item $\sigma_{\max}(\alpha)$ equals the slope of the tangent to the graph
    of $f(\cdot,0)=f(\cdot,1)$, drawn in green, that starts at $Q(\alpha)$.
    As the function $f(\cdot,0)$ is S-shaped
    for all $\alpha$ in some interval $[0,\bar{\alpha}]$,
    we find that
    \begin{align*}
    \sigma_{\max}:=\limsup\limits_{\alpha\in[0,\bar{\alpha}]}
    \sigma_{\max}(\alpha)<\infty.
    \end{align*}
\end{itemize}
\end{remark}

Let us fix $\kappa>0$.
In what follows we suppress the superscript $\kappa$,
writing $\sigma(\alpha)$ instead of $\sigma^{\kappa}(\alpha)$, \emph{etc}.

To prove Theorem~\ref{thm:undercompressive-admissible},
it suffices to prove the next Lemma.
\begin{lemma}
\label{lm:sigman-limit}
\text{$\lim\limits_{\alpha \to 0^+}\sigma(\alpha) = \sigma^u$,
$\lim\limits_{\alpha \to 0^+} U_-(\alpha) = U_-^u$,
and $\lim\limits_{\alpha \to 0^+} U_+(\alpha)= U_+^u$.}
\end{lemma}
Indeed,
Lemma~\ref{lm:sigman-limit} implies that
the undercompressive $c$-shock waves connecting states
$U_-(\alpha)$ and $U_+(\alpha)$
converge to the undercompressive contact connecting states
$U_-^u$ and $U_+^u$ in $L^1_\text{loc}$ as $\alpha\to0^+$.

\par\noindent
\emph{Proof of Lemma~\ref{lm:sigman-limit}.}
First, let us prove that $\sigma(\alpha)\to \sigma^u$ as $\alpha\to0^+$.
We will show that $\sigma(\alpha)$ has a unique accumulation point,
namely $\sigma^u$.
Choose a positive sequence $\{\alpha_n\}$ for $n\in\mathbb{N}$
such $\alpha_n\to0^+$ and $\sigma_n:=\sigma(\alpha_n)$
converges to some limit~$\sigma_\infty$.
By Remark~\ref{rmk:sigmamin-max},
we have that $\sigma_\infty\in[\sigma^u,\sigma_{\max}]$.

Let us prove that $\sigma_\infty=\sigma^u$ by contradiction.
Assume, to the contrary, that $\sigma_\infty>\sigma^u$.
Consider a traveling wave solution of
system~\eqref{eq:main_system_smooth_cap_diff}
that connects $U_-(\alpha_n)$ and $U_+(\alpha_n)$
(which exists by Proposition~\ref{Theorem-nonmonotonicity}).
This solution has the form $s(x,t) := s(\xi)$ and $c(x,t) := c(\xi)$,
where $\xi := \left[x - \sigma(\alpha_n)\,t\right]/\varepsilon_c$;
it satisfies $s(\pm\infty) = s_{\pm}(\alpha_n)$,
$c(-\infty) = 1$, and $c(+\infty) = 0$.
Substituting into~\eqref{eq:main_system_smooth_cap_diff}
and simplifying in the standard way
(see \emph{e.g.}~\cite[Sec.~4.1]{nonmonotonicity2021}),
we obtain the following system of ODEs:
\begin{equation}
\label{eq:dyn_sys_cap_diff}
\begin{cases}
s_\xi = f(s, c) - \sigma_n (s + \alpha_n\,a_1(1)), \\
c_\xi = (\sigma_n\alpha_n/\kappa) \, \left[a_1(1)\,c  - a_1(c)\right].
\end{cases}
\end{equation}
The limit of the system~\eqref{eq:dyn_sys_cap_diff} as $\alpha_n\to0^+$
is the system
\begin{equation}
\label{eq:dyn_sys_cap_diff_alpha=0}
\begin{cases}
s_\xi = f(s, c) - \sigma_\infty s, \\
c_\xi = 0.
\end{cases}
\end{equation}

For all $\sigma_n \in (\sigma_{\min}(\alpha_n),\sigma_{\max}(\alpha_n))$,
the phase portrait has the structure
illustrated in Fig.~\ref{fig:sigma_max_sigma_min}(b).
(See Ref.~\cite[Sec.~4.3]{nonmonotonicity2021}
for a detailed analysis of properties of trajectories.)
In particular,
there exists the saddle-to-saddle connection,
drawn in green,
between the equilibria $U_-$ and $U_+$;
it can be parameterized by a function $s=s_\text{orbit}(c)$.

Following Ref.~\cite{nonmonotonicity2021},
we examine the saturation nullcline of this system,
\emph{i.e.},\ the zero-set of the Eq.~\eqref{eq:dyn_sys_cap_diff},
and its concentation counterpart.
The saturation nullcline $\{(s,c): s_\xi = 0\}$ is drawn as black curves.
Because $a_1$ is concave,
the concentration nullcline $\{(s,c): c_\xi = 0\}$ consists of two lines:
\[
\{(s,c): a_1(1)\,c - a_1(c) = 0\}
= [0,1] \times \{ 0 \} \cup [0,1] \times \{ 1 \}.
\]

In the terminology of Ref.~\cite[Sec.~4.2.2]{nonmonotonicity2021},
the saturation nullcline of system~\eqref{eq:dyn_sys_cap_diff} has Type~II.
In particular,
it satisfies the following properties,
depicted in Fig.~\ref{fig:sigma_max_sigma_min}(b):
\begin{enumerate}
    \item Because $f(\cdot,c)$ is S-shaped for each $c$,
    the saturation nullcline contains at most
    two points for each fixed concentration $c$.
    For $c=1$, there are exactly two critical points,
    the repeller $U_-^1(\alpha_n)$ and the saddle point $U_-(\alpha_n)$;
    and for $c=0$, there are exactly two critical points,
    the saddle point $U_+(\alpha_n)$ and the attractor $U_+^2(\alpha_n)$.

    \item Because $f_c$ changes sign at most once for each $s$
    (see formula~\eqref{eq:f-nonmonotone}),
    the saturation nullcline contains at most two points for each fixed $s$.

    \item Because $\sigma_{n}$ belongs to the interval
    $(\sigma_{\min}(\alpha),\sigma_{\max}(\alpha))$,
    there exist two sequences in $[0,1]$,
    $\{c_{1n}\}$ and $\{c_{2n}\}$ for $n\in\mathbb{N}$,
    such that no point $(s,c)$ on the saturation nullcline
    has $c \in [c_{1n}, c_{2n}]$.
\end{enumerate}

The saturation nullcline $\{s_\xi=0\}$ for the system~\eqref{eq:dyn_sys_cap_diff}
are continuous in $\alpha$.
with respect to $\alpha_n$ as $\alpha_n \to 0^+$.
If $\sigma_{\infty}=\sigma_T$,
then for large enough $n$ the saturation nullcline is close to
the set $\{f(s,c)=\sigma^u s\}$,
which coincides with the stable and unstable manifolds of the equilibrium $A$;
see Fig.~\ref{fig:f_nonmonotone}(b).
This leads to $|c_{2n}-c_{1n}| \to 0^+$ as $\alpha_n \to 0^+$.
Arguing by contradiction,
we assume that $\sigma_{\infty}>\sigma_T$,
which means that the saturation nullcline stays close to the set
$\{f(s,c)=\sigma_\infty s\}$,
which satisfies the properties (1)-(3).
Therefore, there exists an interval $[d_1,d_2] \subseteq (0, 1)$ and $N\in\mathbb{N}$
such that $[d_1,d_2]\subseteq [c_{1n},c_{2n}]$ for all $n>N$.
Thus, there is no point $(s, c)$ on the saturation nullcline with $c \in [d_1,d_2]$.
It follows that
\[
\text{$f(s, c) - \sigma_n(s + \alpha_n a_1(1)) < 0$
for all $s\in(0,1)$, $c\in [d_1,d_2]$, and $n > N$.}
\]
Therefore, for some $\varepsilon>0$,
we have that
\begin{align*}
&(\sigma_n \alpha_n)^{-1} \,\dfrac{s_\xi}{c_\xi}
= \dfrac{ f(s, c) - \sigma_n(s + \alpha_n a_1(1))}
{\kappa^{-1}(a_1(1) c - a_1(c))} > \varepsilon > 0 \\
&\qquad \text{for all $s\in(0,1)$, $c\in [d_1,d_2]$, and $n > N$.}
\end{align*}
Consequently, there exists
$\alpha_* = 2\cdot (\varepsilon \cdot \sigma^u\cdot (d_2 - d_1))^{-1}$
such that for all $\alpha\in(0, \alpha_*)$
and any trajectory $s(c)$ joining the saddle points
$U_-(\alpha_n)$ and $U_+(\alpha_n)$,
\[
\text{$\dfrac{d}{dc}s(c) > \dfrac{1}{d_2 - d_1}$ for all $c\in [d_1,d_2]$.}
\]

\begin{figure}[ht]
\centering
\includegraphics[width=0.53\textwidth]{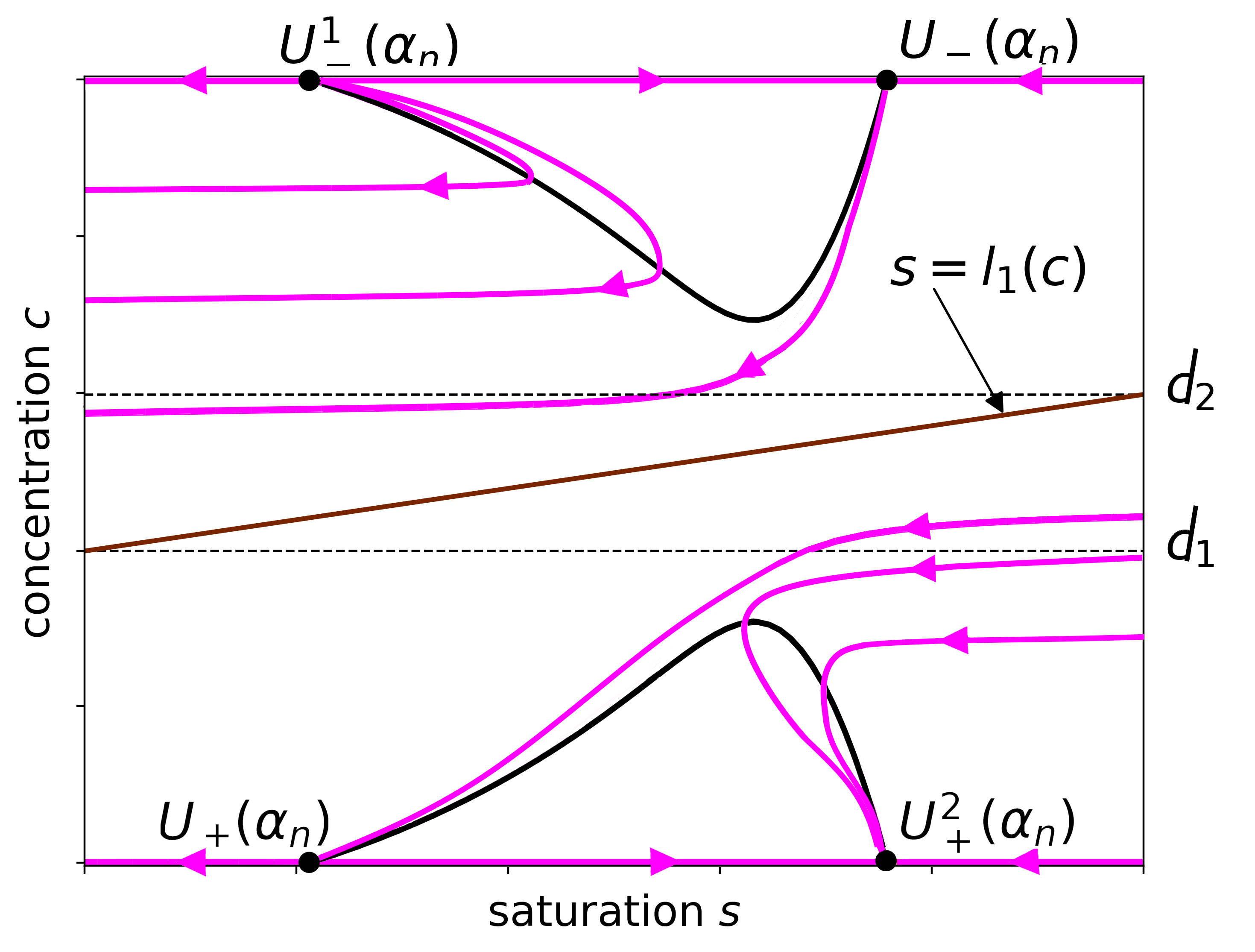}
     \caption{Illustration of the proof of Lemma~\ref{lm:sigman-limit}.}
     \label{fig:lemma43proof}
\end{figure}

Consider the linear function (a barrier)
\[
l_1(c) = \dfrac{c - d_1}{d_2 - d_1}.
\]
There is no trajectory, from either $U_-$ or $U_+$,
that intersects the graph $s = l_1(c)$.
Thus, the trajectory between saddle points does not exist.
This is a contradiction with Prop.~\ref{Theorem-nonmonotonicity}.
We conclude that $\sigma(\alpha)$ converges to $\sigma^u$ as $\alpha \to 0^+$.

Second, let us prove that $U_-(\alpha)\to U_-^u $ and $ U_+(\alpha)\to U_+^u$
as $\alpha \to 0^+$.
By the Rankine-Hugoniot condition~\eqref{eq:Rankine-Hugoniot} for the shock wave
connecting states $U_-(\alpha)$ and $U_+(\alpha)$ with speed $\sigma(\alpha)$,
we obtain that
\begin{equation}
\label{eq:RH-1}
\begin{split}
    \sigma(\alpha)&=\frac{f(s_-(\alpha),1)}{s_-(\alpha)+\alpha\,a_1(1)}
    \quad \text{ and }\quad
    \sigma(\alpha)=\frac{f(s_+(\alpha),0)}{s_+(\alpha)+\alpha\,a_1(1)}.
\end{split}
\end{equation}
Under the restriction that $U_-\in\{\lambda^s<\lambda^c\}$,
the first of these relations defines $s_-(\alpha)$ as a function of $\sigma$.
Moreover, by the implicit function theorem,
this dependence is smooth.
We know that $\sigma(\alpha)$ tends to $\sigma^u$ as $\alpha \to 0^+$;
therefore $s_-(\alpha)$ tends to some value $s_\infty$
that satisfies the equation
\begin{equation}
\label{eq:relations-sinfty}
\text{$\sigma^u=\frac{f(s_\infty,1)}{s_\infty}=\lambda^c_0(U_{\infty})$,
where $U_\infty=(s_\infty,1)\in\{\lambda^s<\lambda^c\}$.}
\end{equation}
This relation implies that the point $U_\infty$
is located on the stable manifold of the equilibrium~$A$,
and due to uniqueness of such a point in the domain
$\{\lambda^s<\lambda^c\}\cap\{c=1\}$,
the state $U_\infty$ coincides with the state $U_-^u$.
By analogous reasoning, $U_+(\alpha)\to U_+^u$ as $\alpha \to 0^+$.
Lemma~\ref{lm:sigman-limit} is proved.

\medskip
Let us make two comments on Theorem~\ref{thm:undercompressive-admissible}.

First,
although we present an explicit example given by formula~\eqref{eq:f-nonmonotone},
we expect that a similar proof works for a more general class of the Buckley-Leverett functions
that change monotonicity in $c$ only once
(\emph{e.g.},\ the models treated in Ref.~\cite{nonmonotonicity2021}).

Second, one can view the adsorption criterion for contact discontinuities
from the following perspective.
The system for the polymer model~\eqref{eq:conslaw}
is a limit of the viscous system
with adsorption~\eqref{eq:main_system_smooth_cap_diff} as $\alpha \to 0^+$,
$\varepsilon_c \to 0^+$, and $\varepsilon_d \to 0^+$. There are different ways of taking this limit, for example,
\begin{align*}
\lim\limits_{\substack{\varepsilon_c\to0^+ \\ \varepsilon_d\to0^+}}\lim\limits_{\alpha\to0^+} 
\qquad\text{or}\qquad
\lim\limits_{\alpha\to0^+} \lim\limits_{\substack{\varepsilon_c\to0^+ \\ \varepsilon_d\to0^+}}.
\end{align*}
If we first take the limit as $\alpha \to 0^+$,
we find no traveling wave solutions (with non-constant $c$),
so that no solution satisfying the vanishing viscosity admissibility criterion exists. 
Interchanging the limits and considering the limit in $\varepsilon_d,\varepsilon_c\to0^+$ first, and then $\alpha\to0^+$, 
resolves the issue, and corresponds to the introduced vanishing adsorption criterion. The proof of Theorem~\ref{thm:undercompressive-admissible}
relies on the assumption that the ratio $\kappa=\varepsilon_d/\varepsilon_c >0$ is maintained fixed, however it is valid under the milder condition
that $\alpha\cdot\varepsilon_c/\varepsilon_d \to 0^+$. Interestingly, the set of admissible contact discontinuities does not depend on the choice of $\kappa>0$. 
The full analysis of the admissible contact discontinuities for all possible limits $\alpha, \varepsilon_d, \varepsilon_c\to0^+$ is an interesting problem to solve, but is beyond the scope of this paper.


\subsection{Example Riemann problem}
\label{sec:numerics}
We now apply our results for the non-monotone ``boomerang'' model
to an interesting Riemann problem
studied by Shen~\cite{shen2017}.
(See Bressan et al.~\cite{bress2019} for further developments
founded on casting system~(\ref{eq:conslaw}) in Lagrangian form
using the transformation of Pires et al.~\cite{pires2006},
as justified rigorously by Wagner~\cite{wagner1987}.)

Shen studied a system of equations closely related to
Eq.~\eqref{eq:main_system_smooth_cap_diff} without adsorption.
In particular, in her Examples~4.5 and~4.6,
she considered Riemann initial data with (i)~$s_L = s_R$
and (ii)~$c_L \ne c_R$ having the property that $\mu(c_L) = \mu(c_R)$.
Because the $c$-wave speeds $f(s_L, c_L)/s_L$ and $f(s_R, c_R)/s_R$ are equal,
the jump between $(s_L, c_L)$ and $(s_R, c_R)$ is a contact discontinuity.
Therefore, one solution of the Riemann problem
for the non-diffusive system
consists of a single moving contact discontinuity.
However,
in numerical simulations of this Riemann problem
for the system with diffusion,
Shen observed that this single-wave solution
decayed into a three-wave solution
that depended sensitively on the diffusion coefficients.

Here we investigate a particular Riemann problem of this type for system~\eqref{eq:main_system_smooth_cap_diff},
again with no adsorption ($\alpha = 0$);
the fractional flow and viscosity ratio functions
are given by Eqs.~\eqref{eq:f-nonmonotone},
but with $\mu(c) = 1 + 4\,c\,(1 - c)$.
The Riemann data has $s_L = s_R = 0.8$,
$c_L = 0.9$, and $c_R = 0.1$, so that $\mu(c_L) = \mu(c_R)$.
We derive an analytical solution for the non-diffusive system
($\varepsilon_c=\varepsilon_d=0$),
and we perform numerical simulations of the diffusive system
with $\varepsilon_c=1$ and $\varepsilon_d=0.5$
using the second-order nonlinear Crank-Nicolson scheme
described in Ref.~\cite{lam2020-RCD}.

For this example, we show that:
\begin{itemize}
    \item there is a second solution of the non-diffusive system,
    which is admissible according to the vanishing-adsorption criterion;
    it comprises a slow $s$-shock wave,
    an undercompressive contact discontinuity,
    and a fast $s$-shock wave;
    
    \item in a numerical simulation of the diffusive system,
    the Riemann data resolves into three waves rather than a single wave; and
            
    \item the numerical simulation agrees
    quantitatively with non-diffusive three-wave solution.
\end{itemize}

\begin{figure}
    \centering
    \includegraphics[width=0.45\textwidth]{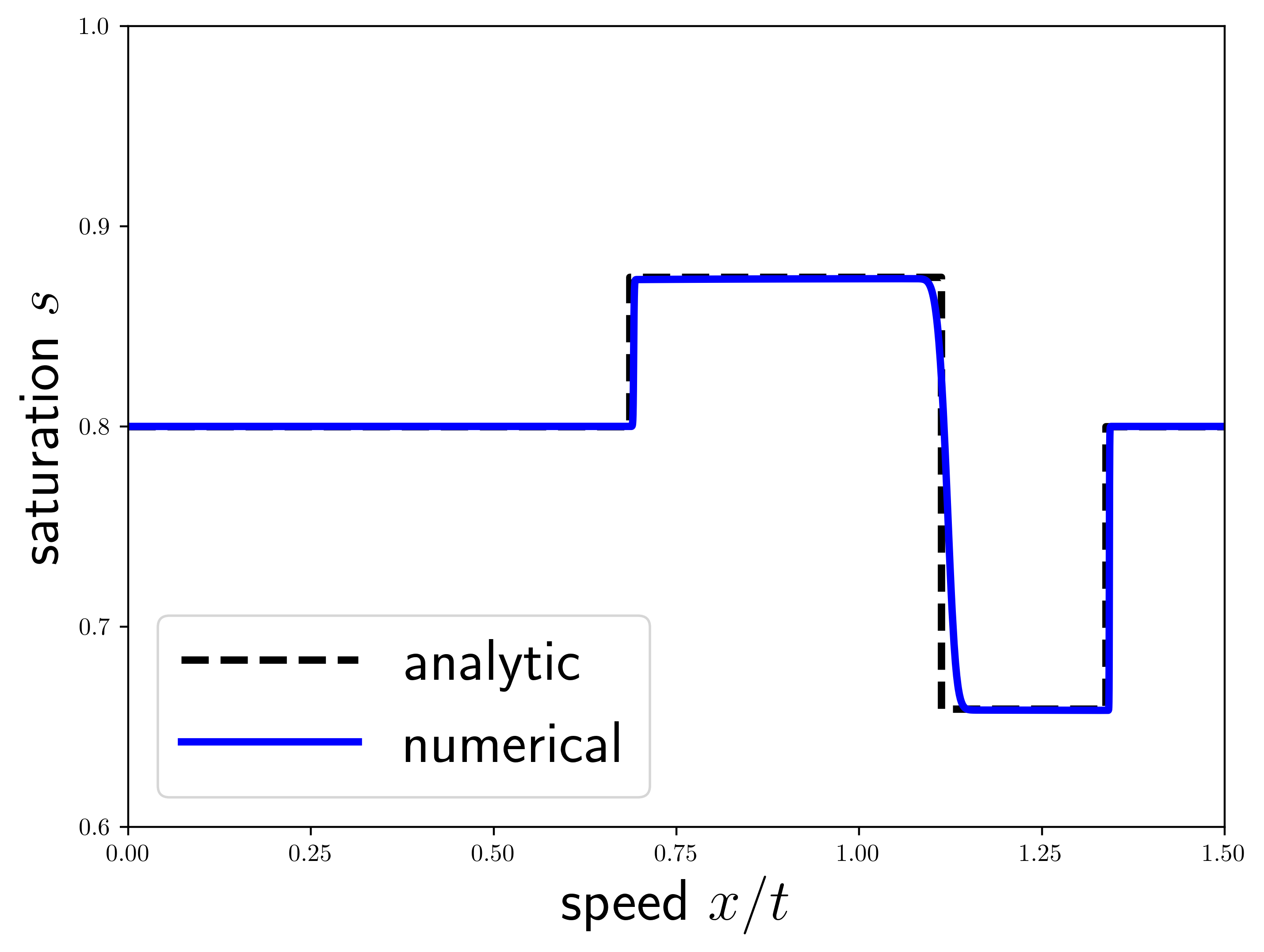}
    \includegraphics[width=0.45\textwidth]{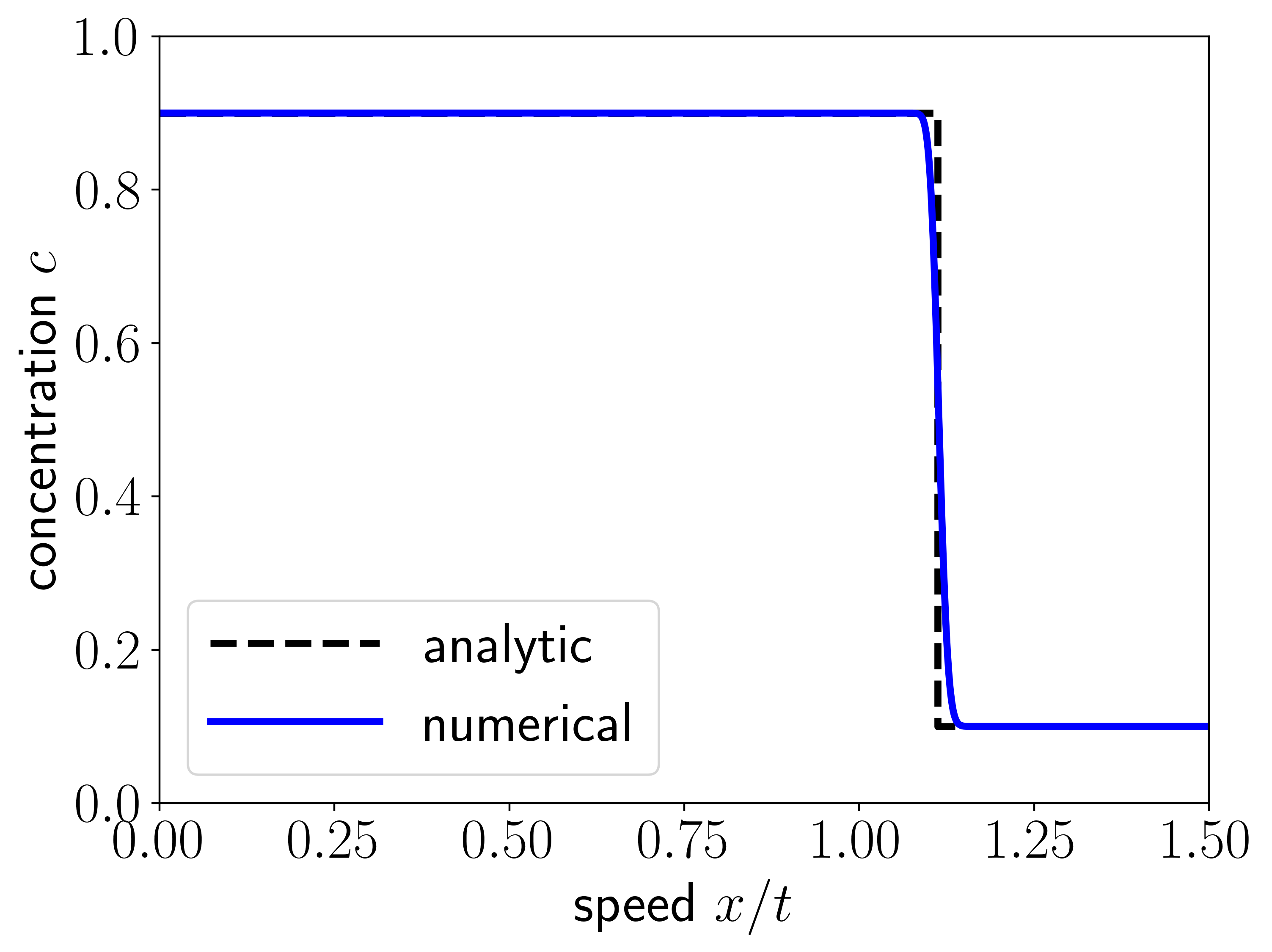}
    \caption{Comparison of the analytic and numerical solutions
    for the Riemann initial-value problem.}
    \label{fig:numerics}
\end{figure}

The Riemann solution~\eqref{eq:analyt-sol} of the non-diffusive system,
which we describe below in more detail,
is depicted as a dashed curve in Fig.~\ref{fig:numerics},
which also shows the result of a numerical simulation
of the diffusive system as a solid curve.
This figure demonstrates excellent agreement between
the numerical simulation and the analytic solution.

The three-wave solution to the Riemann problem
for the non-diffusive system contains the following sequence of waves:
\begin{align}
\label{eq:analyt-sol}
    U_L\xrightarrow{s} U_{M_1}\xrightarrow{c} U_{M_2}\xrightarrow{s} U_R,
\end{align}
where the wave $U_L\xrightarrow{s} U_{M_1}$ is a slow $s$-shock wave, $U_{M_1}\xrightarrow{c} U_{M_2}$ is an undercompressive contact discontinuity, $U_{M_2}\xrightarrow{s} U_{R}$ is a fast $s$-shock wave.
Here $U_{M_1}=(s_{M_1}, c_{L})$ and $s_{M_1}$ is the largest solution $s$
of the equation $\sigma_{\min}(0) = f(s, c_L)/s$,
where $\sigma_{\min}(0)$ is defined in Rem.~\ref{rmk:sigmamin-max}.
Analogously, $U_{M_2}=(s_{M_2}, c_{R})$ and $s_{M_2}$
is the smallest solution $s$ of the equation $\sigma_{\min}(0) = f(s, c_R)/s$.
For our choice of the fractional flow function~\eqref{eq:f-nonmonotone},
the equations for $s_{M_1}$ and $s_{M_2}$ are quadratic.
Moreover, $\sigma_{\min}(0)$ is defined by the requirement that
the quadratic equation $\sigma_{\min}(0) = f(s, c_\text{max})/s$,
where $c_\text{max} = 0.5$ maximizes $\mu(c)$,
has a double root,
i.e., its discriminant vanishes;
this requirement is a quadratic equation for $\sigma_{\min}(0)$.
Thus, we obtain the Riemann solution analytically.

\section{Acknowledgments}
This work was facilitated using the program ELI,
which is an interactive, graphical tool for exploring Riemann problems
(\href{ https://eli.fluid.impa.br/}{https://eli.fluid.impa.br/}). 
YP was supported by Bolsa de Excelência at IMPA, CNPq: 406460/2023-0 and Coordenacao de Aperfeicoamento de Pessoal de Nivel
Superior - Brasil (CAPES) - 23038.015548/2016-06. YP is grateful to Aleksandr Enin for fruitful discussions
as well as help with Python scripts for creating figures. YP also thanks all colleagues from the Fluid Dynamics Laboratory at IMPA, Centro~PI at IMPA and DMAT at PUC-Rio for creating an excellent scientific atmosphere.
DM acknowledges support from 
FAPERJ: E-26/202.764/2017,
FAPERJ: E-26/201.159/2021,
CNPq: 306566/2019-2,
CNPq: 405366/2021-3,
and FAPERJ-PRONEX E-26/010.001267/2016. 
DM is also grateful to Hans Bruining from the Technical University of Delft
for sharing his expertise on flow in porous media.


\newpage

\end{document}